\documentclass[preprint, 12pt]{elsarticle}
  \usepackage{amsmath}
  \usepackage{amssymb}
  \usepackage{bbm}
  \usepackage{graphicx}
  \usepackage{subfig}
\newdefinition{remark}{Remark}
\newcommand{\du}{\, \mathrm{d}}
\newcommand{\myremarkend}{~\hfill$\clubsuit\/$}

\journal{Journal of Computational Physics}
\begin{document}

\begin{frontmatter}

\title{An efficient method for the incompressible Navier-Stokes equations on
       irregular domains with no-slip boundary conditions, high order up to
       the boundary.}

\author{D. Shirokoff}
\ead{shirokof@math.mit.edu}
\author{R. R. Rosales}
\ead{rrr@math.mit.edu}
\address{Department of Mathematics, Massachusetts Institute of Technology,
         Cambridge, MA, USA, 02139}

\begin{abstract}
Common efficient schemes for the incompressible Navier-Stokes equations, such
as projection or fractional step methods, have limited temporal accuracy as a
result of matrix splitting errors, or introduce errors near the domain
boundaries (which destroy uniform convergence to the solution).  In this paper we
recast the incompressible (constant density) Navier-Stokes equations (with the
velocity prescribed at the boundary) as an equivalent system, for the primary
variables velocity and pressure. We do this in the usual way away from the
boundaries, by replacing the incompressibility condition on the velocity by
a Poisson equation for the pressure. The key difference from the usual
approaches occurs at the boundaries, where we use boundary conditions that
unequivocally allow the pressure to be recovered from knowledge of the
velocity at any fixed time. This avoids the common difficulty of an,
apparently, over-determined Poisson problem.  Since in this alternative
formulation the pressure can be accurately and efficiently recovered from the
velocity, the recast equations are ideal for numerical marching methods. The
new system can be discretized using a variety of methods, including semi-implicit treatments of viscosity, and in principle to any
desired order of accuracy. In this work we illustrate the approach with a 2-D
second order finite difference scheme on a Cartesian grid, and devise an
algorithm to solve the equations on domains with curved (non-conforming)
boundaries, including a case with a non-trivial topology (a circular
obstruction inside the domain). This algorithm achieves second order accuracy
in the $L^{\infty}\/$ norm, for both the velocity and the pressure. The scheme
has a natural extension to 3-D.
\end{abstract}


\end{frontmatter}

\section{Introduction.} \label{sec:intro}
A critical issue in the numerical solution of the incompressible Navier-Stokes
equations is the question of how to implement the incompressibility constraint.
Equivalently, how to recover the pressure from the flow velocity, given the
fact that the equations do not provide any boundary condition for the pressure. This 
has been an area of intense research, ever since
the pioneering MAC scheme \cite{HarlowWelch1965} of Harlow and Welch in 1965.
Of course, one can avoid the problem by simultaneously discretizing the
momentum and the divergence free equations, as in the difference scheme
proposed by Krzywicki and Ladyzhenskaya~\cite{KrzywickiLadyzhenskaya1966}, which can 
be shown to converge --- while avoiding the need for any pressure
boundary conditions. Approaches such as these, however, do not lead to efficient schemes.

Generally the dilemma has been that of a trade-off between efficiency, and
accuracy of the computed solution near the boundary. However, many
applications require both efficiency, and accuracy.  For example, 
to calculate fluid solid interactions, both the pressure and gradients of the
velocity are needed at the solid walls, as they appear in the components of
 the stress tensor. Furthermore, these objectives  
must be achievable for ``arbitrary'' geometries, not just
simple ones with symmetries that can be exploited. Unfortunately, these
requirements are not something that current algorithms are generally well
suited for, as the brief review below is intended to show.\footnote{This
   is not intended as a thorough review of the field, and we apologize for
   the many omissions.} However, we believe that algorithms based on a Presure Poisson Equation
(PPE) reformulation of the Navier-Stokes equations --- reviewed towards
the end of this introduction --- offer a path out of the dilemma. The
work presented in this paper is, we hope, a contribution in this direction.

\emph{Projection methods} are very popular in practice because they are
efficient. They achieve this efficiency by: (i) Interpreting the pressure
as effecting a projection of the flow velocity evolution into the set of
incompressible fields. That is, write the equations in the form
$ \mathbf{u}_t = \mathcal{P} \left(\mu\,\Delta\,\mathbf{u} -
  (\mathbf{u} \cdot \nabla )\,\mathbf{u} + \mathbf{f} \right)\/$,
where $\cal{P}\/$ is the appropriate projection operator, $\mu\/$ is the
kinematic viscosity, $\mathbf{u}\/$ is the flow velocity vector, and
$\mathbf{f}\/$ is the vector of applied body forces. (ii) Directly
evolving the flow velocity. The question is then how to compute
$\mathcal{P}\/$.

In their original formulation by Chorin~\cite{Chorin} and Temam~\cite{Temam},
the projection method was formulated as a time splitting scheme in which:
First an intermediate velocity is computed, ignoring incompressibility.
Second, this velocity is projected onto the space of incompressible vector
fields --- by solving a Poisson equation for pressure.
Unfortunately this process introduces numerical boundary layers into the solution,
which can be ameliorated (but not completely suppressed) for simple
geometries --- e.g. ones for which a staggered grid approach can be
implemented \cite{E_Liu}.

The development of \emph{second order projection methods}
\cite{BellCollelaGlaz1989,KarniadakisIsraeliOrszag1991,KimMoin1985,%
OrszagIsraeliDeville1986,Kan1986} provided greater control over the numerical
boundary layers and accuracy in the pressure \cite{brownCortezMinion:2001}. 
These are the most popular schemes used in practice. However,
particularly for moderate or low Reynolds numbers, the effects of the
numerical boundary layers can still be problematic~\cite{Guermond}.
Non-conforming boundaries add an extra layer of difficulty. The search for
means to better control these numerical artifacts is an active area of
research.

The numerical boundary layers in projection methods reflect in the known
convergence results for them (e.g. \cite{Prohl,Rannacher,Shen}). Convergence
is stated only in terms of integral norms, with the main difficulties near
the boundary. There point-wise convergence (and even less convergence of the
flow velocity gradient) cannot be guaranteed --- even if the solution is known
to be smooth. Hence the accurate calculation of wall stresses with these
methods is problematic. Guermond, Minev and Shen~\cite{Guermond} provide
further details on convergence results, as well as an extensive review of
projection methods and the improved pressure-correction schemes.

Two other methods for solving the Navier-Stokes equations are the
\emph{immersed
boundary}~\cite{Linnick,MittalIaccarino2005,Peskin,Peskin2002,Taira},
and the
\emph{vortex-streamfunction}~\cite{BenArtzi, Calhoun,Napolitano} methods.
These also decouple the calculation of the velocity and of the pressure. The
immersed boundary method does so by introducing Dirac forces to replace the
domain walls, which makes obtaining high order implementation of the boundary
conditions difficult. The vortex-streamfunction formulation decouples the
equations, but has dimensional limitations.  An interesting variation of the
vortex-streamfunction approach, using only local boundary conditions, is
presented in reference~\cite{HammoutiAbdelkader2009}.

Closely related to the immersed boundary methods are the \emph{penalty}
(alternatively: \emph{fictitious domain} or \emph{domain embedding}) methods
--- e.g.
see~\cite{AngotBruneauFabrie1999,GlowinskiPanKearsleyPeriaux1995,KhadraAngotParneixCaltagirone2000}.
These methods, effectively, replace solid walls in the fluid by a porous media
with a small porosity $0 < \eta \ll 1\/$. In the limit $\eta \to 0\/$, this
yields no slip and no flow-through at the solid walls. Two important advantages
of this approach are that complicated domains are easy to implement, and that
the total fluid-solid force can be computed using a volume integral, rather
than an integral over the boundary of the solid. Unfortunately, the parameter
$\eta\/$ introduces $\sqrt{\eta}\/$ boundary layers which make convergence
slow and high accuracy computations expensive, since $\eta\/$ cannot be
selected independently of the numerical grid size.

Finally, we mention the algorithms based on a \emph{Pressure Poisson Equation
(PPE) reformulation of the Navier-Stokes
equations}~\cite{Henshaw1994, HenshawKreissReyna1994, HenshawPetersson2003, 
LiuLiuPego, GreshoSani1987,JohnstonLiu2002,JohnstonLiu2004,%
KleiserSchumann1980,Quartapelle,Rempfer,SaniShenPironneauGresho2006},
which
is the class of methods within which the work presented in this paper falls.
In this approach the incompressibility constraint for the flow velocity is
replaced by a Poisson\footnote{The choice of the Poisson equation for the
    pressure is not unique, e.g. see \cite{JohnstonLiu2002}.}
equation for the pressure. This then allows an extra boundary condition ---
which must be selected so that, in fact, incompressibility is maintained
by the resulting system. This strategy was first proposed by
Gresho and Sani~\cite{GreshoSani1987}, who pointed out that adding
$\nabla \cdot \mathbf{u} = 0\/$ as a boundary condition yields a system of
equations that is equivalent to the Navier-Stokes equations. Unfortunately,
their particular PPE formulation incorporates no explicit boundary condition
that can be used to recover the pressure from the velocity, by solving a
Poisson problem --- for a more detailed discussion of this, see
remark~\ref{rem:equivalence1} in this paper. In~\cite{Henshaw1994,  HenshawKreissReyna1994} the issue is resolved at the discrete numerical level, where they demonstrate high order schemes. For instance \cite{Henshaw1994} demonstrates a fourth-order in space and 2nd-order in time implementation using overlapping grids.  Subsequent work at the continuum level was later introduced by Henshaw et. al. \cite{HenshawPetersson2003} and Johnston and Liu~\cite{JohnstonLiu2004}.  Recently, work in PPE formulations have led to interesting improvements and analysis of projection methods \cite{LiuLiuPego}. In this paper --- in
equations (\ref{PPE_Numerical_Splitting1}--\ref{PPE_Numerical_Splitting2})
--- we present another PPE system, also equivalent to the Navier-Stokes
equations, which allows an explicit recovery of the pressure given the
flow velocity.  A comparison with the one in~\cite{JohnstonLiu2004} can be found in remark~\ref{rem:JohnstonLiu2004} in this paper.  Subtle issues can arise with semi-implicit approaches (pressure treated explicitly and viscous term implicitly) leading to time step restrictions \cite{Petersson2001}. In section \ref{Section:SemiImplicitStability} we show that semi-implicit implementations of our scheme do not have time step restrictions of the diffusive type.

PPE reformulations of the Navier-Stokes equations, such as the one in
equations (\ref{PPE_Numerical_Splitting1}--\ref{PPE_Numerical_Splitting2}),
or in reference~\cite{JohnstonLiu2004}, have important advantages over the
standard form of the equations. First, the pressure is not implicitly coupled
to the velocity through the momentum equation and incompressibility. Hence
it can be directly (and efficiently) recovered from the velocity field by
solving a Poisson equation. This allows the velocity field to be marched in
time using the momentum equation, with the pressure interpreted as some
(complicated) function of the velocity. Second, no spurious boundary layers
are generated for neither the velocity, nor the pressure. This because:
\begin{itemize}
 \item[---] \vspace*{-0.7ex}
 There are no ambiguities as to which boundary conditions to use for
    the pressure --- hence errors induced by not-quite-correct boundary
    conditions do not occur. Note that these errors should not be confused
    with the truncation errors that any discretization of the equations will
    produce. Truncation errors are controlled by the order of the
    approximation and, for smooth solutions, are uniformly small.
 \item[---] \vspace*{-0.7ex}
 Incompressibility is enforced at all times.
\end{itemize}
Hence pressure and velocities that are accurate everywhere can be obtained,
in particular: near the boundaries. Finally, PPE formulations allow, at least
in principle, for the systematic generation of higher order approximations.

It follows that PPE strategies offer the promise of a resolution of the
dilemma alluded to in the second paragraph of this introduction. They retain
many of the advantages that have made projection methods popular,
while not suffering from the presence of numerical boundary layers.  On
the negative side, the boundary conditions for PPE systems tend to be more
complicated than the simple ones that the standard form of the equations has.

This paper is organized as follows.
In \S~\ref{Section:PPE} we discuss the Pressure Poisson Equation (PPE)
formulation of the Navier-Stokes equations.
In \S~\ref{Section:TheorySplitting} we present a PPE formulation,
using an alternative form of the boundary conditions, that allows a complete
splitting of the momentum and pressure equations. Namely: the pressure can
be recovered from the flow velocity without boundary condition ambiguities.
In \S~\ref{Section:NumSplitting} we address a numerical question (secular
growth of the error under some conditions) and introduce a corrected
formulation, involving a feedback parameter $\lambda\/$, where this potential
problem is corrected. This system, which is equivalent to the Navier-Stokes
equations, is displayed in equations
(\ref{PPE_Numerical_Splitting1}--\ref{PPE_Numerical_Splitting2}).
A discussion and comparison with the alternative PPE formulation by
Johnston and Liu~\cite{JohnstonLiu2004} is also included in this section.
In \S~\ref{Section:SemiImplicitStability} we introduce and examine a second order in time, semi-implicit in viscosity scheme.  We show the scheme is unconditionally stable with no diffusive time step restriction.
In \S~\ref{Section:NumScheme} we describe a second order solver for the new
system introduced in \S~\ref{Section:NumSplitting}, for irregular domains
embedded within a cartesian staggered grid.
In \S~\ref{Section:Implementation}, we implement and test our proposed
schemes. This section includes convergence plots, indicating \emph{second
order uniform convergence} all the way to the boundary ($L^\infty\/$ norm),
\emph{of the pressure, the flow velocity, and of the derivatives of the
flow velocity.}
Finally, \S~\ref{Section:conclusions} has the conclusions.
The contents of the two appendices is as follows.
In \ref{solvabilityBreakdown} we present an extended system of equations,
valid for arbitrary flows, which for smooth solutions has the Navier-Stokes
equations as an attractor.
In \ref{divfreeBC}, for completeness, we display formulas for the
$\nabla \cdot \mathbf{u} = 0\/$ boundary condition for general conforming
boundaries in curvilinear coordinates.


\section{The Pressure Poisson Equation.} \label{Section:PPE}
In this section we introduce the well-known pressure Poisson equation (PPE), and
use it to construct a system of equations (and boundary conditions) equivalent
to the constant-density (hence incompressible) Navier-Stokes equations, with the
velocity prescribed at the boundaries.  Specifically, consider the
incompressible Navier-Stokes equations in a connected domain
$\Omega \in \mathbbm{R}^D\/$, where $D=2\/$ or $D =3\/$, with a piece-wise
smooth boundary $\partial\,\Omega\/$. Inside $\Omega\/$, the flow velocity field
$\mathbf{u}(\mathbf{x}\/,\,t)\/$ satisfies the equations
\begin{eqnarray} 
 \mathbf{u}_t + (\mathbf{u} \cdot \nabla)\,\mathbf{u} & = &
    \mu\,\Delta\/\mathbf{u} - \nabla\/p + \mathbf{f}\/,
 \label{NS_Momentum} \\
 \nabla \cdot \mathbf{u}                              & = & 0\/,
 \label{NS_Mass}
\end{eqnarray}
where $\mu\/$ is the kinematic viscosity,\footnote{We work in non-dimensional variables, so
   that $\mu = 1/Re\/$ (where $Re\/$ is the Reynolds number) and the fluid
   density is $\rho = 1\/$.}
$p(\mathbf{x}\/,\,t)\/$ is the pressure, $\mathbf{f}(\mathbf{x}\/,\,t)$ are the
body forces, $\nabla\/$ is the gradient, and $\Delta = \nabla^2\/$ is the
Laplacian. Equation (\ref{NS_Momentum}) follows from the conservation of
momentum, while (\ref{NS_Mass}) is the incompressibility condition (conservation
of mass).

In addition, the following boundary conditions apply
\begin{equation} \label{NS_Boundary_Conditions}
   \mathbf{u} = \mathbf{g}(\mathbf{x}\/,\,t) \qquad
   \textrm{for} \hspace*{2mm} \mathbf{x} \in \partial \Omega\/,
\end{equation}
where
\begin{equation} \label{NS_Boundary_ConditionsConstr}
   \int_{\partial\/\Omega} \mathbf{n} \cdot \mathbf{g} \du A = 0\/,
\end{equation}
$\mathbf{n}\/$ is the outward unit normal on the boundary, and $\du A\/$ is the
area (length in 2D) element on $\partial\/\Omega\/$.
Equation~(\ref{NS_Boundary_ConditionsConstr}) is the consistency condition for
$\mathbf{g}\/$, since an incompressible fluid must have zero net flux through
the boundary.

Finally, we assume that initial conditions are given
\begin{eqnarray}
   \mathbf{u}(\mathbf{x}\/,\,0) & = &
      \parbox{8ex}{$\mathbf{u}_0$} \quad
      \textrm{for} \hspace*{2mm} \mathbf{x} \in \Omega\/,
      \label{NS_Initial_Conditions} \\
   \nabla \cdot \mathbf{u}_0    & = &
      \parbox{8ex}{$0$}            \quad
      \textrm{for} \hspace*{2mm} \mathbf{x} \in \Omega\/,
      \label{InitVel_Consistency1}  \\
   \mathbf{u}_0(\mathbf{x})     & = &
      \parbox{8ex}{$\mathbf{g}(\mathbf{x}\/,\,0)$} \quad
      \textrm{for} \hspace*{2mm} \mathbf{x} \in \partial \Omega\/.
      \label{InitVel_Consistency2}
\end{eqnarray}
\begin{remark} \label{rem:gzero}
  Of particular interest is the case of fixed impermeable walls, where no flux
  $\mathbf{u} \cdot \mathbf{n} = 0\/$ and no slip
  $\mathbf{u} \times \mathbf{n} = 0\/$ apply at $\partial \Omega\/$. This
  corresponds to $\mathbf{g} = 0\/$ in (\ref{NS_Boundary_Conditions}). Note
  that the no-slip condition is equivalent to $\mathbf{u}\cdot\mathbf{t}=0\/$
  for all unit tangent vectors $\mathbf{t}\/$ to the boundary.\myremarkend
\end{remark}
\begin{remark} \label{rem:movingboundary}
  In this paper we will assume that the domain $\Omega\/$ is fixed. Situations
  where the boundary of the domain, $\partial \Omega\/$, can move --- either by
  externally prescribed factors, or from interactions with the fluid, are of
  great physical interest. However, to keep the presentation of the new method
  as simple as possible, we postpone consideration of these cases for further
  work. \myremarkend
\end{remark}
Next we introduce the Pressure Poisson Equation (PPE). To obtain the pressure
equation, take the divergence of the momentum equation (\ref{NS_Momentum}),
and apply equation (\ref{NS_Mass}) to eliminate the viscous term and the term
with a time derivative. This yields the following Poisson equation for the
pressure
\begin{equation} \label{eqn:PEP}
   \Delta\/p = \nabla \cdot \left( \mathbf{f} - (\mathbf{u} \cdot \nabla)\,
              \mathbf{u} \right)\/.
\end{equation}
Two crucial questions are now (for simplicity, \emph{assume solutions that are
smooth all the way up to the boundary})
\begin{itemize}
 \item[\ref{Section:PPE}a]
 \emph{Can this equation be used to replace the incompressibility condition
 (\ref{NS_Mass})?} Since --- given (\ref{eqn:PEP}) --- the divergence of
 (\ref{NS_Momentum}) yields the heat equation
 \begin{equation} \label{eqn:heatdivu}
    \phi_t = \mu \Delta\/\phi\/,
 \end{equation}
 for $\,\phi = \nabla \cdot \mathbf{u}\/$ and $\mathbf{x} \in \Omega\/$, it
 would seem that the answer to this question is yes --- provided that the
 initial conditions are incompressible (\emph{i.e.:}~$\phi = 0\/$ for
 $t = 0\/$). However, this works only if we can guarantee that, at all times,
 \begin{equation} \label{eqn:bcdivu}
    \phi = 0\/,
 \end{equation}
 for $\mathbf{x} \in \partial \Omega\/$.
 \item[\ref{Section:PPE}b]
 \emph{Given the flow velocity $\mathbf{u}\/$, can (\ref{eqn:PEP}) be used to
 obtain the pressure $p\/$?} Again, at first sight, the answer to this question
 appears to be yes. After all, (\ref{eqn:PEP}) is a Poisson equation for $p\/$,
 which should determine it uniquely --- given appropriate boundary conditions.
 The problem is: what boundary conditions? Evaluation of (\ref{NS_Momentum}) at
 the boundary, with use of (\ref{NS_Boundary_Conditions}), shows that the flow
 velocity determines the whole gradient of the pressure at the boundary, which
 is too much for (\ref{eqn:PEP}). Further, if only a portion of these boundary
 conditions are enforced when solving (\ref{eqn:PEP}) --- say, the normal
 component of (\ref{NS_Momentum}) at the boundary, then how can one be sure
 that the whole of (\ref{NS_Momentum}) applies at the boundary?
\end{itemize}
\begin{remark} \label{rem:getpfromu}
  From an algorithmic point of view, an affirmative answer to the questions
  above would very useful, for then one could think of the pressure as some
  (global) function of the flow velocity, in which case (\ref{NS_Momentum})
  becomes an evolution equation for $\mathbf{u}\/$, which could then be solved
  with a numerical ``marching'' method.\myremarkend
\end{remark}
The issue in item~\ref{Section:PPE}a can be resolved easily, and we do so
next.  We postpone dealing with the issue in item~\ref{Section:PPE}b till
the next section, \S~\ref{Section:TheorySplitting}. Since the addition of an
equation for the pressure allows the introduction of one extra boundary
condition, we propose to replace the system in (\ref{NS_Momentum}),
(\ref{NS_Mass}), and (\ref{NS_Boundary_Conditions}) by the following Pressure
Poisson Equation (PPE) formulation:
\begin{eqnarray}
   \mathbf{u}_t + (\mathbf{u} \cdot \nabla)\,\mathbf{u} & = &
      \mu\,\Delta\/\mathbf{u} - \nabla\/p + \mathbf{f}\/,
                   \label{PPE_Bulk_1} \\
   \Delta\/p                                          & = &
      \nabla \cdot \left( \mathbf{f} - (\mathbf{u} \cdot \nabla)\,\mathbf{u}
      \right)\/,   \label{PPE_Bulk_2}
\end{eqnarray}
for $\mathbf{x} \in \Omega\/$, with the boundary conditions
\begin{eqnarray}
   \mathbf{u}              & = & \mathbf{g}(\mathbf{x}\/,\,t)\/,
       \label{PPE_Boundary_1} \\
   \nabla \cdot \mathbf{u} & = & 0\/, \label{PPE_Boundary_2}
\end{eqnarray}
for $\mathbf{x} \in \partial \Omega\/$ --- where, of course, the restriction
in (\ref{NS_Boundary_ConditionsConstr}) still applies. The extra boundary
condition is precisely what is needed to ensure that the pressure enforces
incompressibility throughout the flow (see item~\ref{Section:PPE}a)
\begin{remark} \label{rem:equivalence1}
  It can be seen that for smooth enough solutions $(\mathbf{u}, p)\/$, the pair
  of equations (\ref{PPE_Bulk_1}--\ref{PPE_Bulk_2}), accompanied by the
  boundary conditions (\ref{PPE_Boundary_1}--\ref{PPE_Boundary_2}), are
  equivalent to the incompressible Navier-Stokes equations in
  (\ref{NS_Momentum}), (\ref{NS_Mass}), and (\ref{NS_Boundary_Conditions}).
  Of course, this result is not new. This reformulation of the Navier-Stokes
  equations was first presented by Gresho and Sani~\cite{GreshoSani1987} ---
  it can also be found in~reference~\cite{Quartapelle}. However, it should be
  pointed out that Harlow and Welch~\cite{HarlowWelch1965}, in their
  pioneering work, had already noticed that the boundary condition
  $\nabla \cdot \mathbf{u} = 0\/$ was needed to guarantee, within the context
  of their MAC scheme, that $\nabla \cdot \mathbf{u} = 0\/$ everywhere.

  This formulation does not provide any boundary conditions for
  the pressure, which means that one ends up with a ``global'' constraint on
  the solutions to the Poisson equation~(\ref{PPE_Bulk_2}). Hence it does not
  yield a satisfactory answer to the issue in item~\ref{Section:PPE}b, since
  recovering the pressure from the flow velocity is a hard problem with this
  approach.

  Direct implementations of (\ref{PPE_Bulk_1}--\ref{PPE_Boundary_2}) have
  only been proposed for simple geometries --- in particular: grid-conforming
  boundaries.
  In~\cite{KleiserSchumann1980} a spectral algorithm for plane channel flows
  is presented.
  We have already mentioned~\cite{HarlowWelch1965}, where they use
  a staggered grid on a rectangular domain, and the condition
  $\nabla \cdot \mathbf{u} = 0\/$ is used (at the discrete level) to close
  the linear system for the pressure.
  In~\cite{JohnstonLiu2002}, by manipulating the discretization of the boundary
  conditions and of the momentum equation~(\ref{PPE_Bulk_1}), they manage to
  obtain ``local'' approximate Neumann conditions for the pressure --- both on
  rectangular, as well as circular, domains where $\mathbf{u} = 0\/$ on the
  boundary.
  Unfortunately, the approaches in~\cite{HarlowWelch1965,JohnstonLiu2002}
  seem to be very tied up to the details of particular
  discretizations, and require a conforming boundary.  \myremarkend
\end{remark}


\section{Theoretical Reformulation} \label{Section:TheorySplitting}
We still need to deal with the issue raised in item~\ref{Section:PPE}b. In
particular, in order to implement the ideas in remark~\ref{rem:getpfromu}, we
need to split the boundary conditions for (\ref{PPE_Bulk_1}--\ref{PPE_Bulk_2}),
in such a way that: (i) there is a specific part of the boundary conditions
that is used with (\ref{PPE_Bulk_1}) to advance the velocity field in time,
given the pressure.  (ii) The remainder of the boundary conditions is used with
(\ref{PPE_Bulk_2}) to solve for the pressure --- at each fixed time, given the
flow velocity field $\mathbf{u}\/$. This is the objective of this section.

The conventional approach in projection or fractional step methods is to
associate (\ref{PPE_Boundary_1}) with equation (\ref{PPE_Bulk_1}).  This is
reasonable for evolving the heat-like equation (\ref{PPE_Bulk_1}).
Unfortunately, it has the drawback of only implicitly defining the boundary
conditions for the Poisson equation (\ref{PPE_Bulk_2}). Specifically, the
correct pressure boundary conditions are those that guarantee
$\nabla \cdot \mathbf{u} = 0\/$ for $\mathbf{x} \in \partial\/\Omega\/$, and
such boundary conditions cannot easily be known a priori as a function of
$\mathbf{u}\/$. Hence one is left with a situation where the appropriate
boundary conditions for the pressure are not known. This leads to errors in
the pressure, and in the incompressibility condition, which are difficult to
control. In particular, errors are often most pronounced near the boundary
where no local error estimates can be produced, even for smooth solutions.
By the latter we mean that the resulting schemes cannot be shown to be
consistent, all the way up to the boundary, in the classical sense of finite
differences introduced by Lax \cite{LaxFD}.

In this paper we take a different approach, which has a similar spirit to the
one used by Johnston and Liu~\cite{JohnstonLiu2004}
--- see remark~\ref{rem:JohnstonLiu2004}.
Rather than associate all the
$D\/$ components of (\ref{PPE_Boundary_1}) with equation (\ref{PPE_Bulk_1}),
we enforce the $D-1\/$ tangential components only, and complete the set of
boundary conditions for (\ref{PPE_Bulk_1}) with (\ref{PPE_Boundary_2}). Hence,
when evolving equation (\ref{PPE_Bulk_1}) we do not specify the normal
velocity on the boundary, but --- through the divergence condition
(\ref{PPE_Boundary_2}) --- specify the normal derivative of the normal
velocity. Finally, the (as yet unused) boundary condition on the normal
velocity, $\mathbf{n}\cdot (\mathbf{u} - \mathbf{g}) = 0\,$ for
$\,\mathbf{x} \in \partial\/\Omega\/$, is employed to obtain an explicit
boundary condition for equation (\ref{PPE_Bulk_2}). We do this by requiring
that the pressure boundary condition be equivalent to
$\left(\mathbf{n}\cdot (\mathbf{u} - \mathbf{g})\right)_t = 0\,$ for
$\,\mathbf{x} \in \partial\/\Omega\/$ --- which then guarantees that the
normal component of (\ref{PPE_Boundary_1}) holds, as long as the initial
conditions satisfy it. This objective is easily achieved: dotting equation
(\ref{PPE_Bulk_1}) through with $\mathbf{n}\/$, and evaluating at the
boundary yields the desired condition. The equations, with their appropriate
boundary conditions, are thus:
\begin{equation} \label{PPE_Theory_Splitting1}
   \left.
   \begin{array}{rcll}
     \mathbf{u}_t - \mu\,\Delta\/\mathbf{u}      & = &
        - \nabla\/p - (\mathbf{u} \cdot \nabla)\,\mathbf{u} + \mathbf{f}
        \quad & \textrm{for}\;\; \mathbf{x} \in \phantom{\partial}\/\Omega\/,
        \\ \rule{0ex}{2.5ex}
     \mathbf{n} \times (\mathbf{u} - \mathbf{g}) & = & 0
        \quad & \textrm{for}\;\; \mathbf{x} \in \partial\/\Omega\/,
        \\ \rule{0ex}{2.5ex}
     \nabla \cdot \mathbf{u}                     & = & 0
        \quad & \textrm{for}\;\; \mathbf{x} \in \partial\/\Omega\/,
   \end{array}
   \right\}
\end{equation}
and
\begin{equation} \label{PPE_Theory_Splitting2}
   \left.
   \begin{array}{rcll}
     \Delta\/p                  & = &
        - \nabla \cdot \left( (\mathbf{u} \cdot \nabla)\,\mathbf{u} \right)
        + \nabla \cdot \mathbf{f}
        \quad & \textrm{for}\;\; \mathbf{x} \in \phantom{\partial}\/\Omega\/,
        \\ \rule{0ex}{2.5ex}
     \mathbf{n} \cdot \nabla\/p & = &
        \mathbf{n} \cdot \left( \mathbf{f} - \mathbf{g}_t
        + \mu\,\Delta\/\mathbf{u}
        - (\mathbf{u} \cdot \nabla)\,\mathbf{u} \right)
        \quad & \textrm{for}\;\; \mathbf{x} \in \partial\/\Omega\/.
   \end{array}
   \right\}
\end{equation}
Again, for smooth (up to the boundary) enough solutions $(\mathbf{u}\/,\,p)\/$
of the equations: the incompressible Navier-Stokes equations
(\ref{NS_Momentum}--\ref{NS_Mass}), with boundary conditions as in
(\ref{NS_Boundary_Conditions}), are equivalent to the system of equations and
boundary conditions in
(\ref{PPE_Theory_Splitting1}--\ref{PPE_Theory_Splitting2}).

For the sake of completeness, we display now the calculation showing that the
boundary condition splitting in
(\ref{PPE_Theory_Splitting1}--\ref{PPE_Theory_Splitting2}) recovers the normal
velocity boundary condition
$\mathbf{n}\cdot\mathbf{u} = \mathbf{n}\cdot \mathbf{g}\/$. To start, dot the
first equation in (\ref{PPE_Theory_Splitting1}) with the normal $\mathbf{n}\/$,
and evaluate at the boundary. This yields
\begin{equation}
   \mathbf{n} \cdot \mathbf{u}_t = \mathbf{n} \cdot \left(
   \mu\,\Delta\/\mathbf{u} - \nabla\/p -
   (\mathbf{u} \cdot \nabla)\,\mathbf{u} + \mathbf{f} \right)
   \quad \textrm{for}\;\; \mathbf{x} \in \partial\/\Omega\/.
\end{equation}
Next, eliminate $\mathbf{n} \cdot \nabla\/p\/$ from this last equation --- by
using the boundary condition for the pressure in (\ref{PPE_Theory_Splitting2}),
to obtain
\begin{equation} \label{Boundary_Evolution_Drift}
   \left(\mathbf{n}\cdot \left(\mathbf{u} - \mathbf{g}\right)\right)_t = 0
   \quad \textrm{for}\;\; \mathbf{x} \in \partial\/\Omega\/.
\end{equation}
This is a trivial ODE for the normal component of the velocity at each point
in the boundary. Thus, provided that
$\mathbf{n} \cdot (\mathbf{u} - \mathbf{g}) = 0\/$ initially, it holds for all
time.

An important final point to check is the solvability condition for the pressure
problem. Given the flow velocity $\mathbf{u}\/$ at time $t\/$, equation
(\ref{PPE_Theory_Splitting2}) is a Poisson problem with Neumann boundary
conditions for the pressure. This problem has a solution (unique up to an
additive constant) if and only if the ``flux equals source'' criteria
\begin{equation} \label{PoissonConsistency}
   \int_{\partial\/\Omega} \mathbf{n} \cdot \left( \mathbf{f} - \mathbf{g}_t +
   \mu\,\Delta\/\mathbf{u} - (\mathbf{u} \cdot \nabla)\,\mathbf{u} \right)
   \du\/A = \int_{\Omega} \nabla \cdot \left( \mathbf{f} -
   (\mathbf{u} \cdot \nabla)\,\mathbf{u}\right) \du\/V
\end{equation}
applies. This is satisfied because:
\begin{itemize}
 \item[\ref{Section:TheorySplitting}a]
 $\;\int_{\partial\/\Omega} \mathbf{n} \cdot \left( \mathbf{f}
  - (\mathbf{u} \cdot \nabla)\,\mathbf{u} \right)
  \du\/A = \int_{\Omega} \nabla \cdot \left( \mathbf{f} -
  (\mathbf{u} \cdot \nabla)\,\mathbf{u}\right) \du\/V\/$,
 \item[\ref{Section:TheorySplitting}b]
 $\;\int_{\partial\/\Omega} \mathbf{n} \cdot \Delta\/\mathbf{u} \du\/A =
  \int_{\Omega} \Delta\,\left( \nabla \cdot \mathbf{u} \right) \du\/V = 0\/$,
 \item[\ref{Section:TheorySplitting}c]
 $\;\int_{\partial\/\Omega} \mathbf{n} \cdot \mathbf{g}_t \du\/A =
  \frac{\partial}{\partial\/t}\,\int_{\partial\/\Omega} \mathbf{n} \cdot
  \mathbf{g} \du\/A = 0\/$,
\end{itemize}
where we have used Gauss' theorem, incompressibility, and
(\ref{NS_Boundary_ConditionsConstr}).


\section{Modification for Stability} \label{Section:NumSplitting}
For the system of equations in
(\ref{PPE_Theory_Splitting1}--\ref{PPE_Theory_Splitting2}), it is important to
notice that
\begin{itemize}
 \item[\ref{Section:NumSplitting}a]
 The tangential boundary condition on the flow velocity,
 $\mathbf{n}\times (\mathbf{u}-\mathbf{g}) = 0\,$ for
 $\,\mathbf{x} \in \partial\/\Omega\/$, is enforced explicitly.
 \item[\ref{Section:NumSplitting}b]
 The incompressibility condition, $\nabla \cdot \mathbf{u} = 0\,$ for
 $\,\mathbf{x} \in \Omega\/$, is enforced ``exponentially''. By this we mean
 that any errors in satisfying the incompressibility condition are rapidly
 damped, because $\phi = \nabla \cdot \mathbf{u}\/$ satisfies
 (\ref{eqn:heatdivu}--\ref{eqn:bcdivu}). Thus this condition is enforced in
 a robust way, and we do not expect it to cause any trouble for
 ``reasonable'' numerical discretizations of the equations.
 \item[\ref{Section:NumSplitting}c]
 By contrast, the normal boundary condition on the flow velocity,
 $\mathbf{n} \cdot (\mathbf{u} - \mathbf{g}) = 0\,$ for
 $\,\mathbf{x} \in \partial\/\Omega\/$, is enforced in a rather weak fashion.
 By this we mean that errors in satisfying this condition are not damped at
 all by equation (\ref{Boundary_Evolution_Drift}).  Thus, this condition lacks the 
 inherent stability provided by the heat equation.
 In practice, numerical errors add (effectively) noise to equation
 (\ref{Boundary_Evolution_Drift}), resulting in a drift of the normal velocity
 component. This can have de-stabilizing effects on the behavior of a
 numerical scheme. Hence it is a problem that must be corrected.
\end{itemize}
In this section we alter the PPE equations
(\ref{PPE_Theory_Splitting1}--\ref{PPE_Theory_Splitting2}) to address the
problem pointed out in item~\ref{Section:NumSplitting}c. We do this by adding
an appropriate ``stabilizing'' term.  The idea is similar in nature to the feedback 
controller introduced in \cite{GoldsteinHandlerSerovich93} for immersed boundary
methods.  Our goal here is to develop a pair of
differential equations --- fully equivalent to
(\ref{PPE_Theory_Splitting1}--\ref{PPE_Theory_Splitting2}) --- which are
suitable for numerical implementation.
In order to resolve the issue in item~\ref{Section:NumSplitting}c, we add a
feedback term to the equations, by altering the pressure boundary condition.
Specifically, we modify the equations from
(\ref{PPE_Theory_Splitting1}--\ref{PPE_Theory_Splitting2}) to
\begin{equation} \label{PPE_Numerical_Splitting1}
   \left.
   \begin{array}{rcll}
     \mathbf{u}_t - \mu\,\Delta\/\mathbf{u}      & = &
        - \nabla\/p - (\mathbf{u} \cdot \nabla)\,\mathbf{u} + \mathbf{f}
        \quad & \textrm{for}\;\; \mathbf{x} \in \phantom{\partial}\/\Omega\/,
        \\ \rule{0ex}{2.5ex}
     \mathbf{n} \times (\mathbf{u} - \mathbf{g}) & = & 0
        \quad & \textrm{for}\;\; \mathbf{x} \in \partial\/\Omega\/,
        \\ \rule{0ex}{2.5ex}
     \nabla \cdot \mathbf{u}                     & = & 0
        \quad & \textrm{for}\;\; \mathbf{x} \in \partial\/\Omega\/,
   \end{array}
   \right\}
\end{equation}
and
\begin{equation} \label{PPE_Numerical_Splitting2}
   \left.
   \begin{array}{rcll}
     \Delta\/p                  & = &
        - \nabla \cdot \left((\mathbf{u} \cdot \nabla)\,\mathbf{u} \right)
        + \nabla \cdot \mathbf{f}
        \quad & \textrm{for}\;\; \mathbf{x} \in \phantom{\partial}\/\Omega\/,
        \\ \rule{0ex}{2.5ex}
     \mathbf{n} \cdot \nabla\/p & = &
        \mathbf{n} \cdot \left( \mathbf{f} - \mathbf{g}_t
        + \mu\,\Delta\/\mathbf{u}
        - (\mathbf{u} \cdot \nabla)\,\mathbf{u} \right) &
        \\ \rule{0ex}{2.5ex}
                                & + &
        \lambda\,\mathbf{n} \cdot (\mathbf{u} - \mathbf{g})
        \quad & \textrm{for}\;\; \mathbf{x} \in \partial\/\Omega\/,
   \end{array}
   \right\}
\end{equation}
where $\lambda >0\/$ is a numerical parameter --- see \S~\ref{sub:selectLambda}.
This system is still equivalent to the incompressible Navier-Stokes equations
and boundary conditions in (\ref{NS_Momentum}--\ref{NS_Boundary_Conditions}),
for smooth (up to the boundary) enough solutions $(\mathbf{u}\/,\,p)\/$. The
heat equation (\ref{eqn:heatdivu}--\ref{eqn:bcdivu}) for
$\phi = \nabla \cdot \mathbf{u}\/$ still applies, while the equation for the
evolution of the normal velocity at the boundary changes from
(\ref{Boundary_Evolution_Drift}) to:
\begin{equation} \label{Boundary_Evolution_NoDrift}
   \left(\mathbf{n}\cdot \left(\mathbf{u} - \mathbf{g}\right)\right)_t =
   - \lambda\,\mathbf{n}\cdot (\mathbf{u} - \mathbf{g})
   \quad \textrm{for}\;\; \mathbf{x} \in \partial\/\Omega\/.
\end{equation}
Thus, if $\mathbf{n}\cdot (\mathbf{u} - \mathbf{g}) = 0\/$ initially, it
remains so for all times. In addition, this last equation shows that this new
system resolves the issue pointed out in item~\ref{Section:NumSplitting}c.

Finally, we check what happens to the solvability condition for the pressure
problem, given the system change above. Clearly, all we need to do is to
modify equation (\ref{PoissonConsistency}) by adding --- to its left hand
side, the term
\begin{equation} \label{PoissonConsistencyMod}
   \lambda\,\int_{\partial\/\Omega} \mathbf{n} \cdot \left( \mathbf{u}
   - \mathbf{g} \right)\,\du\/A =
   \lambda\,\int_{\Omega} \mathbf{\nabla} \cdot \mathbf{u}\,\du\/V -
   \lambda\,\int_{\partial\/\Omega} \mathbf{n} \cdot \mathbf{g}\,\du\/A
   = 0\/,
\end{equation}
where we have used incompressibility, and equation
(\ref{NS_Boundary_ConditionsConstr}). It follows that solvability remains
valid.
\begin{remark} \label{rem:JohnstonLiu2004} 
 The reformulation of the Navier-Stokes equations in
 (\ref{PPE_Numerical_Splitting1}--\ref{PPE_Numerical_Splitting2}) is similar to
 the one used by Johnston and Liu in~\cite{JohnstonLiu2004}. In their paper the
 authors propose methods of solution to the Navier-Stokes equations based on
 the equivalent system (for simplicity, we set $\mathbf{g} = 0\/$, as done
 in~\cite{JohnstonLiu2004}) where
 \begin{itemize}
  \item[(a)] \vspace*{-0.5ex}
  For $\mathbf{x} \in \phantom{\partial}\/\Omega\/$, the same equations as
  in (\ref{PPE_Numerical_Splitting1}--\ref{PPE_Numerical_Splitting2}) apply.
  \item[(b)] \vspace*{-0.5ex}
  For $\mathbf{x} \in \partial\/\Omega\/$, $\;\mathbf{u} = 0\/$ is used for
  the momentum equation.
  \item[(c)] \vspace*{-0.5ex}
  For $\mathbf{x} \in \partial\/\Omega\/$,
  $\;\mathbf{n} \cdot \nabla p = \mathbf{n} \cdot \left(
   - \mu \, \nabla \times \nabla \times \mathbf{u} + \mathbf{f} \right)\/$
  is used for the Poisson equation.
 \end{itemize}
 For this system the incompressibility condition
 $\phi = \nabla \cdot \mathbf{u} = 0\/$ follows because these equations
 yield
 \begin{itemize}
  \item[(d)] \vspace*{-0.5ex}
  $\phi_t = \mu\,\Delta\,\phi\/$ for
  $\mathbf{x} \in \Omega\/$, with
  $\mathbf{n} \cdot \nabla \, \phi = 0\/$
  for $\mathbf{x} \in \partial\/\Omega\/$.
 \end{itemize}
 Thus, if $\phi\/$ vanishes initially, it will vanish for all times.
 \begin{itemize}
  \item[(e)] \vspace*{-0.5ex}
  The \emph{main advantage} of this reformulation over the one in
  (\ref{PPE_Numerical_Splitting1}--\ref{PPE_Numerical_Splitting2}) is that
  (in general) the boundary condition $\nabla \cdot \mathbf{u} = 0\/$ in
  (\ref{PPE_Numerical_Splitting1}) may couple the components of the flow
  velocity field --- see \S~\ref{subsec:MomentumEquation}. Thus, depending on implementation, an implicit
  treatment of the viscous terms in
  (\ref{PPE_Numerical_Splitting1}--\ref{PPE_Numerical_Splitting2}) may be
  more expensive (and complicated) than for the system in items~a--c above.
  However, it is not clear to us at this moment how much of a problem this
  is. The reason is that the coupling is ``weak'', by which we mean that:
  in the $N_G \times N_G\/$ discretization matrix for the Laplacian --- where
  $N_G\/$ is the number of points in the numerical grid, the coupling induced
  by the boundary condition affects only $O(N_G^{1/2})\/$ entries in 2-D, and
  $O(N_G^{1/3})\/$ entries in 3-D --- at least with the type of discretization
  that we use in \S~\ref{subsec:MomentumEquation}. Hence, it may be possible
  to design algorithms where the extra cost is moderate, and not a
  show-stopper.

  \item[(f)] \vspace*{-0.5ex}
	In cases where the viscosity $\mu$ is small, the velocity divergence $\phi$ may not be adequately damped by the heat equation.  It is possible \cite{Henshaw1994, HenshawPetersson2003} to 
  modify the method by adding a parameter that
  forces an exponential decay of $\phi\/$ to zero. For example: modify the
  Poisson equation for the pressure to
  \[ \Delta\/p  =  - \nabla \cdot \left((\mathbf{u} \cdot \nabla)\,\mathbf{u}
     \right) + \nabla\cdot\mathbf{f} + \lambda\,\nabla \cdot \mathbf{u}\/. \]
  This then changes the system in item~d to
  \[\phi_t = \mu\,\Delta\,\phi - \lambda\,\phi\; \;\mbox{for}\;\;
    \mathbf{x} \in \Omega\/, \;\;\mbox{with}\;\;
    \mathbf{n} \cdot \nabla \, \phi = 0   \;\;\mbox{for}\;\;
    \mathbf{x} \in \partial\/\Omega\/. \]
\myremarkend
 \end{itemize}
\end{remark}
\begin{remark} \label{rem:variablelambda}
  It would be nice to be able to use $\lambda = \lambda(\mathbf{x})\/$, so as
  to optimize the implementation of the condition
  $\mathbf{n}\cdot \mathbf{u} = \mathbf{n}\cdot \mathbf{g}\/$ for different
  points along $\partial\/\Omega\/$.  However, this is not a trivial extension,
  since it destroys the solvability condition for the pressure. Thus, other
  (compensating) corrections are needed as well. We postpone the study of this
  issue for future work.\myremarkend
\end{remark}

\begin{remark} \label{rem:generalIC}
  As show earlier in (\ref{PoissonConsistency}) and (\ref{PoissonConsistencyMod}), the 
  validity of the solvability condition, for the problem in (\ref{PPE_Numerical_Splitting2}), 
  relies on the velocity field satisfying the incompressibility constraint 
  $\nabla \cdot \mathbf{u} = 0\/$.  On the other hand, in the course of a 
  numerical calculation, one may not always start with a numerically divergence free velocity field.  This consideration motivates the following 
  theoretical question: \emph{can the equations in 
  (\ref{PPE_Numerical_Splitting1}--\ref{PPE_Numerical_Splitting2}) be
  modified, so they make sense even for $\nabla \cdot \mathbf{u} \neq 0\/$?}
  In \ref{solvabilityBreakdown} we show that this is possible. \myremarkend
\end{remark}

\begin{remark} \label{rem:singleDerivativeRHS}
    The contribution from the nonlinear advection term to the right hand side of equation (\ref{PPE_Numerical_Splitting2}) involves second derivatives of the velocity field.  For some numerical implementations, such as finite element methods, there are advantages to reducing the number of derivatives in the Poisson equation source terms.  For example, reducing the number of derivatives allows one to choose a finite element basis with less regularity.  As noted by \cite{JohnstonLiu2004, Henshaw1994}, one may rewrite the nonlinear terms to contain one derivative
\begin{equation} \label{PPE_Numerical_SplittingOneDerivative}
   \left.
   \begin{array}{rcll}
     \Delta\/p                  & = &
        \frac{1}{2} (\nabla \cdot \mathbf{u})^2 - (\nabla \mathbf{u}) : (\nabla \mathbf{u})^T
        + \nabla \cdot \mathbf{f}
        \quad & \textrm{for}\;\; \mathbf{x} \in \phantom{\partial}\/\Omega\/,
        \\ \rule{0ex}{2.5ex}
     \mathbf{n} \cdot \nabla\/p & = &
        \mathbf{n} \cdot \left( \mathbf{f} - \mathbf{g}_t
        + \mu\,\Delta\/\mathbf{u}
        - (\mathbf{u} \cdot \nabla)\,\mathbf{u} \right) &
        \\ \rule{0ex}{2.5ex}
                                & + &
        \lambda\,\mathbf{n} \cdot (\mathbf{u} - \mathbf{g})
        \quad & \textrm{for}\;\; \mathbf{x} \in \partial\/\Omega\/,
   \end{array}
   \right\}
\end{equation}
By a similar calculation to the one in the appendix of \cite{JohnstonLiu2004}, we can show that the pressure equation (\ref{PPE_Numerical_SplittingOneDerivative}) is equivalent to Navier-Stokes for smooth enough solutions.

\end{remark}
%
%
\subsection{Selection of the parameter $\lambda\/$.} \label{sub:selectLambda}
For numerical purposes, here we address the issue of how large $\lambda\/$
should be, by using a simple model for the flow's normal velocity drift.
Notice that no precision is needed for this calculation, just order of
magnitude. In actual practice, one can monitor how well the normal velocity
satisfies the boundary condition, and increase $\lambda\/$ if needed. In
principle one should be wary of using large values for $\lambda\/$, since this
will yield stiff behavior in time.  However, the calculation below shows that  $\lambda\/$ does not need to be very large, and does not depend on the grid size $\Delta\/x\/$.

It seems reasonable to assume that one can model how the numerical errors
affect the ODE (\ref{Boundary_Evolution_NoDrift}) for
$\mathcal{E} = \mathbf{n}\cdot (\mathbf{u}-\mathbf{g})\/$, by perturbing the
coefficients of the equation, and adding a forcing term to it. Hence we modify
equation (\ref{Boundary_Evolution_NoDrift}) as follows
\begin{equation} \label{Boundary_Evolution_Forcing1}
 \mathcal{E}_t = - \lambda\,c_p\,\mathcal{E} + \epsilon\,\gamma
 \quad \textrm{for}\;\; \mathbf{x} \in \partial\/\Omega\/,
\end{equation}
where $\epsilon \ll 1\/$ characterizes the size of the errors (determined by
the order of the numerical method), while
$c_p = c_p(\mathbf{x}\/,\,t) = 1 + O(\epsilon)\/$ and
$\gamma = \gamma(\mathbf{x}\/,\,t) = O(1)\/$ are functions encoding the
numerical errors. What exactly they are depends on the details of the
numerical discretization, but for this calculation we do not need to know
these details. All we need is that\footnote{For smooth enough solutions,
   where the truncation errors are controlled by some derivative of the
   solution.}
\begin{equation} \label{Boundary_Evolution_Forcing2}
   0 < C_M \leq c_p \quad \textrm{and} \quad |\gamma| \leq \Gamma\/,
\end{equation}
where $C_M\/$ and $\Gamma\/$ are some positive constants, with
$C_M \approx 1\/$, and $\Gamma = O(1)\/$ --- but not necessarily close to one.

The solution to (\ref{Boundary_Evolution_Forcing1}) is given by
\begin{equation} \label{Boundary_Evolution_Forcing3}
 \mathcal{E} = \mathcal{E}_0\,e^{-\lambda\,I_1} + \epsilon\,\underbrace{
 \int_0^t \gamma(\mathbf{x}\/,\,s)\,e^{-\lambda\,I_2}\,\du\/s}_J\/,
\end{equation}
where $\mathcal{E}_0\/$ is the initial value,
$I_1 = I_1(\mathbf{x}\/,\,t) = \int_0^t c_p(\mathbf{x}\/,\,s)\/\du\/s\/$, and
$I_2 = I_1(\mathbf{x}\/,\,t) - I_1(\mathbf{x}\/,\,s)\/$. The crucial term is
$J\/$, since the first term decreases in size, and starts at the initial
value. However
\begin{equation} \label{Boundary_Evolution_Forcing4}
 |J| \leq \Gamma\,\int_0^t e^{-\lambda\,C_M\,(t-s)}\/\du\/s \leq
 \frac{\Gamma}{\lambda\,C_M}\/.
\end{equation}

Within the framework of a numerical scheme, the normal boundary velocity may
deviate from the prescribed velocity by some acceptable error $\delta\/$. Thus
we require $\epsilon\,J = O(\delta)\/$ or less, which --- given
(\ref{Boundary_Evolution_Forcing4}) --- will be satisfied if
\begin{equation} \label{Lambda_Magnitude}
 \lambda \sim \frac{\epsilon\,\Gamma}{\delta\,C_M} \approx
 \frac{\epsilon\,\Gamma}{\delta}\/, \quad \textrm{or larger.}
\end{equation}
For the second order numerical scheme in \S~\ref{Section:NumScheme}, it is
reasonable to expect that $\epsilon = (\Delta\/x)^2\/$, and to require that
$\delta = (\Delta\/x)^2\/$. Then (\ref{Lambda_Magnitude})
reduces to $\lambda \geq \Gamma\/$. Of course, we do not know (a priori) what
$\Gamma\/$ is; this is something that we need to find by numerical
experimentation --- see the first paragraph in this \S~\ref{sub:selectLambda}.
For the numerical calculations reported in \S~\ref{Section:Implementation}, we
found that values in the range $10 \leq \lambda \leq 100\/$ gave good results.

For this scheme, comparing the time step restriction imposed by feedback $(\Delta t_{\lambda} < O(1/\lambda))$ to the one imposed by diffusion $\Delta t_{\mu} < O((\Delta x)^2/\mu)$, results in a stiffness ratio that scales as
\begin{eqnarray}
	\frac{\Delta t_{\mu}}{\Delta t_{\lambda}} \sim \frac{\lambda}{\mu} \Delta x^2.
\end{eqnarray}
Hence for low to moderate Reynolds numbers, $\lambda$ can be choosen quite large without effecting the stiffness.


\section{Stability of Semi-implicit Schemes} \label{Section:SemiImplicitStability}
In the case of moderate to low Reynolds number flows, the stiff viscosity term $\mu \Delta \mathbf{u}$ requires very small time steps when treated explicitly.  Hence, there is a large practical interest in treating $\mu \Delta \mathbf{u}$ implicitly while keeping the associated pressure explicit.  In this section we write down and analyze semi-implicit schemes for the PPE splitting (\ref{PPE_Theory_Splitting1}--\ref{PPE_Theory_Splitting2}).  Inspired by many of the ideas from \cite{JohnstonLiu2004}, we then show the subsequent schemes are unconditionally stable.  

To analyze the stability of the semi-implicit schemes, we require several properties of the Hodge-Helmholtz decomposition.  Given $\mathbf{w} \in L^2$, then $\mathbf{w}$ has a unique orthogonal decomposition as
\begin{eqnarray}
	\mathbf{w} = \mathbf{a} + \nabla b, 
\end{eqnarray}
where $b$ is determined by $\Delta b = \nabla \cdot \mathbf{w}$ with boundary conditions  $\mathbf{n} \cdot \nabla b = \mathbf{n} \cdot \mathbf{w}$.  Hence $\mathbf{a}$ is divergence free with zero normal boundary component.  Moreover the component $\mathbf{a}$ is orthogonal in the $L^2$ norm to every gradient field (i.e., $\langle\mathbf{a}, \nabla \phi\rangle = \int_{\Omega} \mathbf{a}\cdot \nabla \phi \du V = 0$ is the standard inner product on $L^2(\Omega)$). 
We may therefore write $\nabla b = \mathcal{Q} \mathbf{w}$ and $\mathbf{a} = \mathcal{P} \mathbf{w}$ where $\mathcal{P}$ and $\mathcal{Q}$ are complementary orthogonal projections (ie. $\mathcal{P} + \mathcal{Q} = \mathcal{I}$, $\mathcal{P}^2 = \mathcal{P}$ and $\mathcal{P} = \mathcal{P}^\dagger$).  As a result, the following identities hold for any vector fields $\mathbf{w}$ and $\mathbf{v}$
\begin{eqnarray}
	\langle \mathbf{w}, \mathcal{P}\mathbf{v} \rangle &=& \langle \mathcal{P} \mathbf{w}, \mathcal{P}\mathbf{v} \rangle \\
&=& \langle \mathcal{P} \mathbf{w}, \mathbf{v} \rangle, 
\end{eqnarray}
and
\begin{eqnarray}
	\langle \mathbf{w}, \mathbf{w} \rangle &=& \langle \mathcal{P} \mathbf{w}, \mathcal{P}\mathbf{w} \rangle + \langle \mathcal{Q}\mathbf{w}, \mathcal{Q}\mathbf{w} \rangle.
\end{eqnarray}

To motivate the stability proofs for semi-implicit time discretizations, note that the linear Navier-Stokes equations can also be written in an equivalent projection form
\begin{equation} \label{Projected_Equation}
   \left.
   \begin{array}{rcll}
     \mathbf{u}_t    & = & \mu \mathcal{P} \Delta\/\mathbf{u}
        \quad & \textrm{for}\;\; \mathbf{x} \in \phantom{\partial}\/\Omega\/,
        \\ \rule{0ex}{2.5ex}
     \mathbf{n} \times \mathbf{u} & = & 0
        \quad & \textrm{for}\;\; \mathbf{x} \in \partial\/\Omega\/,
   \end{array}
   \right\}
\end{equation}
Dotting equation (\ref{Projected_Equation}) through by\footnote{To obtain the energy bound (\ref{ContinuousEnergyBound}) by a straightforward calculation, we have kept the $\nabla \cdot \mathbf{u}$ terms throughout the computation even though they are in fact zero.} $-\Delta \mathbf{u} = \nabla \times \nabla \times \mathbf{u} - \nabla (\nabla \cdot \mathbf{u} )$ and integrating yields
\begin{eqnarray}
	\frac{1}{2} \frac{\partial}{\partial t}\Big( \langle \nabla \times \mathbf{u}, \nabla \times \mathbf{u} \rangle + \langle \nabla \cdot \mathbf{u}, \nabla \cdot \mathbf{u} \rangle \Big) &=& - \mu \langle \Delta \mathbf{u} , \mathcal{P} \Delta \mathbf{u} \rangle. 
\end{eqnarray}
Here the first line follows directly from integration by parts, noting that the boundary integrals vanish.  The last equation can be written in a more compact form:
\begin{eqnarray} \label{ContinuousEnergyBound}
	\frac{1}{2}\frac{\partial}{\partial t}\Big(||\nabla \times \mathbf{u}||^2 + ||\nabla \cdot \mathbf{u}||^2 \Big)&=&  -\mu ||\mathcal{P} \Delta \mathbf{u} ||^2. 
\end{eqnarray}
From this it follows that the curl and divergence are bounded in the $L^2$ sense.
As a first example of a semi-implicit scheme, we analyze the backward Euler discretization for the Stokes equation:
\begin{equation} \label{StokesSimpleEuler1}
   \left.
   \begin{array}{rcll}
     \mathbf{u}^{n+1} - \mathbf{u}^n - \Delta t \, \mu\,\Delta\/\mathbf{u}^{n+1}      & = &  - \Delta t \, \nabla\/p^n
        \quad & \textrm{for}\;\; \mathbf{x} \in \phantom{\partial}\/\Omega\/,
        \\ \rule{0ex}{2.5ex}
     \mathbf{n} \times \mathbf{u}^{n+1} & = & 0
        \quad & \textrm{for}\;\; \mathbf{x} \in \partial\/\Omega\/,
        \\ \rule{0ex}{2.5ex}
     \nabla \cdot \mathbf{u}^{n+1}                     & = & 0
        \quad & \textrm{for}\;\; \mathbf{x} \in \partial\/\Omega\/,
   \end{array}
   \right\}
\end{equation}
and
\begin{equation} \label{StokesSimpleEuler2}
   \left.
   \begin{array}{rcll}
     \Delta\/p^n                  & = & 0
        \quad & \textrm{for}\;\; \mathbf{x} \in \phantom{\partial}\/\Omega\/,
        \\ \rule{0ex}{2.5ex}
     \mathbf{n} \cdot \nabla\/p^n & = &
        \mathbf{n} \cdot \mu\,\Delta\/\mathbf{u}^n
        \quad & \textrm{for}\;\; \mathbf{x} \in \partial\/\Omega\/.
   \end{array}
   \right\}
\end{equation}

We note that the scheme (\ref{StokesSimpleEuler1})--(\ref{StokesSimpleEuler2}) preserves the divergence free condition on the velocity exactly.  Assume $\nabla \cdot \mathbf{u}^{n} = 0$ for $\mathbf{x} \in \Omega$.  Taking the divergence of (\ref{StokesSimpleEuler1}), and setting $\phi(\mathbf{x}) = \nabla \cdot \mathbf{u}^{n+1}$ yields
\begin{equation} \label{Divergence}
   \left.
   \begin{array}{rcll}
      - \Delta t \, \mu\,\Delta\/ \phi     & = &  -\phi
        \quad & \textrm{for}\;\; \mathbf{x} \in \phantom{\partial}\/\Omega\/,
        \\ \rule{0ex}{2.5ex}
     \phi                    & = & 0
        \quad & \textrm{for}\;\; \mathbf{x} \in \partial\/\Omega\/,
   \end{array}
   \right\}
\end{equation}
Since $-\Delta$ does not have a negative eigenvalue, $\phi = 0$ for all $\mathbf{x} \in \Omega$ is the only solution to (\ref{Divergence}).  Consequently, the pressure equation for $p^{n+1}$ automatically satisfies the Neumann consistency condition and is therefore well defined.  As a direct consequence of preserving the divergence condition, the pressure as defined by (\ref{StokesSimpleEuler2}) is equivalent to the projection:
\begin{eqnarray} \label{PressureAsHodgeProjection}
	\nabla p^n = \mathcal{Q} \Delta \mathbf{u}^n.
\end{eqnarray}

\begin{remark}
Since the scheme (\ref{StokesSimpleEuler1})--(\ref{StokesSimpleEuler2}) preserves the divergence constraint, there is never the need for a projection step (ie. one may associate $\mathbf{u}^{n+1} = \mathbf{u^*}^{n+1}$, where $\mathbf{u^*}^{n+1}$ is the conventional intermediate velocity field with non-zero divergence). \myremarkend
\end{remark} 

A stability proof for equations (\ref{StokesSimpleEuler1})--(\ref{StokesSimpleEuler2}) closely follows the one for periodic channel flow in \cite{JohnstonLiu2004}.  Specifically, dot both sides of equation (\ref{StokesSimpleEuler1}) by $-\Delta(\mathbf{u}^{n+1}+\mathbf{u}^n)$, and integrate (this is the discrete analog of the steps required to derive equation (\ref{ContinuousEnergyBound})).  The first two terms on the left hand side of equation (\ref{StokesSimpleEuler1}) become
\begin{eqnarray}
	\langle \mathbf{u}^{n+1}-\mathbf{u}^n, -\Delta(\mathbf{u}^{n+1} + \mathbf{u}^{n}) \rangle &=& ||\nabla \times \mathbf{u}^{n+1}||^2 + ||\nabla \cdot \mathbf{u}^{n+1}||^2 \\
 &-& ||\nabla \times \mathbf{u}^n||^2 - ||\nabla \cdot \mathbf{u}^n||^2. 
\end{eqnarray}
In addition, we have
\begin{eqnarray}
	\langle \Delta \mathbf{u}^{n+1}, \Delta (\mathbf{u}^{n+1} + \mathbf{u}^{n})\rangle = \frac{1}{2}||\Delta \mathbf{u}^{n+1}+\Delta \mathbf{u}^{n}||^2 + \frac{1}{2}||\Delta \mathbf{u}^{n+1}||^2 - \frac{1}{2}||\Delta \mathbf{u}^{n}||^2,
\end{eqnarray}
and
\begin{eqnarray}
	\langle \mathcal{Q}\Delta \mathbf{u}^{n}, \Delta (\mathbf{u}^{n+1} + \mathbf{u}^{n})\rangle = \frac{1}{2}||\mathcal{Q}\Delta \mathbf{u}^{n+1}+\mathcal{Q}\Delta \mathbf{u}^{n}||^2 + \frac{1}{2}||\mathcal{Q}\Delta\mathbf{u}^{n}||^2 - \frac{1}{2}||\mathcal{Q}\Delta \mathbf{u}^{n+1}||^2.
\end{eqnarray}
For brevity, we also introduce the energy $\mathcal{E}$ where
\begin{eqnarray}
	\mathcal{E}[\mathbf{u}] = ||\nabla \times \mathbf{u}||^2 + ||\nabla \cdot \mathbf{u}||^2 + \mu \Delta t||\Delta\mathbf{u}||^2 + \mu \Delta t ||\mathcal{Q}\Delta \mathbf{u}||^2.
\end{eqnarray}
Finally, combining everything, we have
\begin{eqnarray}
	\mathcal{E}[\mathbf{u}^{n+1}] - \mathcal{E}[\mathbf{u}^{n}] &=& -\frac{\mu \Delta t}{2}||\Delta \mathbf{u}^{n+1}+\Delta \mathbf{u}^{n}||^2 + \frac{\mu \Delta t}{2}||\mathcal{Q}\Delta \mathbf{u}^{n+1}+ \mathcal{Q}\Delta \mathbf{u}^{n}||^2 \\
&=& -\frac{\mu \Delta t}{2} ||\mathcal{P}\Delta \mathbf{u}^{n+1}+\mathcal{P}\Delta \mathbf{u}^{n}||^2, 
\end{eqnarray}
so that
\begin{eqnarray} \label{DiscreteEnergyBound}
	\mathcal{E}[\mathbf{u}^{n+1}] \leq \mathcal{E}[\mathbf{u}^{n}].
\end{eqnarray}
For simply connected domains, the energy bound (\ref{DiscreteEnergyBound}) implies the scheme  (\ref{StokesSimpleEuler1})--(\ref{StokesSimpleEuler2}) is unconditionally stable.  Unfortunately the bound (\ref{DiscreteEnergyBound}) is not sufficient to prove stability in periodic domains, or domains with holes.  The reason here is that $\mathcal{E}[\mathbf{h}] = 0$ for any vector field $\mathbf{h}$ that has $\nabla \cdot \mathbf{h} = \nabla \times \mathbf{h} = 0$ in $\Omega$, and $\nabla \cdot \mathbf{h} = \mathbf{n} \times \mathbf{h} = 0$ in $\partial \Omega$.  In simply connected domains, the only such vector fields are $\mathbf{h} = 0$.  On the other hand, in periodic domains or domains with holes, there exist vector fields $\mathbf{h}\neq 0$, which satisfy $\nabla \cdot \mathbf{h} = 0$ and $\mathbf{n}\times \mathbf{h} = 0$ on $\partial \Omega$ and $\mathcal{E}[\mathbf{h}] = 0$.  Here the energy bound (\ref{DiscreteEnergyBound}) does not control the growth of these modes.  We note however that such vector fields $\mathbf{h}$ have non zero flux $\mathbf{h}\cdot \mathbf{n} \neq 0$ at the boundary.  Hence, we expect that the modified scheme (\ref{PPE_Numerical_Splitting1})--(\ref{PPE_Numerical_Splitting1}) with addition of the term $\lambda \mathbf{n}\cdot \mathbf{u}$ in the boundary condition for the pressure will stabilize the growth of such modes.

\subsection{Stability of a second order scheme}
Following a procedure analogous to the one in \cite{JohnstonLiu2004}, we discuss the stability of a second order Crank-Nicholson scheme where we treat the pressure with a second order Adams-Bashforth extrapolation:
\begin{equation} \label{StokesCN}
   \left.
   \begin{array}{rcll}
     \mathbf{u}^{n+1} - \mathbf{u}^n & = & \frac{\mu}{2}\Delta t \,\big( \Delta\/\mathbf{u}^{n+1} + \Delta\/\mathbf{u}^{n} \big) - \frac{3}{2} \Delta t \nabla\/p^n + \frac{1}{2} \Delta t \nabla\/p^{n-1}
        \quad & \textrm{for}\;\; \mathbf{x} \in \phantom{\partial}\/\Omega\/,
        \\ \rule{0ex}{2.5ex}
     \mathbf{n} \times \mathbf{u}^{n+1} & = & 0
        \quad & \textrm{for}\;\; \mathbf{x} \in \partial\/\Omega\/,
        \\ \rule{0ex}{2.5ex}
     \nabla \cdot \mathbf{u}^{n+1}                     & = & 0
        \quad & \textrm{for}\;\; \mathbf{x} \in \partial\/\Omega\/,
   \end{array}
   \right\}
\end{equation}
where $p^n$ is given by (\ref{StokesSimpleEuler2}).  Proceeding with a normal mode analysis\footnote{We will not discuss here the question of whether the resulting eigenvectors form a complete basis.}, we set $\mathbf{u}^n = \alpha^n \mathbf{\tilde{u}}$ where $\mathbf{\tilde{u}}$ satisfies the boundary conditions in (\ref{StokesCN}), and $\alpha$ is an eigenvalue of the time stepping operator.  Upon substitution we obtain
\begin{eqnarray} \label{Eigenvalues1}
     (\alpha^2 - \alpha) \mathbf{\tilde{u}} & = & \frac{\mu}{2}\Delta t  (\alpha^2 + \alpha) \Delta\/\mathbf{\tilde{u}} - \frac{3}{2} \Delta t \alpha \mathcal{Q}\Delta \mathbf{\tilde{u}} + \frac{1}{2} \Delta t \mathcal{Q}\Delta \mathbf{\tilde{u}}.
\end{eqnarray}
In analogy with the steps in the previous section, we dot equation (\ref{Eigenvalues1}) through by $-\Delta \mathbf{\tilde{u}}$ and integrate by parts to obtain
\begin{eqnarray} \label{Eigenvalues2}
     2(\alpha^2 - \alpha)\big( ||\nabla \times \mathbf{\tilde{u}}||^2 + ||\nabla \cdot \mathbf{\tilde{u}}||^2 \big) + \mu \Delta t  (\alpha^2 + \alpha) ||\Delta\mathbf{\tilde{u}}||^2 = \mu \Delta t (3\alpha - 1) ||\mathcal{Q}\Delta \mathbf{\tilde{u}}||^2.
\end{eqnarray}
Equation (\ref{Eigenvalues2}) is now a quadratic of the form $a \alpha^2 - b \alpha + c = 0$ for the eigenvalues $\alpha$ where
\begin{eqnarray}
    a &=& 2 ||\nabla \times \mathbf{\tilde{u}}||^2 + 2||\nabla \cdot \mathbf{\tilde{u}}||^2  + \mu \Delta t ||\Delta \mathbf{\tilde{u}}||^2, \\
    b &=& 2 ||\nabla \times \mathbf{\tilde{u}}||^2 + 2||\nabla \cdot \mathbf{\tilde{u}}||^2  - \mu \Delta t ||\Delta \mathbf{\tilde{u}}||^2 + 3\mu \Delta t ||\mathcal{Q}\Delta \mathbf{\tilde{u}}||^2, \\
    c &=& \mu \Delta t ||\mathcal{Q}\Delta \mathbf{\tilde{u}}||^2.
\end{eqnarray}
Since the coefficients $a, b, c$ satisfy the following inequalities $0 \leq c < a$ and $|b| < a + c$, it follows \cite{Strikwerda} that the eigenvalues $\alpha$ lie within the unit circle.  The preceding argument indicates that the normal modes $\mathbf{\tilde{u}}$ remain bounded for the second order Crank-Nicholson scheme (\ref{StokesCN}).  Hence the scheme is stable for domains with simple geometries\footnote{Note that in domains which are not simply connected, and have nontrivial homology groups, there is a family of harmonic modes which satisfy $\Delta \mathbf{\tilde{u}} = 0$.  For these modes one can check that $a = b = c = 0$ and the argument fails to show that $|\alpha| < 1$}. 


\section{A First-Order Explicit Scheme on Irregular Domains} \label{Section:NumScheme}
In this section we outline an explicit numerical scheme for solving the
coupled differential equations
(\ref{PPE_Numerical_Splitting1}--\ref{PPE_Numerical_Splitting2}) on a
two-dimensional irregular domain. 
To achieve such a scheme, we decouple
the pressure and velocity fields, and explicitly treat each term in the time
evolution of (\ref{PPE_Numerical_Splitting1}--\ref{PPE_Numerical_Splitting2}).
Specifically, since both the right hand side and boundary conditions of
equation (\ref{PPE_Numerical_Splitting2}) depend solely on $\mathbf{u}\/$, we
may view the pressure as a computable functional of the velocity,
$p = p[\mathbf{u}]\/$. The computation of $p[\mathbf{u}]\/$ requires the
solution of a Poisson equation with Neumann boundary conditions. With this in
mind, the momentum equation then has the form
$\mathbf{u}_t = \mathbf{F}[\mathbf{u}]\/$, where $\mathbf{F}[\mathbf{u}]\/$
has a complicated, yet numerically computable form:
\begin{equation}\label{eqn:functionalF}
   \mathbf{F}[\mathbf{u}] = \mu\,\Delta\/\mathbf{u} -\nabla\/p[\mathbf{u}] -
   (\mathbf{u} \cdot \nabla)\,\mathbf{u} + \mathbf{f}
   \quad \textrm{for}\;\; \mathbf{x} \in \Omega\/.
\end{equation}
We now use an explicit forward Euler scheme to discretize the time evolution
for (\ref{PPE_Numerical_Splitting1}--\ref{PPE_Numerical_Splitting2}), paired
with an appropriate discretization in space described later in this section.
This yields the scheme
\begin{equation} \label{PPE_Implementation_Splitting1}
   \frac{1}{\Delta\/t}\,\left( \mathbf{u}^{n+1} - \mathbf{u}^n \right) =
        \mu\,\Delta\/\mathbf{u}^n
      - \nabla\/p^n - (\mathbf{u}^n \cdot \nabla)\,\mathbf{u}^n + \mathbf{f}^n
        \quad \textrm{for}\;\; \mathbf{x} \in \Omega\/,
\end{equation}
with boundary conditions
$\mathbf{n} \times \mathbf{u} = \mathbf{n} \times \mathbf{g}\/$ and
$\nabla \cdot \mathbf{u} = 0\/$ for $\mathbf{x} \in \partial\/\Omega\/$,
where the pressure is given by
\begin{equation} \label{PPE_Implementation_Splitting2}
   \Delta\/p^n = - \nabla \cdot \left( (\mathbf{u}^n \cdot \nabla)\,
   \mathbf{u}^n \right) + \nabla \cdot \mathbf{f}^n
   \quad \textrm{for}\;\; \mathbf{x} \in \Omega\/,
\end{equation}
with the boundary condition
\[
   \mathbf{n} \cdot \nabla\/p^n =
   \mathbf{n} \cdot \left( \mathbf{f}^n - \mathbf{g}_t^n
   + \mu\,\Delta\/\mathbf{u}^n - (\mathbf{u}^n \cdot \nabla)\,\mathbf{u}^n
   \right) + \lambda\,\mathbf{n} \cdot (\mathbf{u}^n - \mathbf{g}^n)
\]
for $\mathbf{x} \in \partial\/\Omega\/$. Here, starting with the initial data
$\mathbf{u}^0\/$, a superscript $n\/$ is used to denote a variable at time
$t = n\,\Delta\/t$, where $0 < \Delta\/t \ll 1\/$ is the time step.

\begin{remark} \label{rem:HigherOrderTime}
  Our purpose here is to illustrate the new approach with a simple scheme that
  does not obscure the ideas in the method with technical complications. Hence,
  the scheme here is first order in time (explicit) and second order in space,
  with the stability restriction $\Delta\/t \propto (\Delta\/x)^2\/$.  However,
  unlike projection methods and other approaches commonly used to solve the
  Navier-Stokes equations --- see \S~\ref{sec:intro}, this new formulation does
  not seem to have any inherent order limitations. Unfortunately, the
  Navier-Stokes equations are stiff and nonlinear, which means that the fact
  that higher order extensions are possible does not mean that they are
  trivial. \myremarkend
\end{remark}
%
\subsection{Space grid and discretization.} \label{sub:discretization}
To discretize the equations in space, we use finite differences over a
cartesian, square ($\Delta\/x = \Delta\/y$), staggered grid. The pressure
values are stored at the nodes of the grid, while the horizontal and vertical
components of the velocity are stored at the mid-points of the edges connecting
the grid nodes (horizontal component on the horizontal edges and vertical
component on the vertical edges).

\begin{figure}[htb!]
\centering
\includegraphics[width = \textwidth]{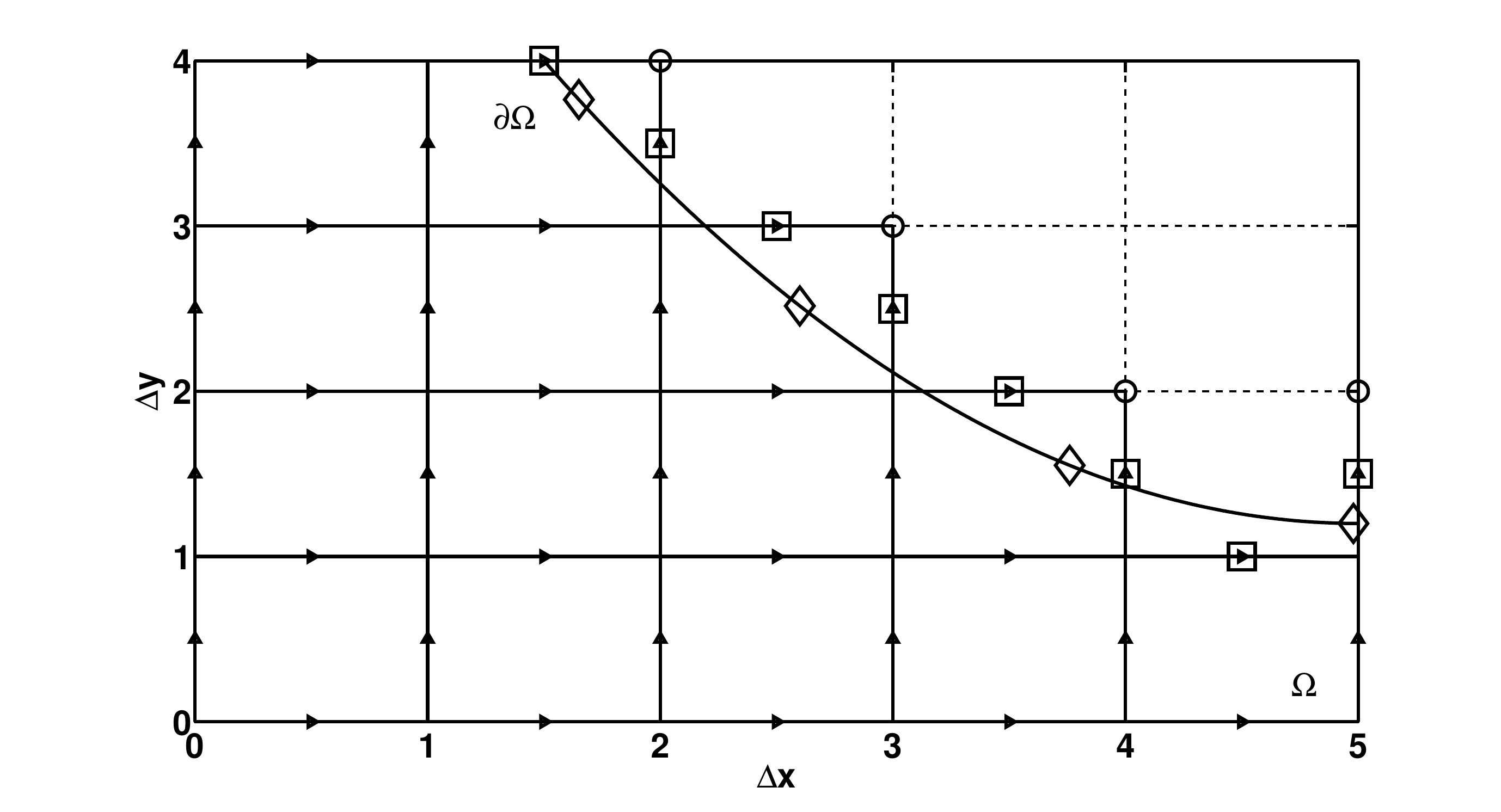}
\caption{This plot shows the staggered grid, and the boundary. The numerical
   pressure values correspond to the graph nodes, while the velocities (arrows)
   correspond to the edge midpoints. Here the circles ($\circ\/$) and squares
   ($\Box\/$) denote ghost pressure points, and boundary velocities,
   respectively. These are used to implement the boundary conditions in
   the Poisson and momentum equations, respectively. The diamonds
   ($\Diamond\/$) denote the points $(x_b\/,\,y_b)_i\/$, used to represent
   the boundary $\partial\/\Omega\/$.}
\label{fig:StaggeredMesh}
\end{figure}

When handling an arbitrary curved boundary, we do not conform the boundary to
the grid, but rather we immerse it within the regular mesh --- see
figure~\ref{fig:StaggeredMesh}.  Then, to numerically describe the domain
boundary, we identify a set $\mathcal{C}_b\/$ of $N_e\/$ points in
$\partial\/\Omega\/$, say
$\mathbf{x}_{b\,j} = (x_b\/,\, y_b)_j$ for $1 \leq j \leq N_e\/$ --- see
item~\ref{Section:NumScheme}c below. These $N_e\/$ points are located at
$O(\Delta\/x)\/$ distances apart, so that the resolution of the boundary is
comparable with that of the numerical grid.
%

For the heat equation $T_t = \mu\,\Delta\/T\/$, the stability restriction for
the standard scheme using a 5 point centered differences approximation for the
Laplacian, and forward Euler in time, is
\begin{equation}\label{eqn:stabilitycrit}
   \Delta\/t \leq C\,\frac{(\Delta\/x)^2}{\mu}\/,
\end{equation}
where $C = \frac{1}{2D}\/$ and $D = 1\/,\,2\/,\,\ldots\/$ is the space
dimension. Since the method described here uses exactly the same approach to
advance the velocity flow field $\mathbf{u}\/$, we expect the same restriction
(with, perhaps, a different constant $C\/$) to apply. For the 2D numerical
calculations presented in \S~\ref{Section:Implementation}, we found that the
algorithm was stable with $C \leq 0.2\/$, while $C \geq 0.3 \/$ generally
produced unstable behavior. Below we separately address the numerical
implementation of equations (\ref{PPE_Implementation_Splitting1}) and
(\ref{PPE_Implementation_Splitting2}).
%
%
\subsection{Poisson Equation} \label{subsec:PoissonSolver}
To solve the pressure Poisson equation (\ref{PPE_Implementation_Splitting2})
with Neumann boundary conditions, we use the ghost point idea with a method 
by Greenspan \cite{Greenspan} to implement the boundary conditions. 
As discussed above, we embed the
domain $\Omega\/$ within a cartesian square grid, and classify the
computational points in the grid into inner and ghost points (see next
paragraph). Then: (i) for the inner points we discretize the Laplace operator
using the standard 5 point centered differences stencil. (ii) For the ghost
points we obtain equations from an appropriate discretization of the Neumann
boundary condition. In particular, if there are $N_a\/$ inner points, and
$N_b\/$ ghost points, then the discretized pressure is represented by a vector
$\widehat{p} \in \mathbbm{R}^N\/$ where $N = N_a + N_b\/$.

The computational domain, $\mathcal{C}_p\/$, for the pressure is comprised
by the \emph{inner} points and the \emph{ghost} points, defined as follows:
\begin{itemize}
 \item[\ref{Section:NumScheme}a]
 The \emph{inner} points are the points in the cartesian grid located inside
 $\Omega\/$.
 \item[\ref{Section:NumScheme}b]
 The \emph{ghost} points are the points in the cartesian grid that are either
 outside $\Omega\/$ or on $\partial\/\Omega\/$, and which are needed to
 complete the 5 point stencil for some inner point.
 \item[\ref{Section:NumScheme}c]
 Construct the set $\mathcal{C}_b\/$ used to track the boundary
 $\partial\/\Omega\/$ as follows: For every pressure ghost point in the
 grid: $\mathbf{x}^{pg}_j$, $1 \leq j \leq N_b\/$, select the closest point in
 the boundary $\mathbf{x}_{b\,j} \in \partial\,\Omega\/$. We will use these
 points to produce equations that approximate the Neumann boundary conditions
 (one equation per point).  Thus, with this choice for $\mathcal{C}_b\/$, $N_e = N_b\/$.
\end{itemize}
It follows that, on the computational domain $\mathcal{C}_p\/$, the Poisson
equation (including the boundary conditions) is discretized by $N_a\/$
equations inside $\Omega\/$, and $N_b\/$ equations derived from the boundary
conditions. Specifically
\begin{itemize}
 \item[\ref{Section:NumScheme}d]
 We use the standard 5 point centered differences stencil for the Laplacian, to
 obtain one equation for each of the $N_a\/$ inner points. Similarly, we use a
 second order approximation for the forcing terms on the right hand side of
 the equation.
 \item[\ref{Section:NumScheme}e]
 We construct one boundary condition equation for each pressure ghost point
 $\mathbf{x}^{pg}_j$, $1 \leq j \leq N_b\/$. This is done by building an extrapolation operator which acts on 6
 nearby points from the computational domain to obtain a second order
 approximation to the normal derivative at the corresponding point
 $\mathbf{x}_{b\,j} \in \partial\,\Omega\/$ --- see
 item~\ref{Section:NumScheme}c above.  The 6 points are comprised of the 5 point stencil
used to approximate the Laplacian, plus one of the four next closest points
to the center point.  The extrapolation operator, which contains 6 coefficients, is then build to linearly interpolate the normal derivative using local values of the pressure.  Hence each boundary point provides one equation coupling six values of the pressure.

We use a similar extrapolation operator to approximate
 the terms on the right hand side of the Neumann boundary conditions.

 Hence, each boundary condition equation involves the ghost point, points
 inside the domain, and (possibly) other nearby ghost points.
\end{itemize}
Thus the Poisson equation can be written in a matrix form with the following
structure
\begin{equation} \label{PoissonMatrixEq1}
 \left(\begin{array}{c} L \\ B \end{array} \right) \,\widehat{p} =
 \left(\begin{array}{c} a \\ b \end{array} \right)\/.
\end{equation}
Here $B \in \mathbbm{R}^{N_b \times N}\/$ and $L \in \mathbbm{R}^{N_a \times N}\/$
are the discrete matrix representations for the derivatives
($\mathbf{n}\cdot \nabla\/$) in the Neumann boundary condition, and Laplacian
($\Delta\/$), respectively.  The terms $a \in \mathbbm{R}^{N_a}\/$ and
$b \in \mathbbm{R}^{N_b}\/$ are the discrete representations of the source
terms and applied boundary conditions for the Poisson equation:
\begin{eqnarray}
 a & \approx & \nabla \cdot \left(\mathbf{f} -
               (\mathbf{u} \cdot \nabla)\,\mathbf{u} \right)\/, \\
 b & \approx & \mathbf{n} \cdot \left( \mathbf{f} - \mathbf{g}_t +
               \mu\,\Delta\,\mathbf{u} - (\mathbf{u} \cdot \nabla)\,\mathbf{u}
               \right) + \nonumber \\
   &         & \lambda\,\mathbf{n} \cdot (\mathbf{u} - \mathbf{g})\/.
\end{eqnarray}
Thus the pair of equations
\begin{equation} 
    L \widehat{p} = a \quad \mbox{and} \quad B \widehat{p} = b
\end{equation}
are the discrete analog of equation (\ref{PPE_Numerical_Splitting2}).
\begin{remark} \label{rem:MatrixSolvability}
 For domains where the velocity is specified on the boundary, the pressure is only determined up to a constant.  Hence, the discrete linear system (\ref{PoissonMatrixEq1}) is singular, and discretization
 errors in $a$ and $b$ may cause solvability for equation (\ref{PoissonMatrixEq1}) to fail.  
 These solvability errors arise from numerical truncation, and are therefore ``small'' (second order). Hence, as it is
 commonly the case with numerical solutions of the Poisson equation with
 Neumann boundary conditions, we project the right hand side of (\ref{PoissonMatrixEq1}) onto 
the column space of the Poisson matrix.  For a conscise approach, on how to solve the system by constructing a square, nonsingular matrix system using the null vector $r = [1 1 \ldots 1]^T$, see \cite{Henshaw1994}. \myremarkend
\end{remark}
%
\subsection{Momentum Equation} \label{subsec:MomentumEquation}
The main numerical difficulty with implementing
(\ref{PPE_Implementation_Splitting1}) is produced by the boundary conditions.
Specifically: on a cartesian staggered grid the boundary conditions
$\nabla \cdot \mathbf{u} = 0\/$ and
$\mathbf{n}\times(\mathbf{u} - \mathbf{g}) = 0\/$ couple the ``horizontal'',
$u\/$, and ``vertical'', $\,v\/$, components of the
flow velocity ---  with the exception of the special case of a boundary aligned
with the grid, where there is no coupling. Hence, in general, the implementation
of the boundary conditions requires the solution (at every time step) of a
linear system of equations that couples all the boundary velocities (these are
defined below).  It may be possible to avoid the difficulties with the boundary
conditions discussed in this section by using a different spatial discretization.  
For example, overlapping grids which conform to the domain boundaries, such as those found
in \cite{Henshaw1994, HenshawPetersson2003}, may avoid the coupling between velocity boundary components.

The \emph{computational domain for the velocities,}
$\mathcal{C}_{\mathbf{u}}\/$,
is defined in terms of the edges in the cartesian grid with which the
velocities are associated (see figure~\ref{fig:StaggeredMesh}). To define
$\mathcal{C}_{\mathbf{u}}\/$, it is convenient to first introduce the
\emph{extended set of pressure nodes}, $\mathcal{E}_p\/$, from which
$\mathcal{C}_{\mathbf{u}}\/$ is easily constructed:
\begin{itemize}
 \item[\ref{Section:NumScheme}f]
 A pressure node in the cartesian grid belongs to $\mathcal{E}_p\/$ if and
 only if it either belongs to $\mathcal{C}_p\/$, or if it is connected (by
 an edge in the grid) to a ghost pressure node that lies on
 $\partial\/\Omega\/$.
 \item[\ref{Section:NumScheme}g]
 A velocity is in $\mathcal{C}_{\mathbf{u}}\/$ if and only if: (a) Its
 corresponding edge connects two points in $\mathcal{E}_p\/$. (b) At least
 one of the two points is in $\Omega\/$ 
 (or $\partial\/\Omega\/$).
\end{itemize}
\begin{remark} \label{rem:TheSetEp}
 Notice that $\mathcal{E}_p\/$ is exactly what $\mathcal{C}_p\/$ becomes if
 $\partial\/\Omega\/$ is modified by an ``infinitesimal'' perturbation that
 turns all the ghost pressure points on $\partial\/\Omega\/$ into inner
 pressure points.
\end{remark}
The computational domain $\mathcal{C}_{\mathbf{u}}\/$ is, in turn, sub-divided
into inner and boundary velocity edges
\begin{itemize}
 \item[\ref{Section:NumScheme}h]
 An edge in $\mathcal{C}_{\mathbf{u}}\/$ is an \emph{inner} velocity if and only
 if $\mathcal{C}_{\mathbf{u}}\/$ includes the four other edges needed to compute
 the Laplacian (either $\Delta\/u\/$ or $\Delta\/v\/$) at the edge mid-point,
 using the 5 point centered differences approximation.
 \item[\ref{Section:NumScheme}i]
 An edge in $\mathcal{C}_{\mathbf{u}}\/$ is a \emph{boundary} velocity if and
 only if it is not an inner velocity. Let $M\/$ \emph{be the number of
 boundary velocities.}
\end{itemize}
The solution of the momentum equation (\ref{PPE_Numerical_Splitting1}) is thus
performed as follows:
\begin{itemize}
 \item[\ref{Section:NumScheme}j]
 At the start of each time step the right hand side in
 (\ref{PPE_Implementation_Splitting1}) is approximated (to second order, using
 centered differences) at the inner velocity locations, which can then be
 updated to their values at the next time.
 \item[\ref{Section:NumScheme}k]
 Next, using the boundary conditions, the values of the boundary velocities at
 the next time step are constructed from the inner velocities. This
 ``extension'' process is explained below.
\end{itemize}
Let $\mathbf{y} \in \mathbbm{R}^M\/$ be the vector of boundary velocities.
Then $\mathbf{y}\/$ is determined by two sets of equations, corresponding to
the discretization on $\partial\/\Omega\/$ of $\nabla \cdot \mathbf{u}=0\/$
and $\mathbf{n}\times(\mathbf{u} - \mathbf{g}) = 0\/$. In our approach the
divergence free criteria is enforced ``point-wise'', while the condition on
the tangential velocity is imposed in a least squares sense.

\medskip \noindent
\emph{Implementation of the divergence free $\nabla \cdot \mathbf{u} = 0$
boundary condition.}
At first sight, this boundary condition appears to be the hardest to
implement, since it is a Neumann condition (essentially) prescribing the
value of the normal derivative of the normal component of the velocity,
in terms of the tangential derivative of the tangential velocity.
However, because (for the exact solution) $\nabla \cdot \mathbf{u} = 0$
everywhere, an implementation of this condition which is second order
consistent (in the classical sense of finite differences introduced by Lax
\cite{LaxFD}) is easy to obtain, as follows:

First: identify the $M_d$ pressure nodes in $\mathcal{C}_p\/$, which have at
least one boundary velocity as an adjacent edge.  These $M_d$ pressure nodes
lie either on $\partial\/\Omega\/$, or inside $\Omega\/$, and are all within
a distance $O(\Delta\/x)\/$ of $\partial\/\Omega\/$ --- definitely no further
away than $\sqrt{2}\,\Delta\/x\/$. Second: for each of these $M_d\/$ points,
use centered differences to approximate the flow divergence at the point, and
set $\nabla \cdot \mathbf{u} = 0$. This provides $M_d\/$ equations that couple
the (unknown) boundary velocities to the (known) inner velocities.  In matrix
form this can be written as
\begin{equation} \label{Boundary_Divergence} 
    D \mathbf{y} = \mathbf{s}
\end{equation}
where $D\/$ is the portion of the discrete divergence operator acting on the
unknown boundary velocities, while $\mathbf{s}\/$ is the associated flux
derived from the known inner velocities. $D\/$ is a rectangular, very sparse,
matrix whose entries are $0$ and $\pm 1$ --- note that $\Delta\/x\/$ can be
eliminated from these equations.
\begin{remark} \label{rem:BoundaryDivergence}
 Notice that $\nabla \cdot \mathbf{u} \equiv 0\/$ for the exact
 solution.  Hence, setting $\nabla \cdot \mathbf{u} = 0$ at the nodes near the
 boundary (as done above) involves no approximation. The second order nature
 of (\ref{Boundary_Divergence}) is caused by the error in computing 
 $\nabla \cdot \mathbf{u}$ --- there is no ``extrapolation'' error. \myremarkend
\end{remark}

\medskip \noindent
\emph{Implementation of the tangential velocity
$\mathbf{n}\times(\mathbf{u} - \mathbf{g}) = 0$ boundary condition.}
It is easy to see that $M_d \approx M/2\/$, since most of the $M_d\/$
pressure nodes in the selected set will connect with two boundary boundary
velocities. Thus (\ref{Boundary_Divergence}) above provides approximately
one half the number of equations needed to recover $\mathbf{y}\/$ from the
inner velocities. It would thus seem natural to seek for $M-M_d\/$
additional equations using the other boundary condition. Namely, find
$M-M_d\/$ points on $\partial\/\Omega\/$, and at each one of them write an
approximation to $\mathbf{n}\times(\mathbf{u} - \mathbf{g}) = 0$ using
nearby points in $\mathcal{C}_{\mathbf{u}}\/$. Unfortunately, this does not
work. It is very hard to do the needed approximations in a fashion that is
robust relative to the way $\Omega\/$ is embedded in the rectangular grid.
Our attempts at this simple approach almost always lead to situations where
somewhere along $\partial\/\Omega\/$ an instability was triggered.

To avoid the problem stated in the previous paragraph, we over-determine the
implementation of the tangential velocity boundary condition, and solve the
resulting system in the least squares sense. The boundary condition is replaced
by the minimization of a (discrete version) of a functional of the form
\begin{equation} \label{Least_Squares_Functional} 
  \int_{\partial \Omega} \left|\mathbf{n} \times (\mathbf{u} - \mathbf{g})
  \right|^2\,w\,\du\/A\/,
\end{equation}
where $w\/$ is some (strictly positive) weight function. This approach yields
a robust, numerically stable, approximation --- fairly insensitive to the
particular details of how $\Omega\/$ is embedded within the cartesian grid.

\begin{figure}[htb!]
\centering
\includegraphics[width = \textwidth]{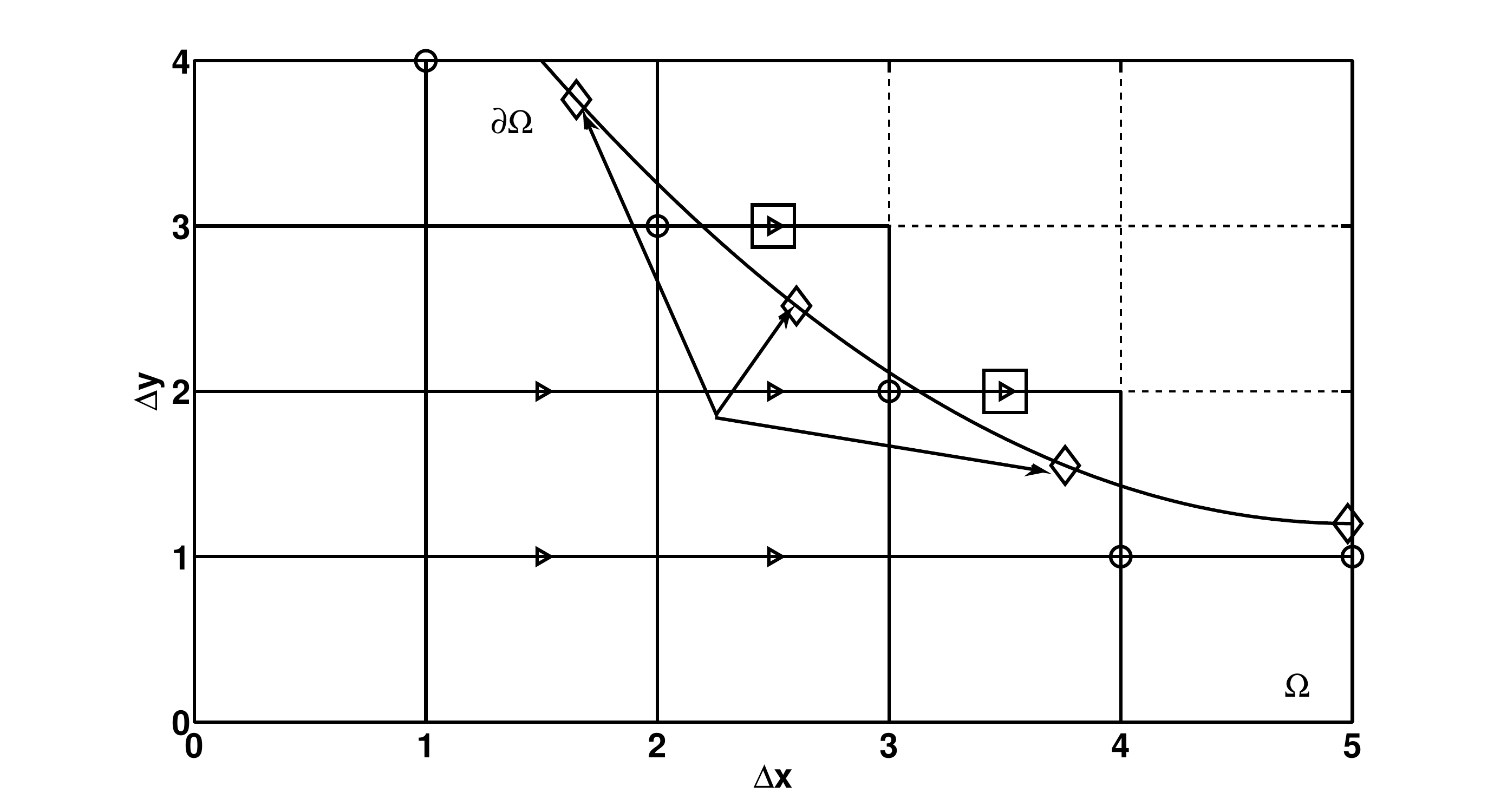}
\caption{This plot illustrates the implementation of the momentum equation
   boundary conditions. The circles ($\circ\/$) indicate the points at which
   the velocity divergence is set to zero. The six arrows indicate the
   horizontal velocity components in the patch, $\mathcal{P}_{j\,\nu}^u\/$,
   used to extrapolate $u\/$ to the three boundary points indicated by the
   diamonds ($\Diamond\/$). The squares ($\Box\/$) denote the boundary
   velocities --- which are part of the patch $\mathcal{P}_{j\,\nu}^u\/$.}
\label{fig:MomentumBC}
\end{figure}

In fact, we do not implement the minimization of a functional as in
(\ref{Least_Squares_Functional}), but a simpler process which is (essentially)
equivalent. To be specific, we start with the set $\mathcal{C}_b\/$ of points
in $\partial\/\Omega$ --- see item~\ref{Section:NumScheme}c. For each
$\mathbf{x}_{b\,j} \in \mathcal{C}_b\/$, where $1 \leq j \leq N_e\/$, we
identify a local horizontal,   $\mathcal{P}_j^u\/$, and vertical,
$\mathcal{P}_j^v\/$, velocity ``patch''. These patches --- see
figure~\ref{fig:MomentumBC} --- have the following properties
\begin{itemize}
\item[\ref{Section:NumScheme}l]
 Each patch contains both inner and boundary velocities.
\item[\ref{Section:NumScheme}m]
Each patch contains 6 velocities, in an appropriate structure, so that these
velocities can be used to extrapolate (with second order accuracy) values of
the corresponding velocity ($u\/$ or $v\/$) to nearby points along the
boundary. For example: the 5 point stencil
used to approximate the Laplacian, plus one of the four next closest points
to the center point.
\item[\ref{Section:NumScheme}n]
The union of all the patches contains all the boundary velocities.
\end{itemize}
These patches are used to (linearly) interpolate/extrapolate the velocities to
nearby points along the boundary.  For example, the horizontal velocity at point $\mathbf{x}_{b\,j}$ is extrapolated using patch $\mathcal{P}_j^u$ by
\begin{eqnarray}
    u |_{\mathbf{x}_{b\,j} } = \sum_{\mathbf{x}_k \in \mathcal{P}_j^u } \gamma_k u_k
\end{eqnarray}
where $u_k$ is the velocity at $\mathbf{x}_k$, and $\gamma_k$ are precomputed coefficients which linearly extrapolate $u$ to second order.

Each patch is used to extrapolate the
velocity to the three ``closest'' boundary points to the patch.  In this
fashion, for every $\mathbf{x}_{b\,j}\/$, $1 \leq j \leq N_e\/$, three
different approximations to the tangential velocity at the boundary follow.
These involve different (linear) combinations of velocities at the nearby
edges in $\mathcal{C}_{\mathbf{u}}\/$, including both boundary and inner
velocities. In this fashion, a set of $3\,N_e\/$ additional linear equations
--- beyond those in equation (\ref{Boundary_Divergence}) --- for the boundary
velocities follow.
\begin{remark} \label{rem:Patches}
 The idea in the process above is to ``link'' the patches used to extrapolate
 the velocity to the boundary. This is so that no ``gaps'' occur in the
 implementation of the tangential boundary conditions, as the positioning of
 $\partial\,\Omega\/$ relative to the grid changes.  We elected to use minimal coupling
 (each patch linked to its two neighbors). In principle
 one could increase the amount of over-determination of the tangential boundary
 condition. This would require using a set of points in $\partial\/\Omega\/$
 that is  ``denser'' than $\mathcal{C}_b\/$, but \emph{in our
 calculations we found that this was not needed.}

 Finally, an interesting question is: why is a process similar to this one not
 needed for the implementation of the boundary conditions for the Poisson
 equation in \S~\ref{subsec:PoissonSolver}? The obvious answer is that the
 Poisson equation itself couples everything, so that no extra coupling needs
 to be added. \myremarkend
\end{remark}

Putting everything together, an overdetermined system of equations for the
boundary velocity vector $\mathbf{y}\/$ follows, which can be written in the
form
\begin{eqnarray} 
   D \mathbf{y} & = & \mathbf{s}\/,  \label{BC_Div}\\
   E \mathbf{y} & = & \mathbf{t}\/.  \label{BC_TanVel}
\end{eqnarray}
Here $E \in \mathbbm{R}^{3\,N_e \times M}$, and the second equation is to be
solved in the least squares sense, subject to the constraint imposed by the
first equation.
One way to do this is as follows: write
\begin{equation} \label{Divergence_Constraint}
    \mathbf{y} = \mathbf{y}_p + P \mathbf{c}
\end{equation}
where $P\/$ is a matrix whose columns are a basis for the kernel of $D\/$ ---
hence $D\,P = 0\/$, $\;\mathbf{y}_p\/$ is a particular solution of
(\ref{BC_Div}), and $\mathbf{c}\/$ is a constant vector parameterizing the
space of (numerical) divergence free boundary velocities.  Substituting the
ansatz (\ref{Divergence_Constraint}) into (\ref{BC_TanVel}) yields
\begin{equation} \label{Reduced_Least_Squares}
    (E\,P)\,\mathbf{c} = \mathbf{t} - E\,\mathbf{y}_p\/,
\end{equation}
which is a constraint-free least squares problem for the vector $\mathbf{c}\/$.
Thus we can write
\begin{equation} \label{SlnRedLSq}
    \mathbf{c} = (E\,P)^{+}\,\left(\mathbf{t} - E\mathbf{y}_p\right)\/,
\end{equation}
which, together with (\ref{Divergence_Constraint}), gives the solution
$\mathbf{y}\/$.

\begin{remark} \label{rem:NoZeroSV}
 In (\ref{SlnRedLSq}), $(E\,P)^{+}\/$ is the pseudo-inverse, defined by the
 singular value decomposition of the matrix $F = E\,P$. In fact, since
 $F^T\,F\/$ should be invertible (see below),  $F^+ = (F^T\,F)^{-1}\,F^T\/$.

 One of the objectives of the above construction is to ensure
 that the boundary velocities are completely determined by the boundary
 conditions (and the inner velocities). If the problem in (\ref{BC_TanVel}) 
 for $\mathbf{y}\/$, with the constraint in (\ref{BC_Div}), were to have
 more than one solution that minimizes the $L^2\/$ norm of 
 $\mathbf{t}-E\,\mathbf{y}\/$, then this would be a sure
 sign that the boundary condition implementation is flawed, and there is
 missing information (i.e.: more equations are needed, see
 remark~\ref{rem:Patches}).

 Of course, the requirement in the prior paragraph is equivalent to the
 statement that $F^T\,F\/$ is invertible.  Finally, notice that in a 
 domain with a fixed boundary, one may preprocess the matrices $E\/$, $D\/$ 
 and $P\/$.
 \myremarkend
\end{remark}
\begin{remark} \label{rem:consistent}
 It should be clear that the scheme developed in this section, is second order
 consistent\footnote{For the system of PDE in
    (\ref{PPE_Numerical_Splitting1}--\ref{PPE_Numerical_Splitting2}).}
 (all the way up to the boundary) in the classical sense of finite
 differences introduced by Lax \cite{LaxFD}. \myremarkend
\end{remark}
%
%
\subsection{Comparison with the Projection Method}
The scheme proposed here --- which involves the implementation of the equations
in (\ref{PPE_Implementation_Splitting1}--\ref{PPE_Implementation_Splitting2}),
appears to be quite similar to fractional step methods, such as
Chorin's \cite{Chorin} original projection method.  Specifically: advancing one
time step requires both the evolution of a diffusion equation for the flow
velocity, and the inversion of a Poisson equation for the pressure --- which is
the same as in projection methods. But there are important differences, mainly
related to the implementation of the boundary conditions, and of the
incompressibility condition. Below we highlight the similarities and
differences between the pressure Poisson approach proposed here and, for
simplicity, the projection method in its original formulation~\cite{Chorin}.
A (fairly recent) thorough review of projection methods, exploring 
the improvements to the approach since~\cite{Chorin}, as well as their
drawbacks, can be found in~\cite{Guermond}.

First, the starting point for the method here is the discretization of a
reformulation of the equations:
(\ref{PPE_Numerical_Splitting1}--\ref{PPE_Numerical_Splitting2}), not the
``original'' Navier-Stokes equations (\ref{NS_Momentum}--\ref{NS_Mass}). In
this equivalent set: (i) there is a natural way to recover the
pressure from the flow field at any given time. (ii) The time evolution
automatically enforces incompressibility. As a consequence, there is
(in-principle) no limitation to the order in time to which the reformulated
equations can be numerically discretized. In contrast, the projection method
is equivalent to an approximate $L\,U\/$ matrix
factorization~\cite{Guermond, Perot, Taira} of the discrete differential operators
coupling $\mathbf{u}\/$ and $p\/$. This approximate factorization yields a
splitting error in time, which is very hard to circumvent in order to achieve
higher accuracy. 

Second, the pressure Poisson formulation used here ensures both that, at
every point in the time integration: (i) the normal and tangential velocity
boundary conditions are accurately satisfied. (ii) The zero divergence
condition is accurately satisfied. On the other hand, in the original
projection method, the step advancing the flow velocity forward in time
enforces the velocity boundary conditions, but cannot guarantee that
incompressibility is maintained. In the other step, the pressure
is used to recover incompressibility --- by projecting the flow velocity field
onto the space of divergence free fields. Unfortunately, while removing the
divergence, this second step does not necessarily preserve the correct
velocity boundary values~\cite{Guermond}.

Finally the modified equation (\ref{PPE_Implementation_Splitting1}) implemented
here differs from the projection's method velocity forward step primarily by
the boundary conditions imposed.
Specifically, the projection method imposes Dirichlet boundary conditions for
each component of the velocity field when propagating the flow velocity.
Dirichlet boundary conditions have the obvious advantage of decoupling the
velocity components. Numerically this means that an implicit treatment of the
stiff viscosity operator $\partial_t - \mu\,\Delta\/$ is straightforward. By
contrast, on a regular grid, the boundary conditions in
(\ref{PPE_Implementation_Splitting1}) couple the boundary velocities. Although
this coupling does not pose a serious difficulty for explicit schemes, it adds
a (potentially serious) difficulty to the implementation of any scheme that
treats the viscosity operator in the equations implicitly.


\section{Implementation} \label{Section:Implementation}
The objective of this section is to study the convergence and accuracy of the
schemes proposed in \S~\ref{Section:SemiImplicitStability} and \S~\ref{Section:NumScheme}, and to test the methods for a simple physical problem.  In particular: to verify the theoretical prediction that the temporal splitting in \S~\ref{Section:SemiImplicitStability} is second order in time, while the spatial discretization in \S~\ref{Section:NumScheme} is second order in space, all the way up to the boundary. 

To this end, here we present the results from numerical
computations for
(i) Externally forced flow on a square domain in \S~\ref{subsec:FlowSqD}, 
(ii) Computation of flow in a lid driven cavity in \S~\ref{subsec:drivencavity}, and 
(iii) Externally forced flow on an irregular domain in \S~\ref{subsec:FlowIrrD}.
In examples (i) and (iii), the external forcing is selected so that an exact
(analytic) solution to the equations can be produced, following the same procedure as in \cite{Guermond}.\footnote{That is,
   proceed backwards: first write an incompressible flow, and then compute
   the forcing, and boundary conditions, it corresponds to.}  In what follows the numerical errors are measured using the discrete
$L^{\infty}$ grid norm defined by
\begin{equation} \label{eqn:LinfGridNorm}
   \|g\|_{\infty} = \mathrm{max} \{ |g_{i\,j}| \}\/.
\end{equation}
Here the maximum is over all the indexes $i\/$ and $j\/$ --- corresponding to
the appropriate grid coordinates $(x_i\/,\, y_j)\/$ in the computational domain
--- for the field $g\/$ in the staggered grid. Namely: nodes for the pressure
$p\/$, mid-points of the horizontal edges for $u\/$, and mid-points of the
vertical edges for $v\/$. Note that neither ghost pressure points, nor boundary
velocities, are included within the index set in (\ref{eqn:LinfGridNorm}) ---
we consider these as auxiliary numerical variables, introduced for the purpose
of implementing the boundary conditions, but not part of the actual solution.
Finally, for the exact (analytic) solution, grid values are
obtained by evaluation at the points $(x_i\/,\, y_j)\/$. For example, to
compute the error in the pressure, first define the discrete field
$e_{ij} = \hat{p}_{i\,j} - p(x_i\/,\,y_j)\/$ --- where $\hat{p}\/$ is the
numerical pressure field and $p\/$ is the exact continuous field,
and then compute the norm above for $e\/$.

As should be clear from the prior sections,
the main issue we aim to address in this paper, is that of how to effectively
implement the incompressibility condition and the boundary conditions for the
pressure, avoiding the difficulties that projection and fractional step
methods have. These are problems that are not related to the nonlinearities
in the Navier-Stokes equations, and occur even in the absence of the advective term $(\mathbf{u}\cdot \nabla)\mathbf{u}$.  Hence, for simplicity, the calculations presented in
subsections \S~\ref{subsec:FlowSqD} and \S~\ref{subsec:FlowIrrD} do not contain the advective term.  This is the same practice employed in \cite{Guermond}.  Subsection  \S~\ref{subsec:drivencavity}, however, investigates a physical flow and therefore includes the advective term.
%
%
\subsection{Flow on a square domain} \label{subsec:FlowSqD}
\begin{figure}[htb!]
\centering
\includegraphics[width=0.7\textwidth]{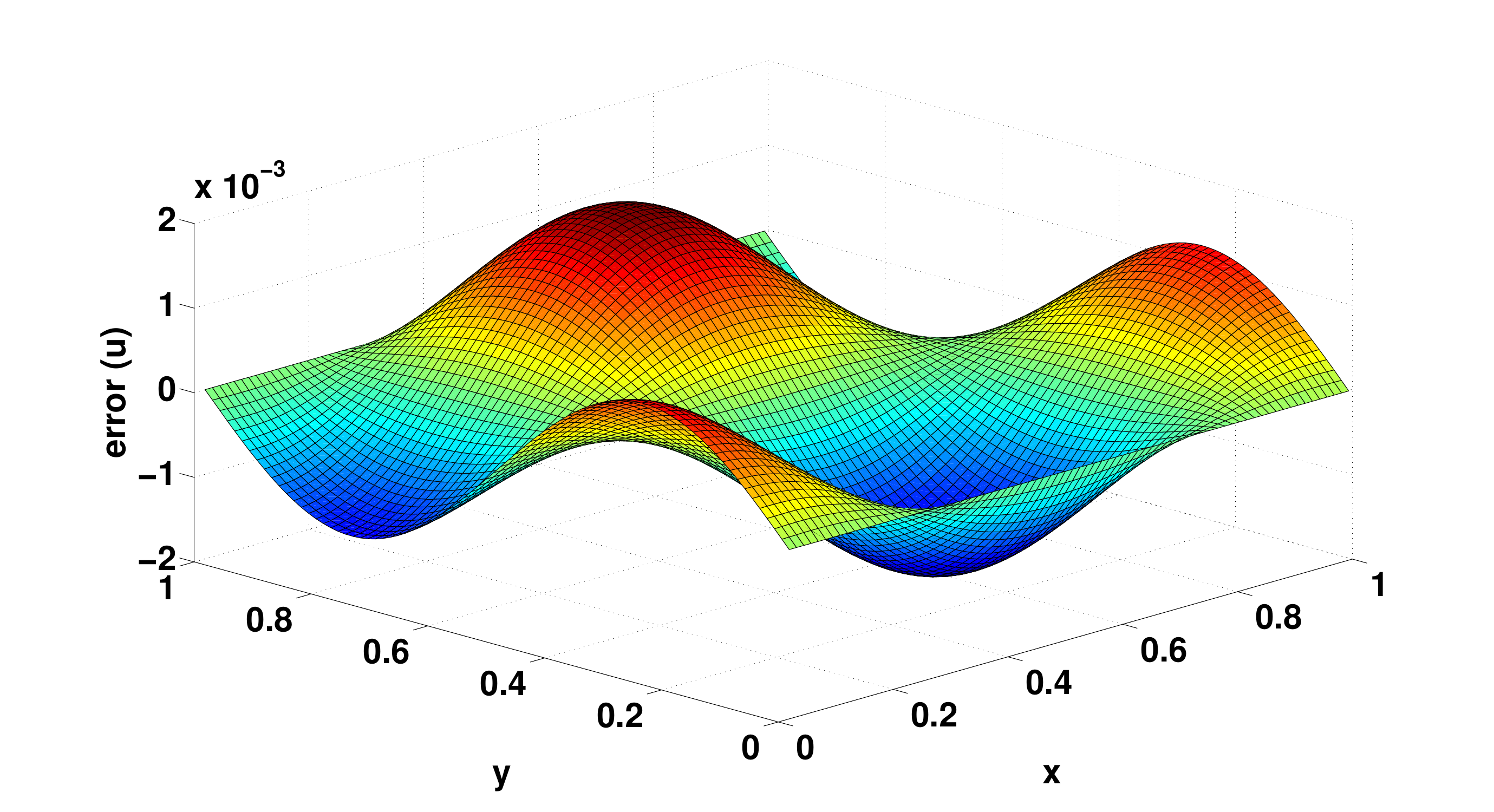}
\includegraphics[width=0.70\textwidth]{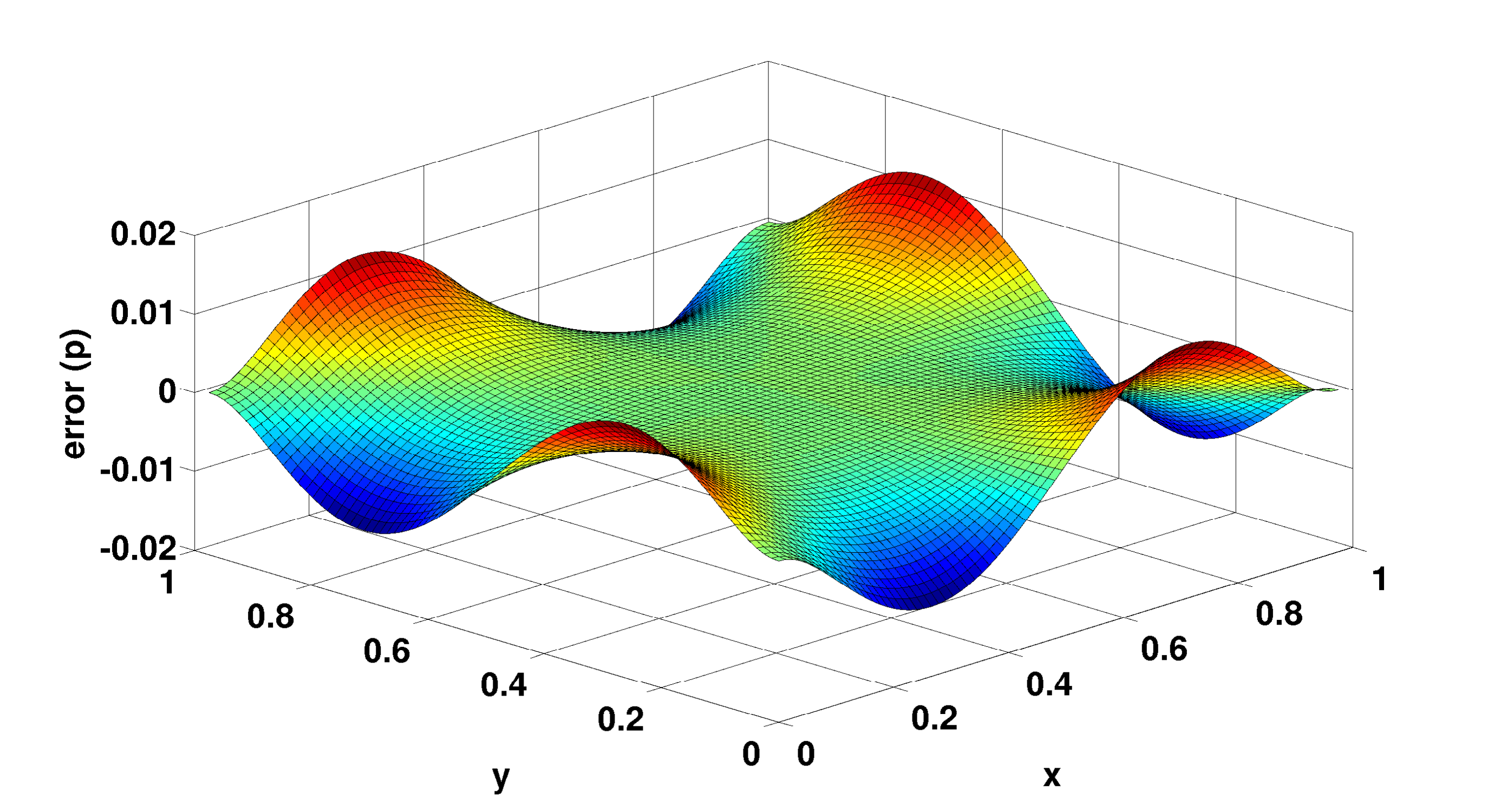}
\caption{Error fields for the second order semi-implicit numerical scheme. 
   The plots are at time $t = 4\,\pi\/$, for an $80\times 80\/$ grid, with
   $\Delta\,t=0.2\,(\Delta\,x)\/$ and $\lambda = 100$.
   The horizontal velocity $u\/$ (top) and the pressure $p\/$ (bottom)
   error fields are shown. The error is uniform in size across
   the domain. There are neither numerical boundary layers, nor
   numerical corner layers.}
\label{fig:SquareErrorFields}
\end{figure}

Here we test the semi-implicit scheme outlined in \S~\ref{Section:SemiImplicitStability}.  
Specifically, we present the results of solving the linear Navier-Stokes
equations on the unit square $0 \leq x\/,\, y \leq 1\/$, with no-slip and no
flux boundary conditions (i.e.: $u = v = 0\/$ on $\partial\,\Omega\/$), and
viscosity $\mu = 1\/$. We select the forcing function
$\mathbf{f} = \mathbf{u}_t + \nabla\,p - \mu\,\Delta\,\mathbf{u}\/$
to yield the following (incompressible) velocity and pressure fields:
\begin{eqnarray}
 u (x\/,\,y\/,\,t) & = & \phantom{-}\pi\,\cos(t)\,\sin(2\,\pi\,y)\,
    \sin^2(\pi\,x)\/, \label{UField_SquareForcing}\\
 v (x\/,\,y\/,\,t) & = & -\pi\,\cos(t)\,\sin(2\,\pi\,x)\,\sin^2(\pi\,y)\/,
    \label{VField_SquareForcing}\\
 p (x\/,\,y\/,\,t) & = & -\cos(t)\,\cos(\pi\,x)\,\sin(\pi\,y)\/.
    \label{pField_SquareForcing}
\end{eqnarray}
We numerically evolve the equations with varying grid 
and time steps. We then compare the numerical results with the exact fields in
(\ref{UField_SquareForcing}--\ref{pField_SquareForcing}).  For instance,
figure~\ref{fig:SquareErrorFields} shows the error fields for the velocity and
pressure at time $t = 4\,\pi\/$, for an $80 \times 80\/$
grid, time step $\Delta t = 0.2 \, \Delta x$ and feedback $\lambda = 30$. 
Note that this test is essentially the same as the one used in
  \cite{Guermond}\footnote{See \S~3.7.1, equation (3.30).},
  to test projection and fractional step methods. Hence the results
  in this subsection can be used as a comparison basis for the approach
  in this paper versus projection methods. Unlike the common situation with projection methods, in our case there
are neither numerical boundary layers, nor numerical corner layers.
\begin{figure}[htb!]
\centering
\includegraphics[width=0.85\textwidth]{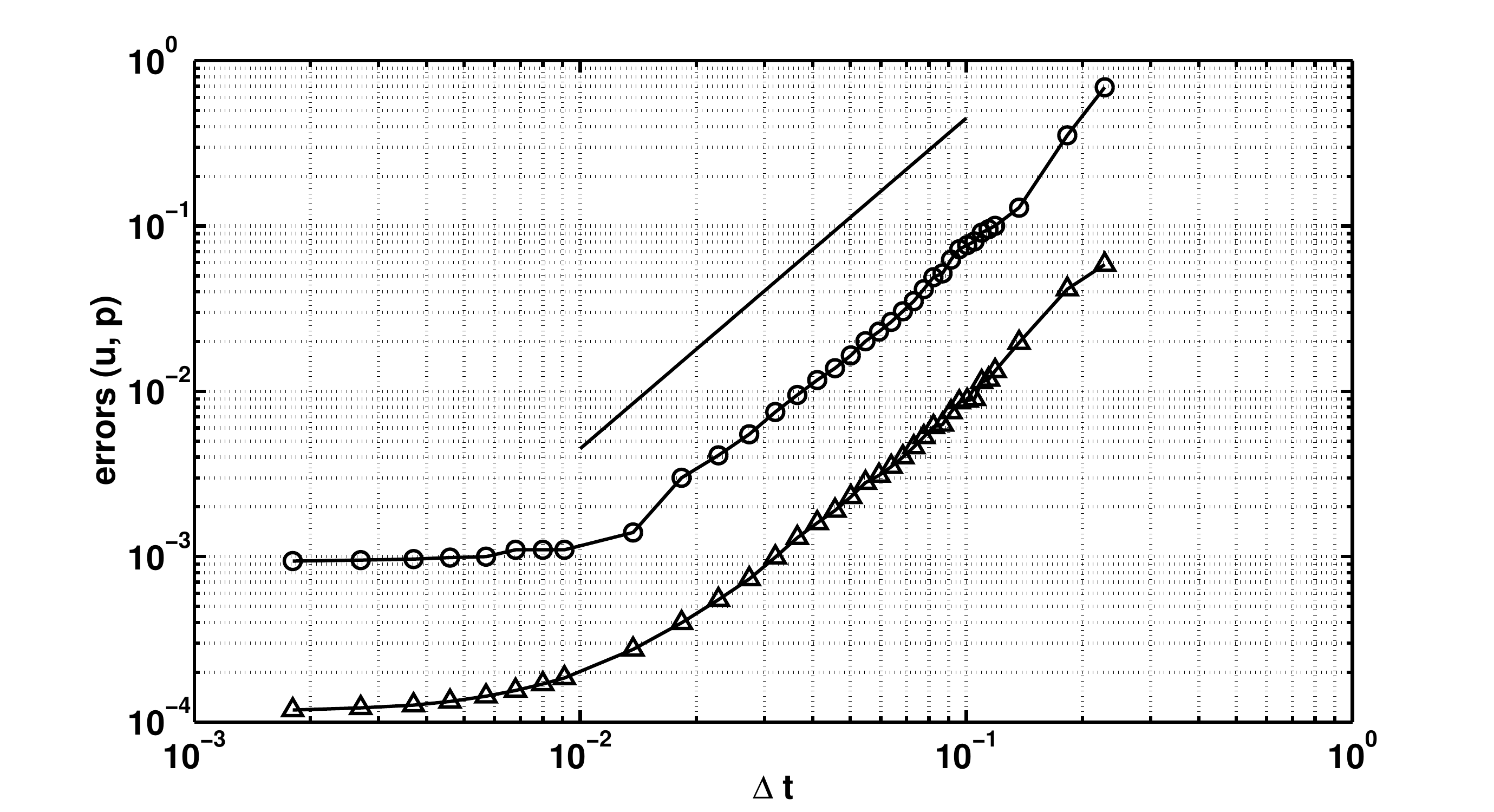}
\caption{Convergence plot for semi-implicit numerical scheme.  The plot corresponds to a fixed $220\times 220$ grid and varying time steps ($\lambda = 30$).  The straight line corresponds to second order convergence, while the saturation occurs around $\Delta t = \Delta x$ when the temporal errors become the same order as the spatial errors.  The circles ($\;\circ\/$) and squares ($\;\square\/$) correspond to the $L^\infty$ norm of the pressure and horizontal velocity respectively.}\label{fig:ImplicitConvergence}
\end{figure}

To verify the scheme is second order in time, we fix a $220 \times 220$ grid and evolve the solution from $t = 0$ to $t = 1$ for varying $\Delta t$.  Figure (\ref{fig:ImplicitConvergence}) shows second order convergence, while the saturation occurs when the spatial and temporal errors become comparable at $\Delta t \sim \Delta x$.

\begin{figure}[htb!]
\centering
\includegraphics[width=0.85\textwidth]{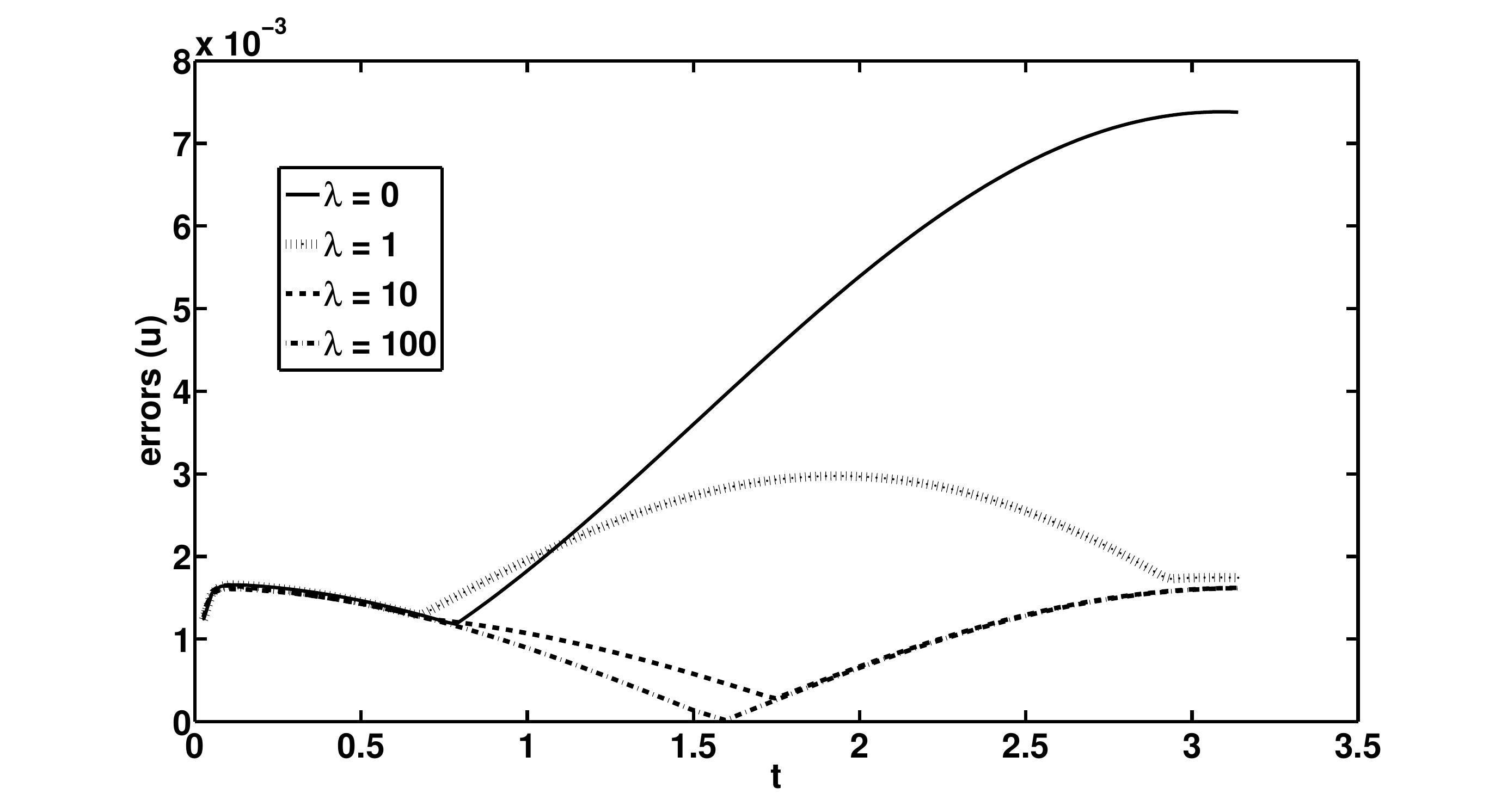}
\caption{Time evolution of the $L^\infty\/$ error in $u\/$ on an $80 \times 80\/$ mesh. 
   The flow is computed using the second order semi-implicit scheme with different
   values of the parameter $\lambda\/$. For $\lambda = 0\/$ the
   error grows in time. Values larger than $\lambda \sim 10\/$
   control the error.}
\label{fig:LambdaTest_U_Errors}
\end{figure}

We now illustrate the role of the
feedback parameter $\lambda\/$. For this purpose we present here the results of
two tests.  In the first test, we evolved the numerical solutions (with a fixed
$\Delta\,x\/$) for different values of $\lambda\/$. \footnote{With an $80 \times 80$ grid, viscosity $\mu = 1\/$ and a a time step
   $\Delta\,t = 0.2\,(\Delta\,x)\/$.}
Figure~\ref{fig:LambdaTest_U_Errors} illustrates the results of this first
test, with plots of the $L^{\infty}\/$ error in the horizontal velocity $u\/$,
as a function of time, for various values of $\lambda\/$. The errors in the
pressure $p\/$, and in the vertical velocity $v\/$, exhibit the same qualitative
behavior. In the absence of feedback --- i.e. $\lambda = 0\/$, the errors grow
steadily in time, as can be expected from equation
(\ref{Boundary_Evolution_Drift}) when influenced by the presence of numerical
noise. Even though this effect does not constitute an instability (in the sense
of exponential growth), there appears to be no bound to the growth. Thus, after
a sufficiently long time, the errors can become substantial. The figure shows
also that moderate values of $\lambda\/$ --- i.e. $\lambda \sim 10$ or
bigger, are enough to control the errors.

Physically, the errors that occur when $\lambda\/$ is too small correspond to
fluid leakage through the domain walls. However, these errors cause little net mass
loss or gain, since $\nabla \cdot \mathbf{u} = 0\/$ applies --- actually,
$\nabla \cdot \mathbf{u} = O\left((\Delta\,x)^2\right)\/$, as we show later:
see figure~\ref{fig:Divergence}. Any positive flow across some
part of the boundary must be compensated by a negative flow elsewhere.

In the second test, we introduce artificial random errors to the normal
velocity along the boundary. The purpose of this second test is to study, under
a controlled situation, the error sensitivity of the normal velocity boundary
condition implementation --- the purpose for which the parameter $\lambda\/$
was introduced in \S~\ref{Section:NumSplitting}. Specifically, the main concern
are the errors that are introduced by the implementation of the Neumann
boundary condition for the pressure on non-conforming boundaries (as happens
for the example in \S~\ref{subsec:FlowIrrD}). 
\begin{figure}[htb!]
\centering
\includegraphics[width=0.85\textwidth]{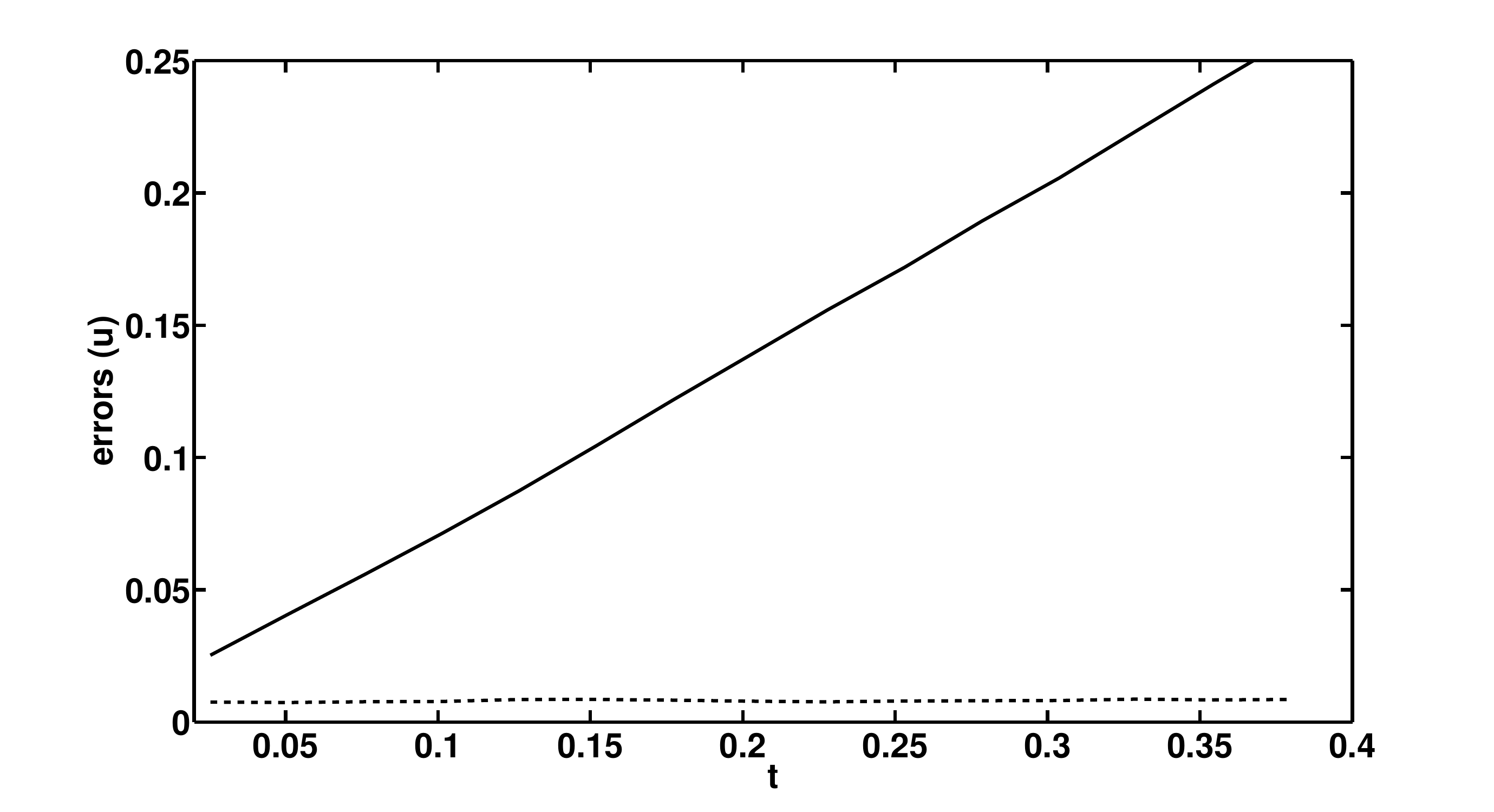}
\caption{Time evolution of the $L^\infty\/$ error in $u\/$ on an $80 \times 80$ mesh.
    Flow is computed using a second order semi-implicit scheme and random forcing on
   the boundary. Without feedback
   ($\lambda = 0\/$, solid line) the error grows steadily. With feedback
   ($\lambda = 100\/$, dashed line) the error growth saturates at a small
   value.}
\label{fig:FeedbackControlErrors}
\end{figure}

To perform the test, we introduce independent random errors at each point on
the domain boundary, at the start of each Euler time step. To be precise, the
boundary condition for equation~(\ref{PPE_Implementation_Splitting2}) is
taken as
\begin{equation} \label{perturbedPPE_BC}
   \mathbf{n} \cdot \nabla\/p^n = \mathbf{n} \cdot \left( \mathbf{f}^n
   + \Delta\/\mathbf{u}^n + \lambda\,\mathbf{u}^n \right) + r\/,
\end{equation}
where, at each boundary point $r\/$ is sampled randomly from the interval
$[0, 1]\/$. This represents an $O(1)\/$ random perturbation to the normal
component of the desired velocity, $\mathbf{g} = 0\/$, at the
boundary.\footnote{Recall that in this section we take $\mu = 1\/$, and neglect
   the nonlinear terms.}
Hence, it is a very demanding test of how insensitive the boundary condition
implementation is to errors. 

Figure~\ref{fig:FeedbackControlErrors} shows the error in
the horizontal velocity $u\/$ as a function of time, both for the
feedback controlled solution with $\lambda = 100\/$, and the undamped
solution with $\lambda = 0\/$.  The errors saturate in the feedback solution,
but they grow steadily for the undamped solution.  Without damping, the drift
errors contribute sizable effects to the solution after a short time period.

\begin{figure}[htb!]
\centering
\includegraphics[width=0.85\textwidth]{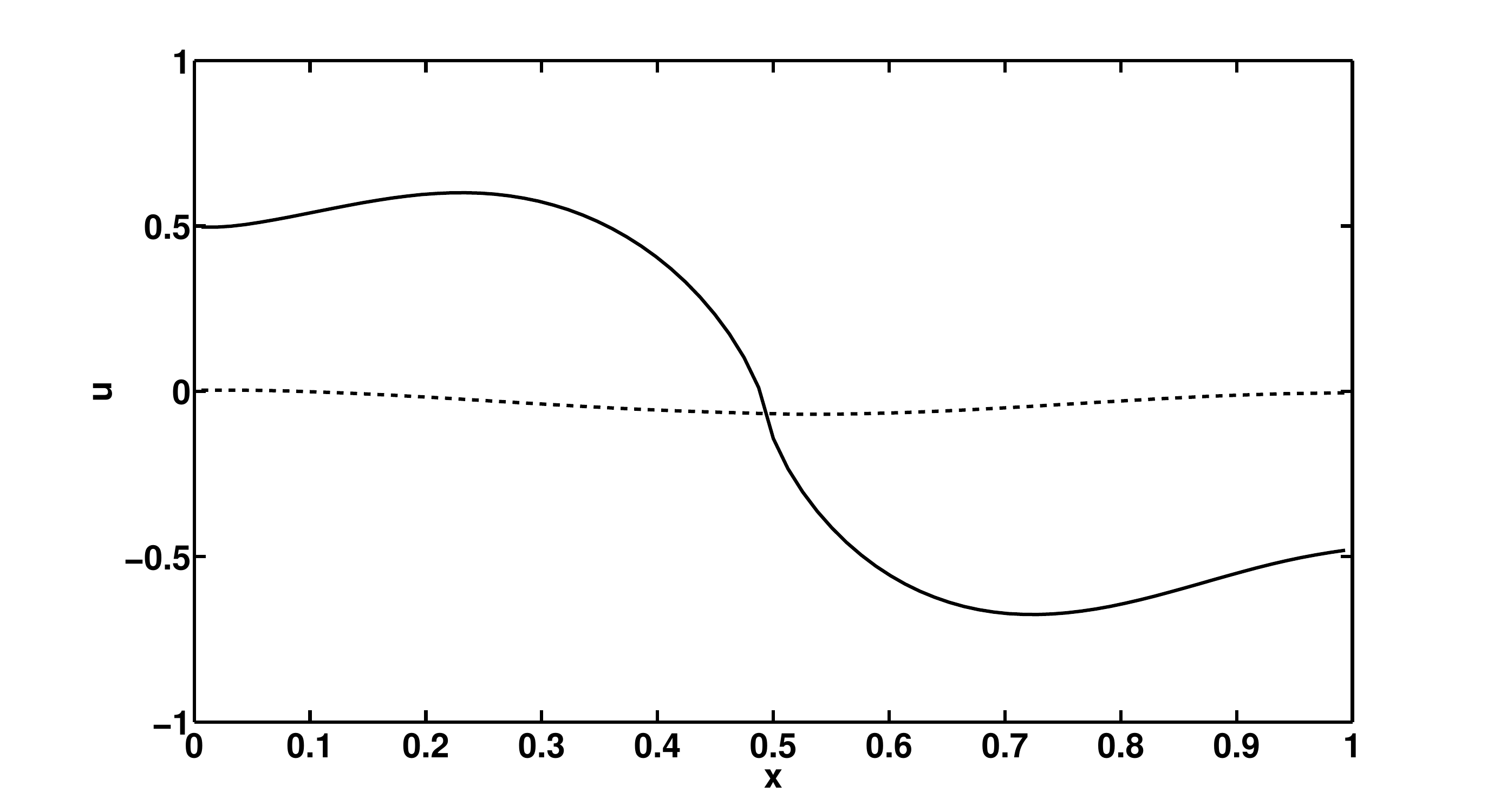}
\caption{Velocity cross-section $u(x\/,\,y = 0.4937)\/$,
   at time $t = 1\/$, for the second order semi-implicit scheme, 
   on an $80 \times 80$ mesh, with a random forcing on
   the boundary. The dashed line ($\lambda = 100\/$) reproduces the
   correct field with zero flux at the boundaries. The solid line
   ($\lambda = 0\/$) shows a non-zero flow through the boundary.}
\label{fig:FeedbackControl_U}
\end{figure}

Finally, figure~\ref{fig:FeedbackControl_U} illustrates the flow leakage
produced by the random errors. This figure shows a cross section of the
horizontal velocity $u\/$ after time $t = 1$, for two values of the
control parameter: $\lambda = 0\/$ and $\lambda = 100\/$. The first value
shows a significant flow through the boundary, while the second does not.
This plot also illustrates the point made earlier: there is no significant
mass loss (gain) even when $\lambda = 0\/$, with positive flows compensated
by negative flows elsewhere.
%
%

\subsection{Driven Cavity} \label{subsec:drivencavity}
\begin{figure}[htb!] 
\centering
\includegraphics[width=0.90\textwidth]{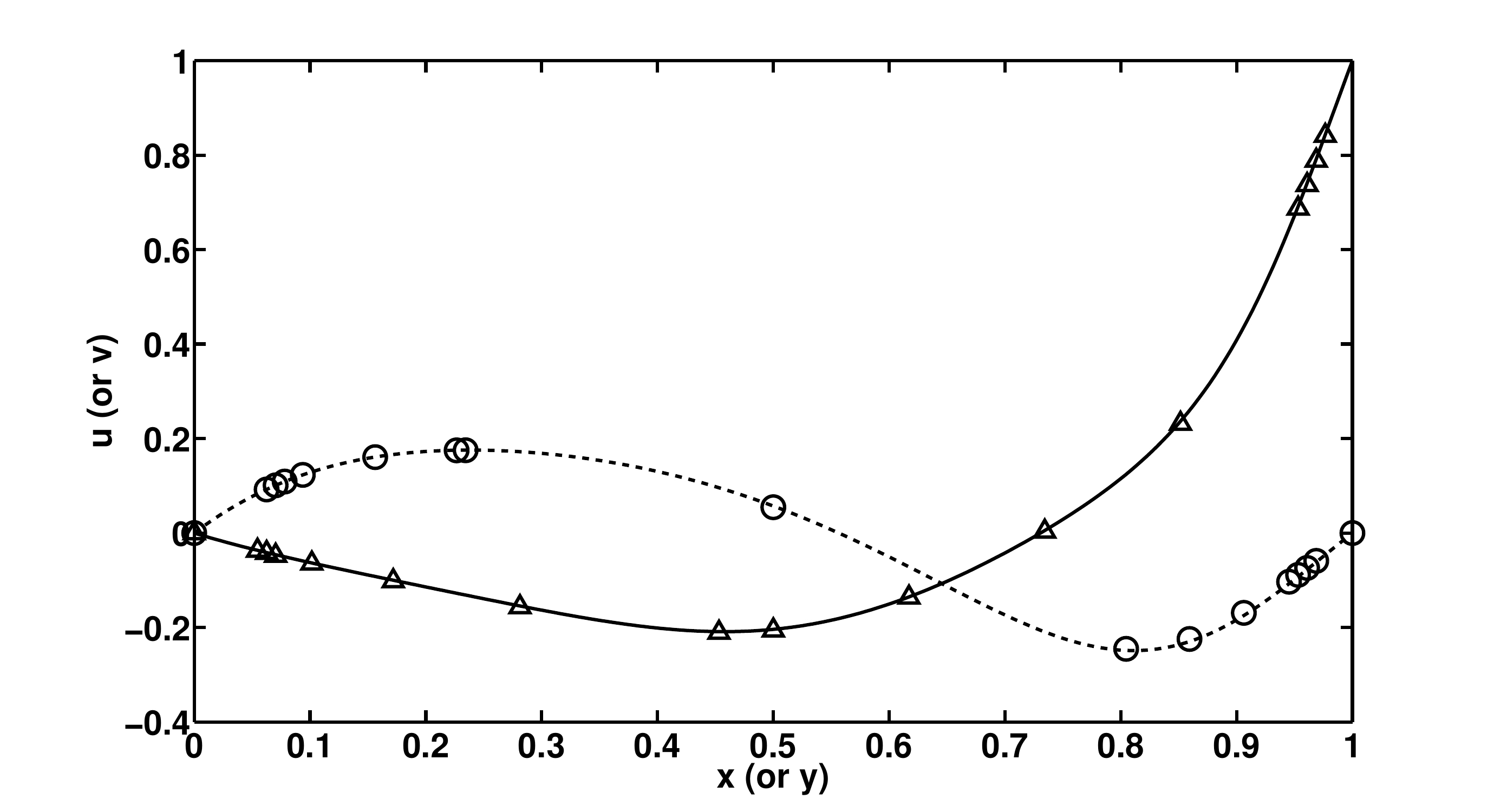}
\caption{Midpoint steady state velocities (after $t = 20$) for $Re = 100$ and a $220\times 220$ grid.  For $\mathbf{u} = (u, v)$, the 
solid line shows $u(0.5, y)$, the dashed line shows $v(x, 0.5)$, while the circles and triangles 
correspond to data from \cite{GhiaGhiaShin1982}.}
\label{DrivenCavityMidVelocitiesRe100}
\end{figure}

\begin{figure}[htb!] 
\centering
\includegraphics[width=0.90\textwidth]{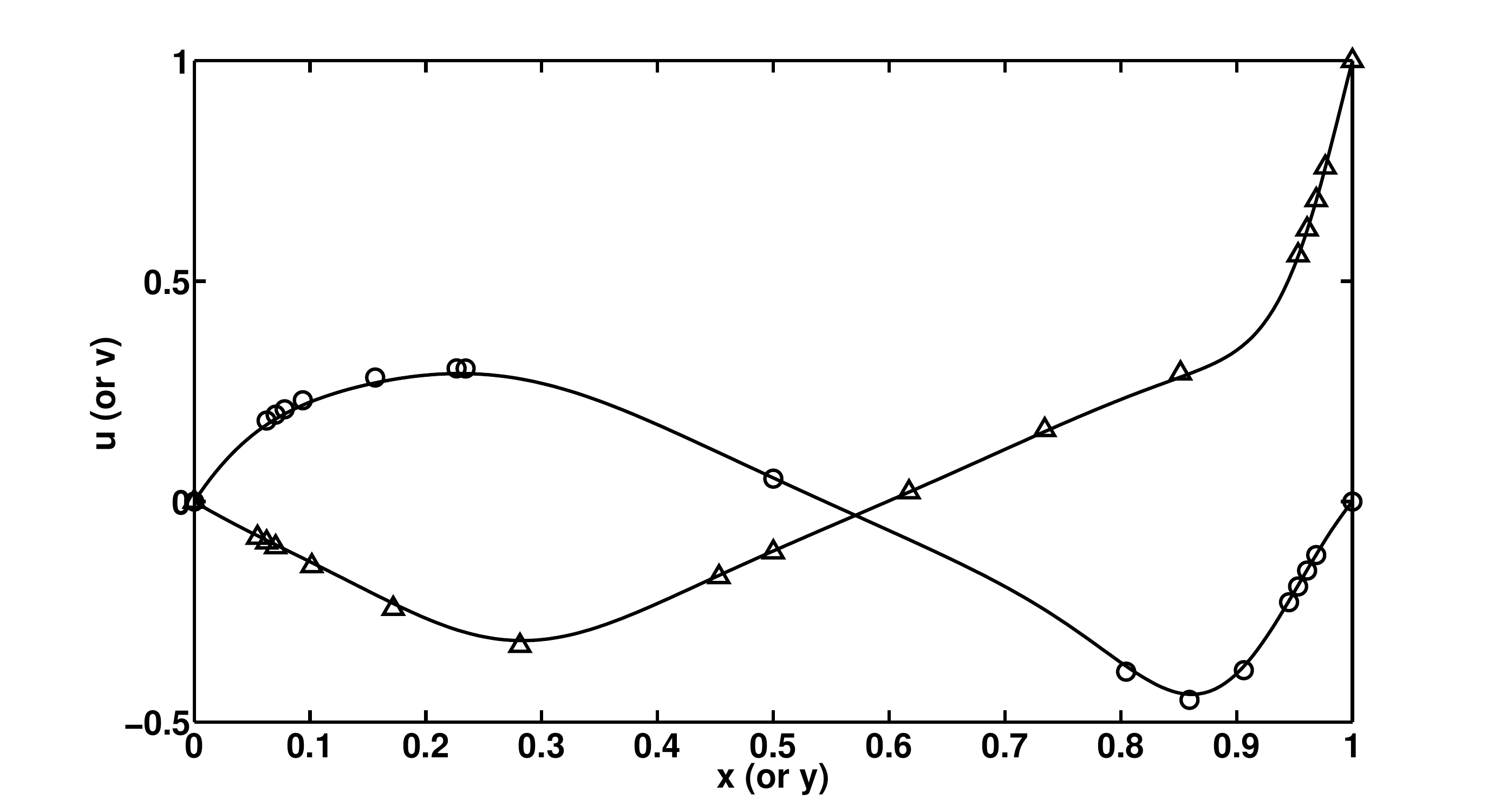}
\caption{Midpoint steady state velocities (after $t = 30$) for $Re = 400$ and a $220\times 220$ grid.  For $\mathbf{u} = (u, v)$, the 
solid line shows $u(0.5, y)$, the dashed line shows $v(x, 0.5)$, while the circles and triangles 
correspond to data from \cite{GhiaGhiaShin1982}.}
\label{DrivenCavityMidVelocitiesRe400}
\end{figure}

\begin{figure}[htb!] 
\centering
\includegraphics[width=\textwidth]{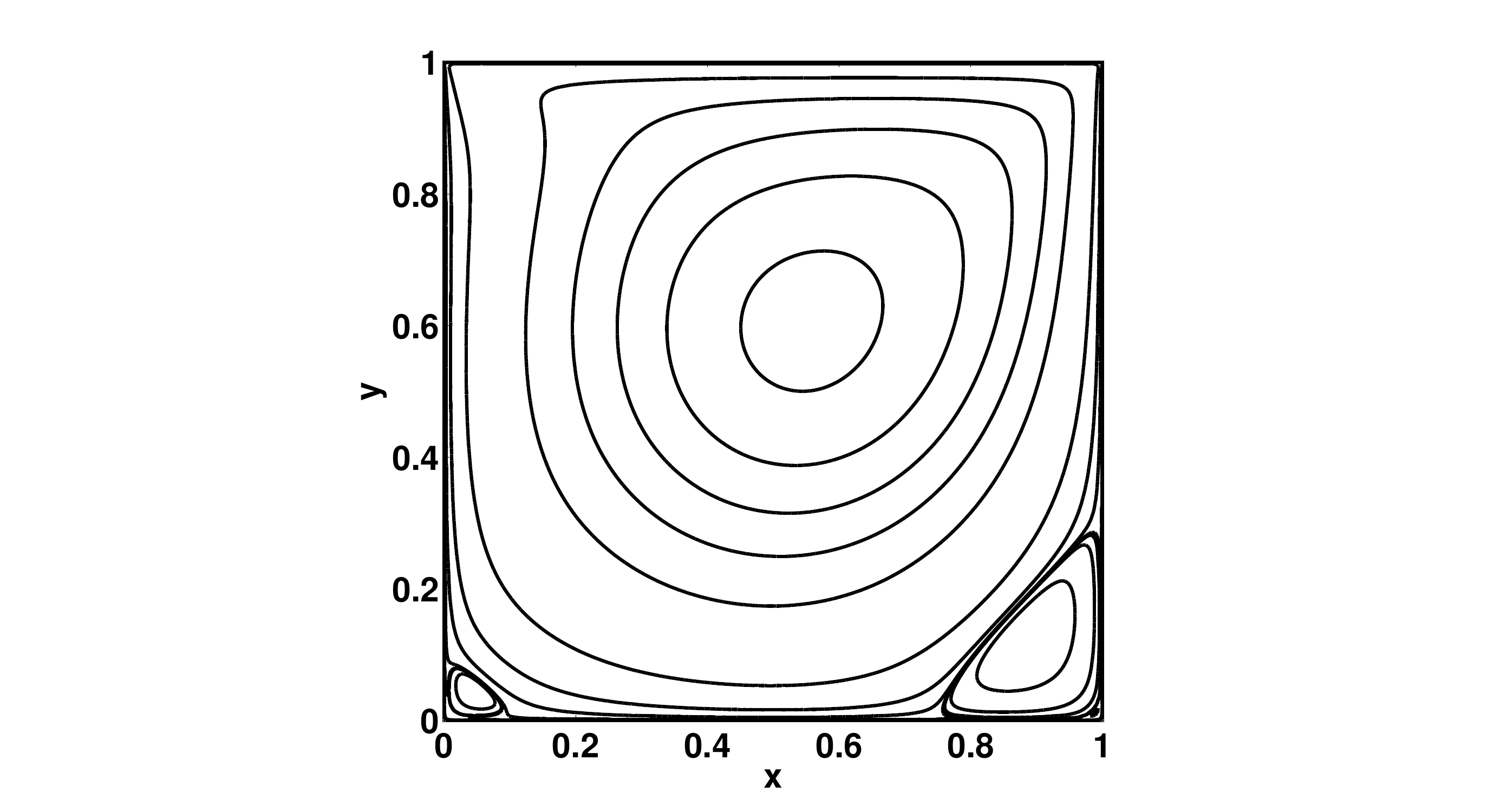}
\caption{Plots shows steady state streamfunction contours for the driven cavity with $Re = 400$.}
\label{fig:DrivenCavityStreamfunction}
\end{figure}

In this subsection, we solve for the fluid velocity field in a driven 
cavity.  We use the second order semi-implicit scheme (\ref{StokesCN}) and the regular staggered 
grid described in \S~\ref{Section:NumScheme}.  In addition, this test also demonstrates the application of a non-vanishing 
boundary condition where $\mathbf{n}\times \mathbf{u} \neq 0$.  Here we include the nonlinear terms and treat them explicitly in the numerical 
scheme.  To discretize them in space, we use a second order upwind (one sided) finite difference approximation.  In this test, we take the 
velocity of the moving wall $U_{wall} = 1$ and
consider different Reynolds numbers.  Starting with rest initial conditions, we evolve the fields until they reach a steady state flow.  
Figures \ref{DrivenCavityMidVelocitiesRe100} and \ref{DrivenCavityMidVelocitiesRe400} show the midpoint steady state velocities 
for $\mathbf{u}$ at Reynolds numbers $Re = 100$ and $Re = 400$ respectively.  The plots are done with a $220 \times 220$ grid and also
compare the computation to the standard benchmark of \cite{GhiaGhiaShin1982}.  In addition figure \ref{fig:DrivenCavityStreamfunction}
shows the streamfunction contours for a Reynolds number of $Re = 400$.


\subsection{Flow on an Irregular Domain} \label{subsec:FlowIrrD}
In this section we implement the numerical scheme described in
\S~\ref{Section:NumScheme}, for an example with an externally forced velocity
field on an irregular shaped domain.  We perform the same test, and take the same forcing from 
\S~\ref{subsec:FlowSqD} using equations (\ref{UField_SquareForcing}--\ref{VField_SquareForcing}) however, instead of the unit square, the selected
domain $\Omega\/$ is the $1 \times 1\/$ square $0 < x\/,\,y < 1\/$, with
the disk of radius $1/4\/$ centered at $(1/2\/,\, 1/2)\/$ removed. In addition,
we use periodic boundary conditions in the $y$ direction,
and impose a nonzero flux ($\mathbf{u} = \mathbf{g} \neq 0\/$) on the rest of
the boundary.

In this example we set the feedback parameter $\lambda\/$ to $\lambda=100\/$,
considerably larger than the minimum needed to get good behavior in the
calculations in \S~\ref{subsec:FlowSqD}. This is to be expected from the
results of the last test in \S~\ref{subsec:FlowSqD} --- see
figures~\ref{fig:FeedbackControlErrors} and \ref{fig:FeedbackControl_U}.
Because of the curved inner boundary in the current example, errors in the
implementation of the Neumann boundary condition for the pressure occur, which
trigger a drift in the normal boundary condition for the velocity --- unless
a large enough $\lambda\/$ is used. Thus, for example, values where
$\lambda = O(10)\/$ did not allow the feedback to quickly track fluctuations
in the normal velocity component of the boundary conditions.

\goodbreak
\begin{center} \emph{Second order convergence.} \end{center}

To illustrate the second order convergence of the scheme, below we present the
results of evolving the velocity and pressure fields (for the problem described
above), for varying grid sizes $\Delta\,x\/$, up to the fixed time
$t = 0.0657\/$, using $\lambda = 100\/$ and $\Delta\,t = 0.2\,(\Delta\,x)^2\/$.
Note that the selected final time appears small, but it is large enough to
require a number of time steps in the range $O(10^3)\/$ to $O(10^4)\/$ --- for
the grid sizes we used. This is large enough to provide a reliable measure of
the order of the scheme.

\begin{figure}[htb!]
\centering
\includegraphics[width=0.85\textwidth]{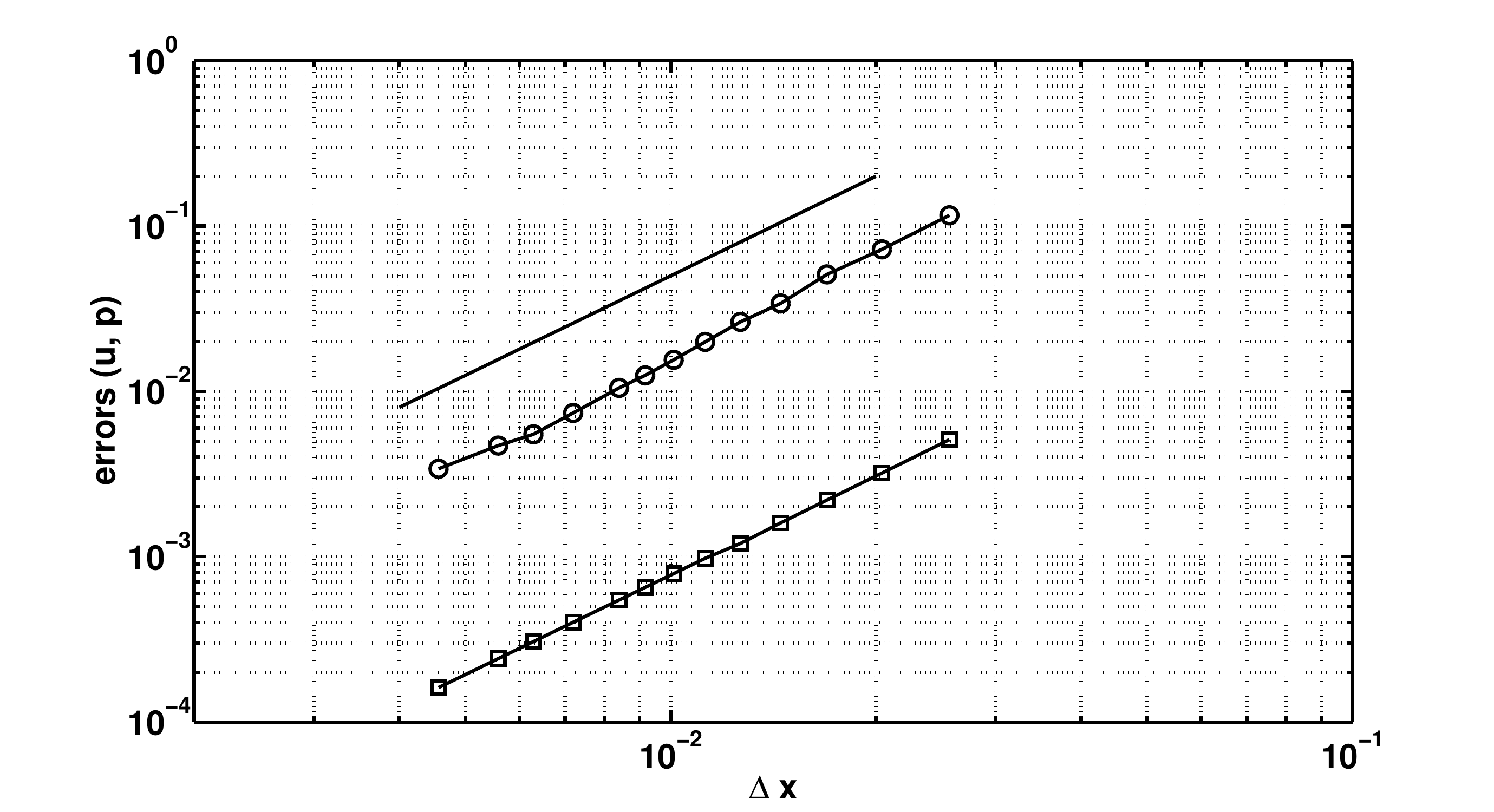}
\caption{(Flow on an irregular domain example). Convergence plot for the
   pressure (circles: $\;\circ\/$) and the horizontal velocity (squares:
   $\;\square\/$). The errors (in the $L^{\infty}\/$ norm) are computed at the
   fixed time $t = 0.0657\/$, for different grid resolutions, with
   $\Delta t = 0.2\,(\Delta x)^2\/$. The slope of the solid straight line
   corresponds to the second order scaling error $\propto (\Delta x)^2\/$.}
\label{fig:upConvergence}
\end{figure}
\begin{figure}[htb!]
\centering
\includegraphics[width=0.85\textwidth]{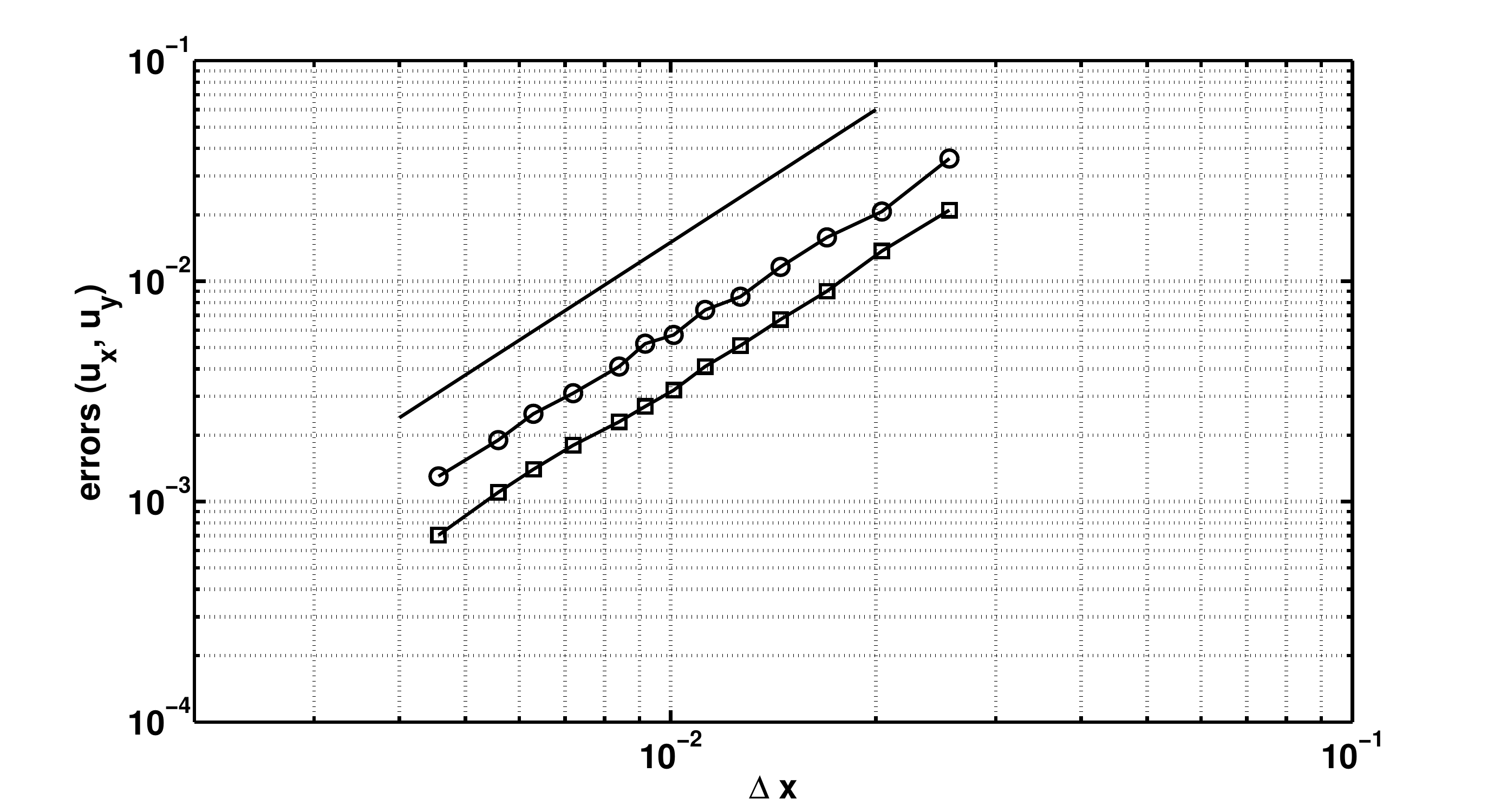}
\caption{(Flow on an irregular domain example). Convergence plots for the
   partial derivatives of the horizontal velocity: $u_y\/$
   (circles: $\;\circ\/$) and $u_x\/$ (squares: $\;\square\/$).
   The errors (in the $L^{\infty}\/$ norm) are computed at the fixed time
   $t = 0.0657\/$, for different grid resolutions, with
   $\Delta t = 0.2\,(\Delta x)^2\/$. The slope of the solid straight line
   corresponds to the second order scaling error $\propto (\Delta x)^2\/$.
   The fact that the errors for the gradient of the velocity are also second
   order is important --- see \S~\ref{subsec:DerivativesOrder}.}
\label{fig:derivativeConvergence}
\end{figure}

Figure~\ref{fig:upConvergence} shows the convergence of velocity and pressure,
while figure~\ref{fig:derivativeConvergence} shows the convergence of the
partial derivatives of the horizontal velocity $u\/$.  The figures indicate
second order $L^{\infty}\/$ convergence for the velocity $u\/$, the pressure
$p\/$, and even the velocity gradient $(u_x\/,\, u_y)\/$ --- see
\S~\ref{subsec:DerivativesOrder}.

\begin{center} \emph{Typical error behavior.} \end{center}

We illustrate the (typical) error behavior --- both in time and in
space, by showing the results of a calculation done at a fixed resolution.
Specifically, we take an $80 \times 80\/$ grid, with
$\Delta t = 0.2 (\Delta x)^2\/$ and $\lambda = 100$, and solve the forced
system of equations in this ``flow on an irregular domain'' example --- recall
that $\mu = 1$.

\begin{figure}[htb!]
\centering
\includegraphics[width=0.75\textwidth]{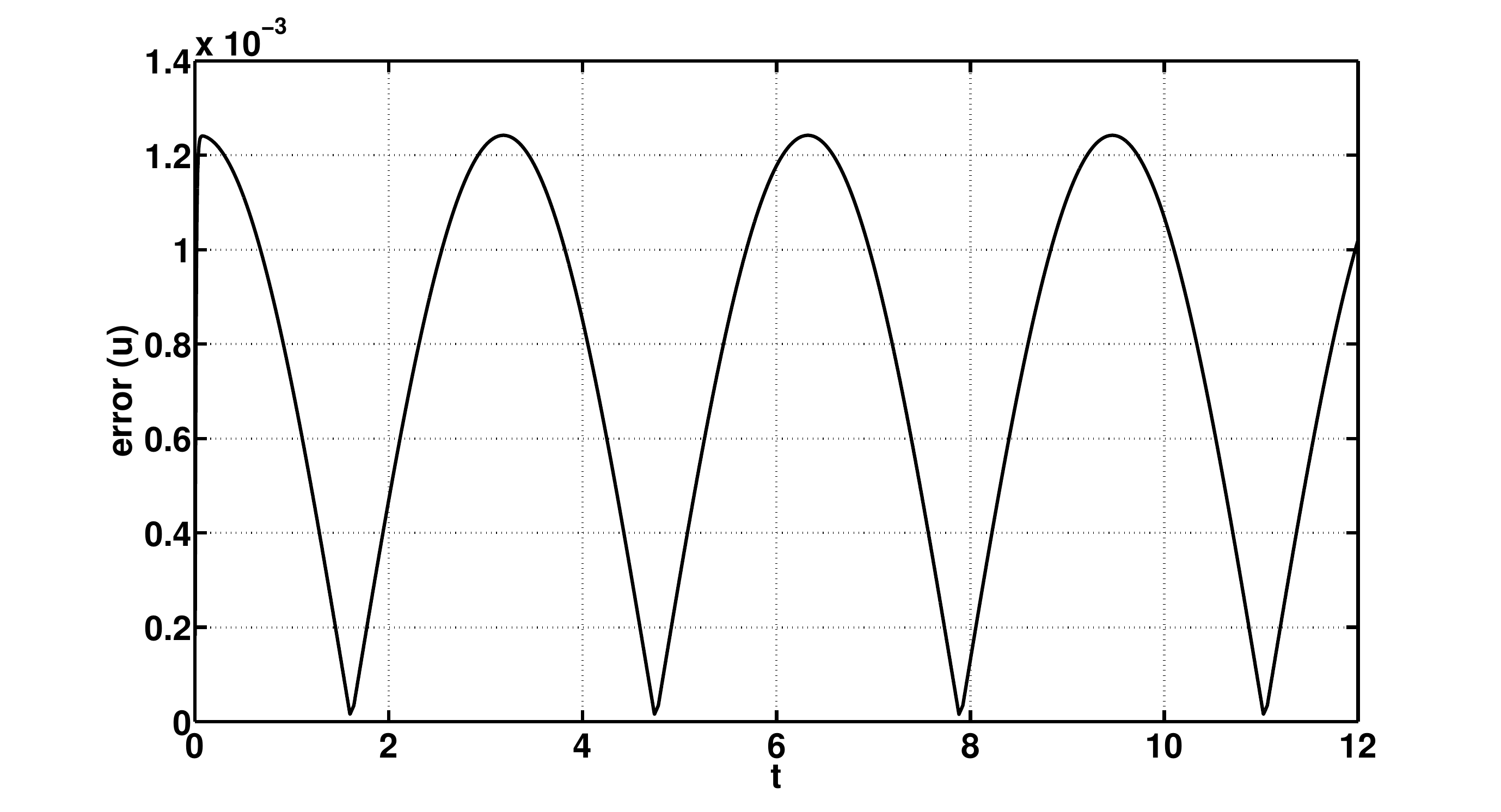}
\caption{(Flow on an irregular domain example, on an $80 \times 80$ grid).
   Evolution in time of the $L^{\infty}\/$ norm of the (spatial) error for
   the horizontal velocity.}
\label{fig:LinfErrorsVsTime}
\end{figure}

Figure~\ref{fig:LinfErrorsVsTime} shows the time evolution of the
$L^{\infty}\/$ (spatial) error in the horizontal velocity.
It should be clear that, while the errors oscillate (over almost one decade in
amplitude), they do not exhibit any measurable growth with time.

\begin{figure}[htb!]
\centering
\includegraphics[width=0.80\textwidth]{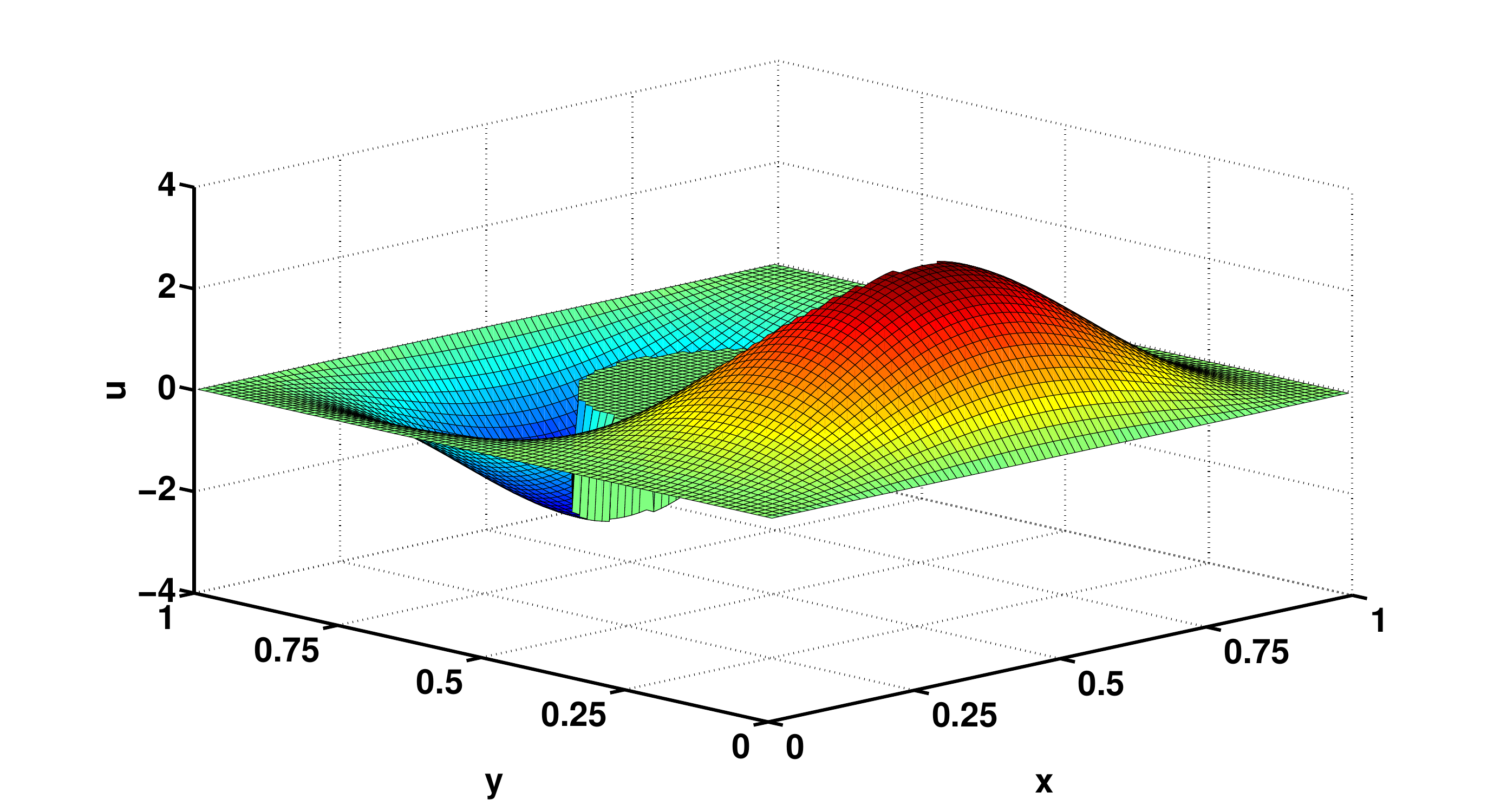}
\caption{(Flow on an irregular domain example, on an $80 \times 80$ grid).
   Numerical horizontal velocity field $u\/$ at $t = 4\,\pi$.}
\label{fig:USurf}
\end{figure}

\begin{figure}[htb!]
\centering
\includegraphics[width=0.80\textwidth]{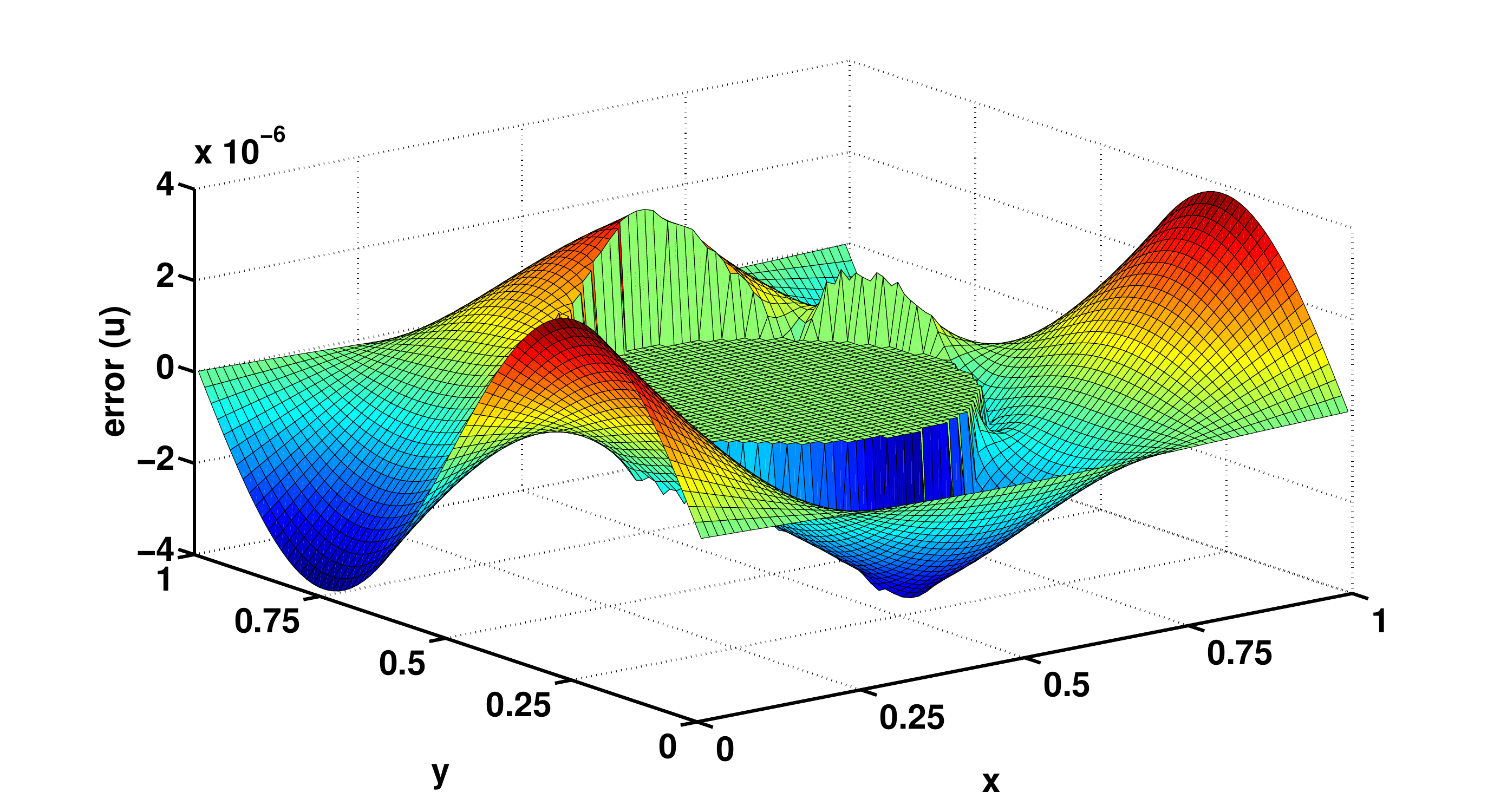}
\caption{(Flow on an irregular domain example, on an $80 \times 80$ grid). Error
   field for the horizontal velocity $u\/$ at $t = 4\,\pi$. }
\label{fig:USurfError}
\end{figure}

Figure~\ref{fig:USurf} shows the horizontal velocity field at time
$t = 4\,\pi$, 
while figure~\ref{fig:USurfError} shows the associated error field. Notice that,
while the error is largest near the boundary --- as expected from the
difficulties that a curved, non-conforming, boundary causes --- it is fairly
well behaved, without abrupt transitions on the scale of the grid size. This
explains why the error in the gradient of the velocity is also second order
accurate --- see figure~\ref{fig:derivativeConvergence} and
remark~\ref{fig:derivativeConvergence}, a feature that should be particularly
useful for the calculation of fluid-solid forces (stresses) along domain
boundaries.

\begin{figure}[htb!]
\centering
\includegraphics[width=0.50\textwidth]{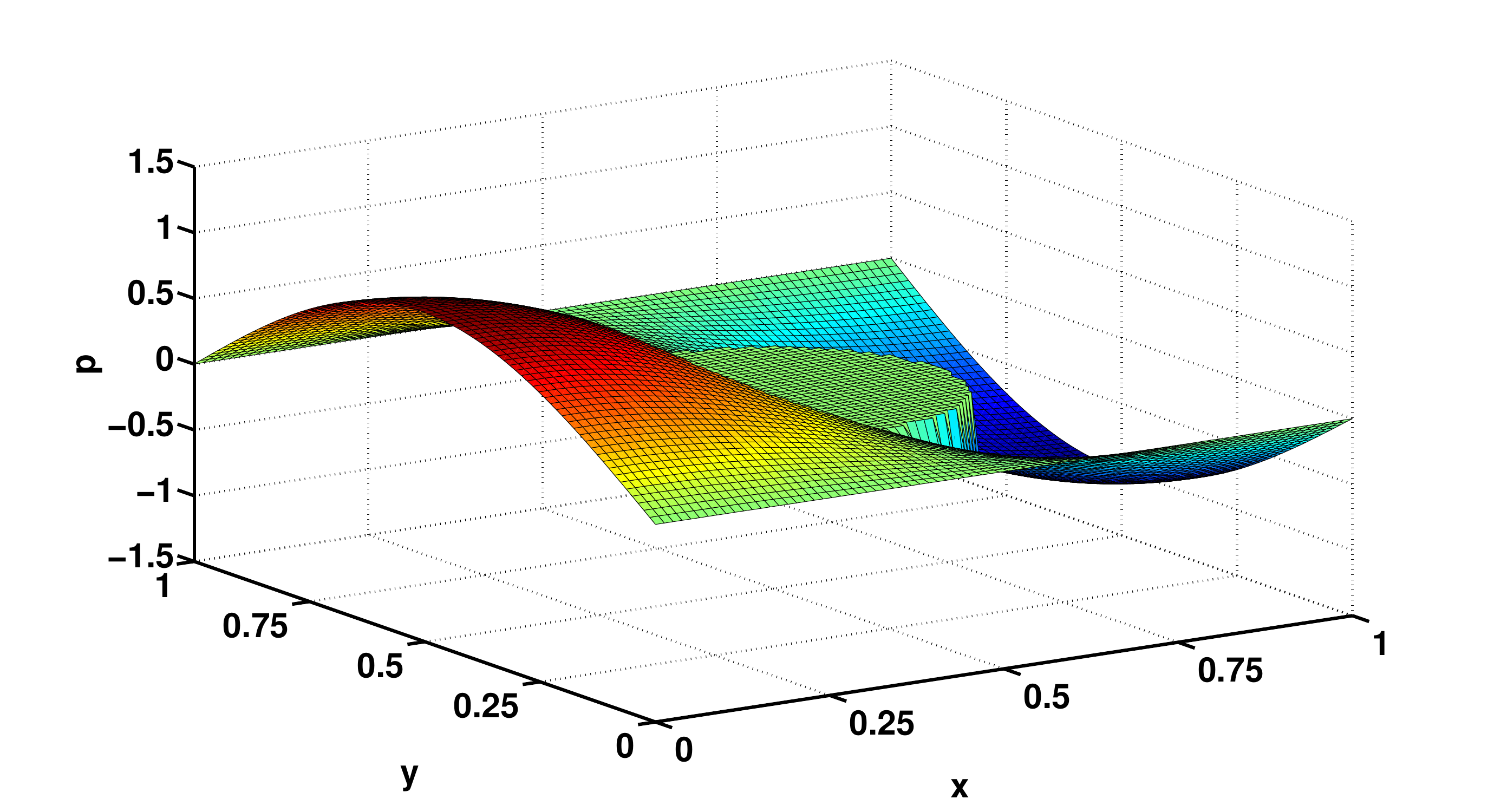}
\hspace*{-3ex}
\includegraphics[width=0.50\textwidth]{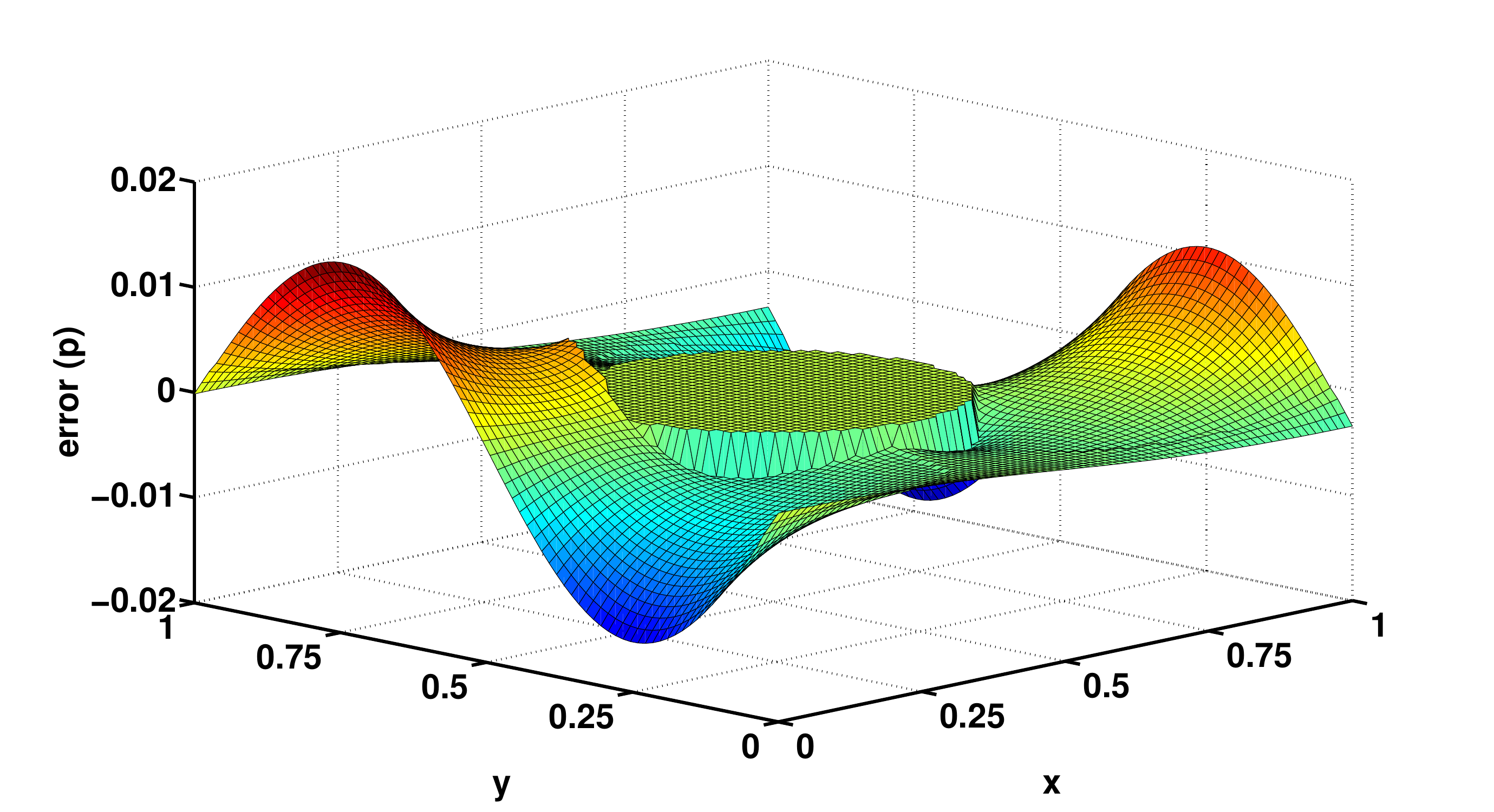}
\caption{(Flow on an irregular domain example, on an $80 \times 80$ grid).
   Left: numerical pressure field $p\/$ at $t = 4\,\pi\/$.
   Right: associated error field.}
\label{fig:pSurf}
\end{figure}

Similarly, figure~\ref{fig:pSurf} shows plots of the pressure and of the
associated pressure error field, again for the time $t = 4\,\pi$. Just as
in the case of the velocity, the errors are dominated at the 
boundary. We did not check the order of
convergence for the gradient of the pressure since one does not
need the gradient of the pressure to compute forces.

\begin{figure}[htb!]
\centering
\includegraphics[width=0.80\textwidth]{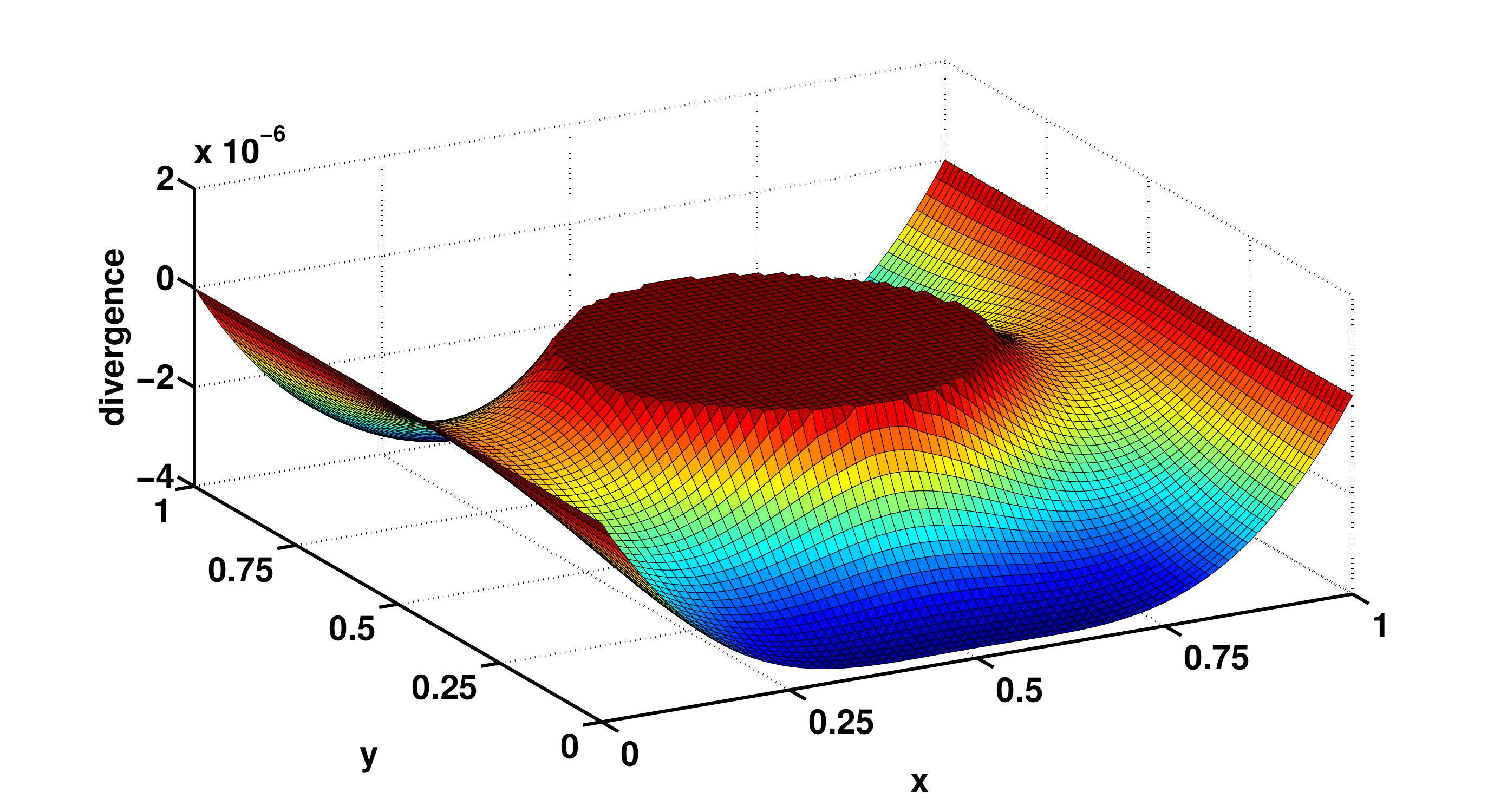}
\caption{(Flow on an irregular domain example, on an $80 \times 80$ grid).
   Numerical error in the divergence of the velocity at $t = 4\,\pi$.
   The $O(10^{-6})$ amplitude is the result of the $O(\Delta x^2)$
   order of the scheme.}
\label{fig:Divergence}
\end{figure}

Finally, figure~\ref{fig:Divergence} shows a plot of the error in the numerical
divergence of the flow field, also for time $t = 4\,\pi$. The errors in the
divergence are also second order. However: notice that they are much smaller
than the errors in the flow field or in the pressure. This is indicative of the
rather strong enforcement of incompressibility that equations (\ref{eqn:heatdivu})
and (\ref{eqn:bcdivu}) guarantee.
%
%
\subsection{Convergence of the derivatives} \label{subsec:DerivativesOrder}
A key reason to worry about uniform convergence, up to the boundary, of the
numerical solution is that this is a necessary condition for the accurate
modeling of solid-fluid interactions.
The evaluation of fluid-solid stresses requires both the pressure, as well as
the derivatives of the velocities, to be accurate at the boundaries. Hence, in
this paper we investigated the behavior of the errors for not just the
pressure $p\/$ and the velocity fields $(u\/,\,v)\/$, but for the gradient
of the flow velocity as well.

An \emph{important point} to notice is that
figure~\ref{fig:derivativeConvergence} shows that
\begin{equation} \label{eqn:2ndOrderStress}
   \left.
   \parbox{0.70\textwidth}{The derivatives of the velocities exhibit
      errors that appear to be second order as well--- in the
      $L^\infty\/$ norm.}\    \right\}
\end{equation}
At first sight this may seem surprising: since the velocities are second order
accurate, one would expect their derivatives to be first order only. However,
this is a worse case scenario, based on the assumption of random errors. On
the other hand, the errors for finite differences approximations (on regular
grids) are typically \emph{not} random: If the solutions are smooth enough,
the errors can be expanded in powers of $\Delta\/x\/$, with coefficients that
involve derivatives of the solution. Hence, in this case, taking low order
derivatives\footnote{At high order, round off errors dominate. There is an
   irreducible $O(\epsilon\,(\Delta\/x)^{-q})\/$ contribution from them,
   where $\epsilon\/$ is the round off error, and $q\/$ is the  order of
   the derivative.}
does not degrade the order of convergence. This, of course, is an important advantage of
finite differences over other approaches (e.g. projection) which do not have
this property.

The astute reader may have noticed that the argument in the prior paragraph
avoids the issue of boundary condition implementation, which can ruin
the smooth error expansions available for finite differences. This requires
some extra considerations:
\begin{itemize}
 \item[\ref{Section:Implementation}a]
 On conforming boundaries, it is usually possible to approximate the boundary
 conditions in such a way that smooth expansions (in powers of $\Delta\,x\/$)
 remain available for the truncation errors.
 \item[\ref{Section:Implementation}b]
 On non-conforming boundaries, on the other hand, it is very easy to ruin the
 smooth error expansions. The local stencils used to approximate the boundary
 conditions along the boundary will, typically, experience abrupt changes in
 response to the placement of the boundary relative to the regular grid.
\end{itemize}
This is where the global approach that we used in \S~5.2 to implement the
boundary conditions comes to the rescue: by linking each local stencil to its
neighbors, the non-smoothness caused by abrupt stencil changes is smeared.
Thus a better behaved error is obtained, which then explains why the result
in (\ref{eqn:2ndOrderStress}) occurs. Notice that, as of now, we neither know
if the smearing process is enough to make the errors $C^{1}\/$ at leading
order --- which is what is needed to make the errors in the velocity gradient
second order, nor if the errors in the velocity gradient are actually second
order. However, the numerical evidence seems to point in this direction.

A final point is that, the ``simplest'' way to get around the issue in
item~\ref{Section:Implementation}b above, is to implement the boundary
conditions at a higher order than required. This has the disadvantage that it
can lead to a messier than needed algorithm, but (generally) it should not
significantly increase the computational cost --- since it only involves the
boundary points of the grid.


\section{Conclusions} \label{Section:conclusions}
Through the introduction of the pressure Poisson equation with consistent
boundary conditions, we give an equivalent formulation of the incompressible
Navier-Stokes equations. In this formulation the momentum equation takes the
form of a vector heat equation with unconventional boundary conditions, while
the pressure Poisson equation can be used to explicitly describe the pressure
as a function of velocity at any fixed time. It follows that the reformulated
system of equations is ideal for using efficient numerical marching methods,
where there are no particular theoretical limitations to the order of accuracy
or the method of implementation.  

In addition, we devise and implement a second order discretization (uniform up
to the boundary) of the equations on irregular domains.  We address the issue of
numerical stability for the normal boundary velocity by adding an appropriate
feedback term to the equations.  For the scheme, we use finite differences on
a regular, staggered grid, so that the domain boundary can be immersed within
an $N\times M\/$ mesh.  We discretize to second order spatial accuracy, and
verify the order of accuracy in $L^{\infty}\/$, both for the velocity and the
pressure on an irregular domain.  Hence, the solutions\footnote{Note that we
   have explored only the case where the solution is at least $C^4\/$.}
converge $O((\Delta x)^2)\/$ uniformly (all the way up to the boundary) as the
grid spacing $\Delta x \rightarrow 0\/$.  Although, we formulate and implement
the scheme in two dimensions, the algorithm has a natural extension to three
dimensions.

There are several issues and extensions that we hope to address in future work.
First, for irregular domains on a regular grid, the proposed momentum equation
boundary conditions implicitly couple all the components of the velocity field.
Therefore, the resulting vector heat operator, $\partial_t - \mu\,\Delta$,
cannot be solved ``component by component''. This makes the implementation of
schemes that do a naive implicit treatment of the stiff viscosity term
computationally expensive, is there a better way? Second, the
formulation in this paper is for fixed domains, but extensions to deformable
and or moving geometries seem possible. Third, can the ideas in this paper be
extended to flows with variable densities (e.g.: stratified flows)? Lastly,
there is the issue of obtaining implementation of the boundary conditions that
yield smooth errors. As explained in \S~\ref{subsec:DerivativesOrder}, this
guarantees that the derivatives of the flow velocity have the same order of
accuracy as the velocity itself, allowing an accurate calculation of the fluid
forces on objects immersed in the flow (as well as the domain boundaries).


\section*{Acknowledgments}
We would like to thank Matt Ueckermann for many helpful discussions and
suggestions. In addition, we greatly appreciate conversations with Alex Marques
and J. C. Nave regarding the numerical implementation of the method, and both
Alex Chorin and Benjamin Seibold for various comments regarding theoretical
considerations. The authors would like to acknowledge the support of the
National Science and Engineering Research Council of Canada, as well as the
National Science Foundation. This research was partially supported by an
NSERC PGS, and by NSF grant DMS-0813648.

\appendix

\section{Further Modification for Solvability} \label{solvabilityBreakdown}
In this appendix we address the question posed and motivated in
remark~\ref{rem:generalIC}. Namely: \emph{Can the equations in
(\ref{PPE_Numerical_Splitting1}--\ref{PPE_Numerical_Splitting2}) be modified
in such a way that they make sense even for initial conditions that are not
incompressible?} In fact, in such a way that if a solution starts with
$\nabla \cdot \mathbf{u} \neq 0\/$ and
$\mathbf{n}\cdot (\mathbf{u} - \mathbf{g}) \neq 0\/$, then
(as $t \to \infty\/$) $\nabla \cdot \mathbf{u} \to 0\/$ and
$\mathbf{n}\cdot (\mathbf{u} - \mathbf{g}) \to 0\/$ --- so that the solution
converges towards a solution of the Navier-Stokes equation.

As pointed out in remark~\ref{rem:generalIC}, the problem with
(\ref{PPE_Numerical_Splitting1}--\ref{PPE_Numerical_Splitting2}) is that the
solvability condition for the Poisson equation (\ref{PPE_Numerical_Splitting2})
is not satisfied when $\nabla \cdot \mathbf{u} \neq 0\/$. Hence the equations
become ill-posed, as they have no solution. The obvious answer to this dilemma
is to interpret the solution to (\ref{PPE_Numerical_Splitting2}) in an
appropriate least squares sense, which is equivalent to modifying the
non-homogeneous terms in the equation by \emph{projecting} them onto the space
of right hand sides for which the Poisson equation has a solution.
Symbolically, write (\ref{PPE_Numerical_Splitting2}) in the form
\begin{equation} \label{PPE_UnProjectedPoisson}
   \left.
   \begin{array}{rcll}
     \Delta\/p                  & = &  g
        \quad & \textrm{for}\;\; \mathbf{x} \in \phantom{\partial} \Omega\/,
        \\ \rule{0ex}{2.5ex}
     \mathbf{n} \cdot \nabla\/p & = &  h
        \quad & \textrm{for}\;\; \mathbf{x} \in \partial \Omega\/,
   \end{array}
   \right\}
\end{equation}
where $g\/$ and $h\/$ are defined in (\ref{PPE_Numerical_Splitting2}). Then
modify the equation to
\begin{equation} \label{PPE_ProjectedPoissonSim}
   \left.
   \begin{array}{rcll}
     \Delta\/p                  & = &  g_p
        \quad & \textrm{for}\;\; \mathbf{x} \in \phantom{\partial} \Omega\/,
        \\ \rule{0ex}{2.5ex}
     \mathbf{n} \cdot \nabla\/p & = &  h_p
        \quad & \textrm{for}\;\; \mathbf{x} \in \partial \Omega\/,
   \end{array}
   \right\}
\end{equation}
where $(g_p\/,\,h_p) = \mathbbm{P}\,(g\/,\,h)\/$ for some projection operator
$\mathbbm{P}\/$ such that
\begin{equation} \label{ProjSolv}
  \int_\Omega g_p\,\du\/V = \int_{\partial\,\Omega} h_p\,\du\/A\/.
\end{equation}
The question, however, is: which projection?

The discretization failures in solvability arising during the course of a
numerical calculation are small, since it should be
$(g_e\/,\,h_e) - (g\/,\,h) = O\left((\Delta\,x)^q\right)\/$ --- where $q\/$
is the order of the method, and $(g_e\/,\,h_e)\/$ is the exact right hand
side. Thus one can argue that, as long as
$(g_p\/,\,h_p) - (g\/,\,h) = O((\Delta\,x)^q)\/$, the resulting numerical
solution will be accurate to within the appropriate order
--- see remarks~\ref{rem:generalIC} and \ref{rem:MatrixSolvability}. In
this section, however, the aim is to consider situations where there is no
small parameter (i.e. $\Delta\,x\/$) guaranteeing that solvability is
``almost'' satisfied. In particular, we want to consider situations where
$\nabla \cdot \mathbf{u} \neq 0\/$, and $|\nabla \cdot \mathbf{u}| \ll 1\/$
does not apply --- leading to errors in solvability which are not small. It
follows that here we must be careful with the choice of the projection.

Obviously, a very desirable property of the selected projection is that it
should preserve the validity of
equations~(\ref{eqn:heatdivu}--\ref{eqn:bcdivu}) --- so that the time
evolution drives $\nabla \cdot \mathbf{u}\/$ to zero. Hence:
\emph{it must be that} $g_p = g\/$, with only $h\/$ affected by the
projection. Further: since the solvability condition involves $h\/$ only via
its mean value over $\partial\,\Omega\/$, the simplest projection that works
is one that appropriately adjusts the mean of $h\/$, and nothing else. Thus
we propose to modify the equations in
(\ref{PPE_Numerical_Splitting1}--\ref{PPE_Numerical_Splitting2}) as follows:
leave (\ref{PPE_Numerical_Splitting1}) as is, as well as the imposed boundary
condition constraint (\ref{NS_Boundary_ConditionsConstr}), but replace
(\ref{PPE_Numerical_Splitting2}) by
\begin{equation} \label{PPE_ProjectionSplitting}
   \left.
   \begin{array}{rcll}
     \Delta\/p                  & = &
        - \nabla \cdot \left( (\mathbf{u} \cdot \nabla)\,\mathbf{u} \right)
        + \nabla \cdot \mathbf{f}
        \quad & \textrm{for}\;\; \mathbf{x} \in \phantom{\partial} \Omega\/,
        \\ \rule{0ex}{2.5ex}
     \mathbf{n} \cdot \nabla\/p & = &
        \mathbf{n} \cdot \left( \mathbf{f} - \mathbf{g}_t
        + \mu\,\Delta\,\mathbf{u}
        - (\mathbf{u} \cdot \nabla)\,\mathbf{u} \right) &
        \\ \rule{0ex}{2.5ex}
                                & + &
        \lambda\,\mathbf{n} \cdot (\mathbf{u} - \mathbf{g}) - \mathcal{C}
        \quad & \textrm{for}\;\; \mathbf{x} \in \partial \Omega\/,
   \end{array}
   \right\}
\end{equation}
where
\begin{equation} \label{eqn:CorrectForDivu}
   \mathcal{C} = \frac{1}{S}\,\int_{\partial\,\Omega} \mathbf{n} \cdot
   \left(\mu\,\Delta\,\mathbf{u} + \lambda\,\mathbf{u} \right)\,\du\,A\/,
\end{equation}
and $S = \int_{\partial\,\Omega} \du\,A\/$ is the surface area of the boundary.
In terms of (\ref{PPE_UnProjectedPoisson}--\ref{PPE_ProjectedPoissonSim})
this corresponds to the projection
\begin{equation}\label{Projection1}
  g_p = g \quad \mbox{and} \quad h_p = h - \mathcal{C}\/,
\end{equation}
where
\begin{equation}\label{Projection2}
  \mathcal{C} = \frac{1}{S}\,\left( \int_{\partial\,\Omega} h\,\du\,A -
  \int_{\Omega} g\,\du\,V \right)\/.
\end{equation}
This is clearly a projection, since $\mathcal{C} = 0\/$ for $(g_p\/,\,h_p)\/$,
so that $\mathbbm{P}^2 = \mathbbm{P}\/$. Further, since the solvability
condition for (\ref{PPE_UnProjectedPoisson}) is precisely $\mathcal{C} = 0\/$
--- see equation (\ref{ProjSolv}), the solvability condition for
(\ref{PPE_ProjectionSplitting}) is satisfied --- even if
$\nabla \cdot \mathbf{u} \neq 0$, though (of course) $\mathcal{C} = 0$
when $\nabla\cdot\mathbf{u} = 0\/$.

Finally, we remark (again) that the projection in  (\ref{Projection1}) is
not unique. In particular, numerical implementations of the Poisson equation
with Neumann conditions often use least squares projections, which alter both
the boundary condition $h\/$ and the source term $g\/$.
This makes sense if the solvability errors are small. However, in general it
seems desirable to not alter the source term, and keep
(\ref{eqn:heatdivu}--\ref{eqn:bcdivu}) valid. This still does not make
(\ref{Projection1}) unique, but it makes it the simplest projection. Others
would also alter (in some appropriate eigenfunction representation) the zero
mean components of $h\/$. This seems not only un-necessarily complicated, but
it may also affect the validity of the result stated below equation
(\ref{Boundary_Evolution_NoDrift2}), defeating the whole purpose of the
reformulation in this appendix.

The system of equations in (\ref{PPE_Numerical_Splitting1}) and
(\ref{PPE_ProjectionSplitting}) makes sense for arbitrary flows
$\mathbf{u}\/$, which are neither restricted by the incompressibility
condition $\phi = \nabla\cdot\mathbf{u} = 0\/$ in $\Omega\/$, nor the normal
velocity boundary condition
$\mathcal{E} = \mathbf{n} \cdot (\mathbf{u} - \mathbf{g}) = 0\/$. Furthermore:
\emph{this system, at least in bounded domains $\Omega\/$, includes the
Navier-Stokes equations as a global attractor for the smooth solutions.} This
is easy to see as follows:

First, because of equations (\ref{eqn:heatdivu}--\ref{eqn:bcdivu}), $\phi\/$
decays exponentially, at a rate controlled by the smallest eigenvalue of
$L = -\Delta\/$ in $\Omega\/$, with Dirichlet boundary conditions. In
particular:
\[
 \mathcal{C} = \frac{1}{S}\,\int_{\partial\,\Omega} \mathbf{n} \cdot
   \left(\mu\,\Delta\,\mathbf{u} + \lambda\,\mathbf{u} \right)\,\du\,A =
   \frac{1}{S}\,\int_{\Omega} \left(\mu\,\Delta\,\phi +
   \lambda\,\phi\right)\,\du\,V
\]
vanishes exponentially.

Second, it is easy to see that, for the system in
(\ref{PPE_Numerical_Splitting1}) and (\ref{PPE_ProjectionSplitting}),
equation~(\ref{Boundary_Evolution_NoDrift}) is modified to
\begin{equation} \label{Boundary_Evolution_NoDrift2}
  \mathcal{E}_t = -\lambda\,\mathcal{E} + \mathcal{C}
  \quad \textrm{for}\;\; \mathbf{x} \in \partial \Omega\/,
\end{equation}
where
$\mathcal{E} = \mathbf{n} \cdot (\mathbf{u} - \mathbf{g})\/$. Hence
$\mathcal{E}\/$ also vanishes exponentially.

Of course, if $\phi = 0\/$ and $\mathcal{E} = 0\/$ initially, then they
remain so for all times, and the evolution provided by
(\ref{PPE_Numerical_Splitting1}) and (\ref{PPE_ProjectionSplitting})
is, exactly, the Navier-Stokes evolution.

In conclusion, the formulation in this section is not only an interesting
theoretical fact. It also provides a robust framework within which numerical
solvers for the incompressible Navier Stokes equations can be developed,
without having to worry about the (potentially deleterious) effects that
discretization (or initial condition) errors, can cause when they violate
mass conservation --- because either $\nabla\cdot\mathbf{u} = 0\/$ in
$\Omega\/$, or $\mathbf{n} \cdot (\mathbf{u} - \mathbf{g}) = 0\/$ in
$\partial\,\Omega\/$, fail. In addition, the formulation eliminates the 
necessity of having to enforce the condition $\nabla\cdot\mathbf{u} = 0\/$ 
directly, which is a core difficulty for the solution of the
Navier-Stokes equations.


\section{The $\nabla \cdot \mathbf{u} = 0\/$ Boundary Condition}\label{divfreeBC}
In the main body of the paper, the $\nabla \cdot \mathbf{u} = 0\/$ boundary
condition is implemented using a regular Cartesian grid. In Cartesian
coordinates, this boundary condition relates the horizontal and vertical
velocities (for two dimensions) via the standard formulas
\begin{eqnarray}
  \mathbf{u} & = &
     u\;\boldsymbol{\hat{\imath}} + v\;\boldsymbol{\hat{\jmath}}\/,
     \label{eqn:CartesianSplit} \\
  \nabla \cdot \mathbf{u} & = &
     \frac{\partial\,u}{\partial\,x} + \frac{\partial\,v}{\partial\,y}\/,
     \label{eqn:CartesianDiv}
\end{eqnarray}
where $\boldsymbol{\hat{\imath}}\/$ and $\hat{\boldsymbol{\jmath}}\/$ are the
coordinate unit vectors.

For some applications, it may be more convenient to implement the numerical
computation on a coordinate system where the boundaries are conforming. In
general curvilinear coordinate systems, the divergence acquires a more
complicated form than the one above. Hence, it is our purpose here
to display the form that the  $\nabla \cdot \mathbf{u} = 0\/$ boundary
condition takes for a conforming boundary in a general orthogonal curvilinear
coordinate system. We do the 2-D case only --- the 3-D case is quite similar.
%
%
\subsection{Curvilinear coordinates and a conforming boundary}
In some region near the boundary $\partial\,\Omega\/$ of the domain of
integration, assume that an orthogonal, curvilinear, set of coordinates
$(\eta\/,\,\xi)\/$ has been selected --- such that the boundary is given by
$\eta = c\/$, for some constant $c\/$.  In terms of the coordinates 
$(\eta\/,\,\xi)\/$, the vector field $\mathbf{u}\/$ can then be written in the form
\begin{equation} \label{eqn:GeneralSplit} \\
    \mathbf{u} = u_{\eta}\;\boldsymbol{\hat{\eta}} +
                 u_{\xi} \;\boldsymbol{\hat{\xi}}\/,
\end{equation}
where $u_{\eta}\/$ and $u_{\xi}\/$ are the velocity components in the coordinate
directions given by the unit vectors $\boldsymbol{\hat{\eta}}\/$
and  $\boldsymbol{\hat{\xi}}\/$, respectively. Let now the functions
$x = x(\eta\/,\,\xi)\/$ and $y = y(\eta\/,\,\xi)\/$ describe the relationship
of the curvilinear coordinates with a Cartesian system.
Then (e.g. see \cite{Jackson})
\begin{equation} \label{eqn:CurvilinearUV}
  \boldsymbol{\hat{\eta}} = \frac{1}{s_{\eta}}\,\left(
    \frac{\partial\,x}{\partial\,\eta}\;\boldsymbol{\hat{\imath}} +
    \frac{\partial\,y}{\partial\,\eta}\;\boldsymbol{\hat{\jmath}} \right)
  \quad \mbox{and} \quad
  \boldsymbol{\hat{\xi}} = \frac{1}{s_{\xi}}\,\left(
    \frac{\partial\,x}{\partial\,\xi}\;\boldsymbol{\hat{\imath}} +
    \frac{\partial\,y}{\partial\,\xi}\;\boldsymbol{\hat{\jmath}} \right)\/,
\end{equation}
where $s_{\eta}\/$ and $s_{\xi}\/$ are local normalization factors
\begin{equation} \label{eqn:CurvilinearNF}
  s_{\eta}^2 = \left(\frac{\partial\,x}{\partial\,\eta}\right)^2 +
              \left(\frac{\partial\,y}{\partial\,\eta}\right)^2
  \quad \mbox{and} \quad
  s_{\xi}^2  = \left(\frac{\partial\,x}{\partial\,\xi}\right)^2 +
              \left(\frac{\partial\,y}{\partial\,\xi}\right)^2\/.
\end{equation}
Further, re-interpreting $s_{\eta}\/$ and  $s_{\xi}\/$ as functions of 
$(\eta\/,\,\xi)\/$, we have
\begin{equation} \label{eqn:GeneralDiv}
  \nabla \cdot \mathbf{u} = \frac{1}{s_{\eta}\,s_{\xi}}\,\left(
     \frac{\partial}{\partial\,\eta} \left(s_{\xi} \,u_{\eta}\right) +
     \frac{\partial}{\partial\,\xi}  \left(s_{\eta}\,u_{\xi} \right)
     \right)\/.
\end{equation}
The $\nabla \cdot \mathbf{u} = 0\/$ boundary condition is accompanied by the
Dirichlet boundary condition
$\mathbf{n} \times \mathbf{u} = \mathbf{n} \times \mathbf{g}\/$ for the vector
momentum equation. For a conforming boundary, this is just
$u_{\xi} = g_{\xi}\/$ along the boundary curve where $\eta = c\/$. Thus the
$\nabla \cdot \mathbf{u} = 0\/$ boundary condition reduces to a \emph{Robin}
boundary condition on the normal velocity (provided that the boundary is
piecewise smooth)
\begin{equation} \label{eqn:RobinBC}
  \frac{\partial}{\partial\,\eta}\left(s_{\xi}\,u_{\eta}\right) +
  \frac{\partial}{\partial\,\xi} \left(s_{\eta}\,g_{\xi} \right) = 0
\end{equation}
along $\eta = c\/$.
%
%
\subsection{Example: polar coordinates}
Consider the special case of a polar coordinate system $x = r\/\cos\/\theta\/$
and $y = r\/\sin\/\theta\/$, with the boundary at $r = 1\/$. Let the angular
velocity at the boundary be given by
$u_{\theta}(1\/,\,\theta) = g_{\theta}(\theta)\/$. Then since (in this case)
$s_{\theta} = r\/$ and $s_{r} = 1\/$, the zero divergence boundary condition
becomes
\begin{equation}
  \frac{\partial}{\partial\,r} \left(r\,u_{r}\right) +
  \frac{\partial}{\partial\,\theta}\,g_{\theta} = 0
  \quad \mbox{for} \;\; r = 1\/.
\end{equation}


\bibliographystyle{elsarticle-num}

\end{document}